\documentclass{article}

\usepackage{amssymb,latexsym,amsmath, amscd}

\usepackage[pdftex]{graphicx}

\newtheorem{theorem}{Theorem}[section]
\newtheorem{corollary}[theorem]{Corollary}
\newtheorem{lemma}[theorem]{Lemma}
\newtheorem{proposition}[theorem]{Proposition}
\newtheorem{definition}[theorem]{Definition}

\newtheorem{remark}[theorem]{Remark}

\begin{document}

\newcommand{\lrightarrow}{\stackrel{L}{\rightarrow}}

\newcommand{\rrightarrow}{\stackrel{R}{\rightarrow}}

\renewcommand{\contentsname}{ }

\title{Sub-riemannian geometry from intrinsic viewpoint}

\author{Marius Buliga \\
\\
Institute of Mathematics, Romanian Academy \\
P.O. BOX 1-764, RO 014700\\
Bucure\c sti, Romania\\
{\footnotesize Marius.Buliga@imar.ro}}

\date{This version: 14.06.2012}

\maketitle

\begin{abstract}
Gromov proposed to extract the (differential) geometric content of a sub-riemannian space 
exclusively from its Carnot-Carath\'eodory distance. One of the 
most striking features of a regular sub-riemannian space is that it has at any
point a metric tangent space with the algebraic structure of  a Carnot group, hence a homogeneous Lie group.   Siebert characterizes homogeneous Lie
groups as locally compact groups admitting a contracting and continuous
one-parameter group of automorphisms. Siebert result has not a metric character.

In these notes I show that sub-riemannian geometry may be described by 
about 12 axioms, without using any a priori given differential structure, but using dilation structures instead.   Dilation structures bring forth the other intrinsic ingredient,
namely the dilations, thus blending Gromov metric point of view with Siebert algebraic 
one.
\end{abstract}

{\bf MSC2000:} 51K10, 53C17, 53C23

%\newpage 

\tableofcontents

%\newpage

\section{Introduction}

In these notes I show that sub-riemannian geometry may be described intrinsically, in terms of dilation structures, {\bf without using any a priori given differential structure}.

A complete riemannian manifold is a length metric space  by the  Hopf-Rinow theorem. The problem of intrinsic characterization of riemannian spaces  asks for the recovery of the manifold structure and of the riemannian metric from the distance function (associated to the length functional).

For 2-dim riemannian manifolds the problem has been solved by A. Wald in 1935 \cite{wald}. 
In 1948 A.D. Alexandrov \cite{alexandrov} introduces its famous curvature (which uses comparison triangles)  and proves that, under mild smoothness conditions  on this curvature, one is capable to recover the differential structure and the metric of the 2-dim riemannian manifold. In 1982 Alexandrov proposes as a conjecture that a characterization of a riemannian manifold (of any dimension) is possible in terms of metric (sectional) curvatures (of the type introduced by Alexandrov) and weak smoothness assumptions formulated in metric way (as for example H\"{o}lder smoothness). 

The problem has been solved by Nikolaev \cite{nikolaev} in 1998. He proves that every locally compact length metric space $M$, not linear at one of its points,  with $\alpha$-H\"{o}lder continuous metric sectional curvature of the “generalized tangent bundle” $T^{m}(M)$ (for some $m=1,2,…$, which admits local geodesic extendability, is isometric to a $C^{m+2}$ smooth riemannian manifold. We shall remain vague about what is the meaning of: not linear, metric sectional curvature, generalized tangent bundle. Please read the excellent paper by Nikolaev for grasping the precise meaning of the result. Nevertheless, we may summarize the solution of Nikolaev like this: 
\begin{enumerate}
\item[-] He constructs a (family of) intrinsically defined tangent bundle(s) of the metric space, by using a generalization of the cosine formula for estimating a kind of a  distance between two curves emanating from different points. This will lead him to a generalization of the tangent bundle of a riemannian manifold endowed with the canonical Sasaki metric. 
\item[-] He defines a notion of sectional curvature at a point of the metric space, as a limit of a function of nondegenerated geodesic triangles, limit taken as these triangles converge (in a precised sense) to the point. 
\item[-] The sectional curvature function thus constructed is supposed to satisfy a smoothness condition formulated in metric terms. 
\end{enumerate}
He proves that under the smoothness hypothesis he gets the conclusion.

In this paper we prove (see theorem \ref{tchariem}): 

\begin{theorem}
The dilation structure associated to  a riemannian manifold, as in proposition \ref{priemann}, is tempered (definition \ref{dtempered}), has the Radon-Nikodym property (definition \ref{defrn}) and is a length dilation structure (definition \ref{deflds}).

 If  $(X, d, \delta)$ is a strong dilation structure which is 
tempered,  it has  the Radon-Nikodym property and moreover for any $x \in X$ the tangent space in the sense of dilation structures (definition \ref{deftdil}) is a commutative local group, then any open, with compact closure subset of $X$ can be endowed with a $\displaystyle C^{1}$ riemannian structure which gives a distance $d'$ which is bilipschitz equivalent with $d$. 
\end{theorem}

Sub-riemannian spaces are length metric spaces as well, why are them different? First of all, any riemannian space is a sub-riemannian one, therefore sub-riemannian spaces are more general than riemannian ones. It is not clear at first sight why the characterization of riemannian spaces does not extend to sub-riemannian ones. In fact, there are two problematic steps for such a program for extending Nikolaev result to sub-riemannian spaces: the cosine formula, as well as the Sasaki metric on the tangent bundle don't have a correspondent in sub-riemannian geometry (there is, basically, no statement canonically corresponding to Pythagoras theorem); the sectional curvature at a point cannot be introduced by means of comparison triangles, because sub-riemanian spaces do not behave well with respect to this comparison of triangle idea.

The problem of intrinsic characterization of sub-riemannian spaces has been formulated by Gromov in \cite{gromovsr}.  Gromov takes  the Carnot-Carath\'eodory distance as  the only intrinsic object of a sub-riemannian space. Indeed, in \cite{gromovsr}. section 0.2.B. he writes: 

"{\em If w live inside a Carnot-Carath\'eodory metric space $V$ we may know nothing whatsoever about the (external) infinitesimal structures (i.e. the smooth structure on $V$, the subbundle $H \subset T(V)$ and the metric $g$ on $H$) which were involved in the construction of the CC metric.}"

He then formulates two main problems: 
\begin{enumerate}
\item[(1)] "{\em Develop a sufficiently rich and robust internal C-C language which would enable us to capture the essential external characteristics of our C-C spaces}". (he proposes as an example to recognize the rank of the horizontal distribution, but in my opinion this is, say, something much less essential than to "recognize" the "differential structure", in the sense proposed here as the equivalence class under local equivalence of dilation structures, see definition \ref{dilequi})
\item[(2)] "\em Develop external (analytic) techniques for evaluation of internal invariants of $V$."
\end{enumerate}

Especially the problem (2) raises a big question: what is, in fact, a sub-riemannian space? Should it be defined only in relation with the differential geometric construction using horizontal distributions and CC distance? In this paper we propose that a sub-riemannian space is a particular case of  a pair of spaces, one looking down on another. The "sub-riemannian"-ess is a relative notion. 

As for the problem (1), a solution is proposed here, by using dilation structure. The starting point is to remark that a  regular sub-riemannian space  has metric tangent spaces with the structure of a  Carnot group. All known proofs are using indeed the intrinsic CC distance, but also the differential structure of the manifold and the differential geometric definition of the CC distance. The latter are not intrinsic, according to Gromov criterion,  even if the conclusion (the tangent space has an algebraic Carnot group structure) looks intrinsic. But is this Carnot group structure intrinsic or is an artifact of the method of proof which was used?  

Independently, there is an interest into the characterization of contractible topological groups. A result of Siebert \cite{siebert} characterizes homogeneous Lie
groups as locally compact groups admitting a contracting and continuous
one-parameter group of automorphisms. This result is relevant because, we argued before, the Carnot group structure comes from the self-similar metric structure of a tangent space, via the result of  Siebert. 

If we enlarge the meaning of "intrinsic", such as to contain the CC distance, but also the approximate self-similar structure of the sub-riemannian space, then we are able to give a characterization of these spaces. This approximate self-similarity is modeled by dilation structures (initially called "dilatation structures" \cite{buligadil1}, a nod to the latin origin of the word). Dilation structures  bring forth the other intrinsic ingredient,
namely the dilations, which are generalizations of Siebert' contracting 
group of automorphisms. 

According to the characterization given in this paper, a regular sub-riemannian space is one which can be constructed from: a tempered dilation structure with the Radon-Nikodym property and commutative tangent spaces,   and from a coherent projection   which satisfies (Cgen), (A) and (B) properties. As a corollary, we recover the known fact that sub-riemannian spaces have the Radon-Nikodym property and we learn the new fact that they are length dilation structures, which provides a characterization of the behaviour of the rescaled length functionals induced by the Carnot-Carath\'eodory distance.

As it is maybe to be expected from a course notes paper, these notes are based on previous papers of mine, mainly \cite{buligachar} (section 12 follows almost verbatim  the section 10 of \cite{buligachar}), \cite{buligadil1}, \cite{buligadil2}, 
\cite{buligasr}, \cite{buligaultra} and also from a number of arxiv papers of mine, mentioned in the bibliography. Many clarifications and theorems are added, in order to construct over the foundations laid elsewhere.   I hope that the unitary presentation will help the understanding of the subject. 

\paragraph{Acknowledgements.}These are  the notes prepared for the course "Metric spaces with dilations and sub-riemannian geometry from intrinsic point of view",  CIMPA research school on sub-riemannian geometry (2012). Unfortunately I have not been able to attend the school. I want 
to express my thanks to the organizers for inviting me and also my excuses for not being there. 
This work was partially  supported by a grant of the Romanian National Authority for Scientific Research, CNCS – UEFISCDI, project number PN-II-ID-PCE-2011-3-0383.

\section{Metric spaces, groupoids, norms}

Metric spaces  have been introduced by Fr\' echet (1906) in the paper \cite{frechet}. 

\begin{definition}
A metric space $(X,d)$ is a pair formed by a set $X$ and a function called  distance,   
$d: X \times X \rightarrow [0,+\infty)$, which satisfies the following: (i) $d(x,y) = 0$ if and only if $x=y$; (ii) (symmetry) for any $x, y \in X$ we have  $d(x,y) = d(y,x)$;  (iii) (triangle inequality) for any $x, y, z \in X$ we have  $d(x, z) \leq d(x,y) + d(y,z)$. 
The ball of radius $r>0$ and center $x \in X$ is the set $\displaystyle  B(x,r) \, = \, \left\{ y \in X \mbox{  :  }  d(x,y) < r \right\}$. 
Sometimes we may  use the notation $\displaystyle B_{d}(x,r)$,  which indicates  the dependence on the distance $d$. 

The topology and uniformity on the   metric space $(X,d)$ is the one generated by balls, respectively by preimages of the distance fonction. 
\label{dmetricspace}
\end{definition}

\subsection{Normed groups and normed groupoids}

Starting from the observation that normed trivial groupoids are in bijective correspondence with metric spaces, proposition \ref{enoughgen}, we think it is  interesting to extend the theory of metric spaces to normed groupoids. This is explained in detail  in \cite{groupoids}. 

Groups ae groupoids with one object, so let us start with normed groups, then pass to normed groupoids.

\begin{definition}
A normed group  is a pair $(G, \rho)$, formed by a group $G$, with operation $(x, y) \in G \times G \mapsto xy$, the inverse  $\displaystyle x \in G \mapsto x^{-1}$  and neutral element denoted by $e$, and a norm  $\rho : G \rightarrow [0, +\infty)$ which satisfies the following:  
(i)  $\rho(x) = 0$ if and only if $x = e$;  (ii) (symmetry) for any $x \in G$ we have $\displaystyle \rho(x^{-1} ) = \rho(x)$; (iii) (sub-additivity) for any $x, y \in G$ we have  
$\rho(xy) \leq \rho(x) + \rho(y)$. 
\label{dnormedgroup}
\end{definition}

A normed group $(G, \rho)$ can be seen as a metric space. Indeed, as expected, the norm $\rho$ induce distances, left or right invariant:  
$$d_{L}(x, y) \, = \, \rho(x^{-1} y)  \quad , \quad d_{R}(x, y) \, = \, \rho(x y^{-1}) \quad .$$

Groupoids are generalization of groups. We shal model them by looking to  the set of arrows, which is a set with a partially defined binary operation and a unary operation (the inverse function). A groupoid norm will be   a function defined on the set of arrows, with properties similar with the ones of a norm over a group.

\begin{definition} 
A normed groupoid $(G, \rho)$ is a pair formed by: 
\begin{enumerate}
\item[-]  a groupoid  $G$, i.e. a set with two partially defined operations:  the composition 
$\displaystyle m: G^{(2)} \subset G \times G \rightarrow G$, denoted multiplicatively $\displaystyle m(a,b) = ab$, and the inverse  $\displaystyle inv: G \rightarrow 
G$, denoted $\displaystyle inv(a) = a^{-1}$. The operations satisfy: 
 for any $a,b,c \in G$ 
\begin{enumerate}
\item[(i)] if $\displaystyle (a,b) \in G^{(2)}$ and 
$\displaystyle (b,c) \in G^{(2)}$ then $\displaystyle (a,bc) \in G^{(2)}$ and 
$\displaystyle (ab, c) \in G^{(2)}$ and we have $a(bc) = (ab)c$, 
\item[(ii)] $\displaystyle (a,a^{-1}) \in G^{(2)}$ and 
$\displaystyle (a^{-1},a) \in G^{(2)}$, 
\item[(iii)] if $\displaystyle (a,b) \in G^{(2)}$ then $\displaystyle a b b^{-1} = a$ and 
$\displaystyle a^{-1} a b = b$. 
\end{enumerate}
The set of objects of the groupoid $X= Ob(G)$ is formed by all products 
$\displaystyle a^{-1} a$, $a \in G$. For any $a \in G$ we let $\alpha(a) = 
a^{-1} a \in X$ to be the source (object) of $a$ and and $\omega(a) = a a^{-1} \in X$ to be the target of $a$.
\item[-]   a (groupoid) norm $d : G \rightarrow [0, +\infty)$ which satisfies: 
\begin{enumerate}
\item[(i)] $d(g) = 0$ if and only if $g \in Ob(G)$, 
\item[(ii)] (symmetry) for any $g \in G$,  $d(g^{-1}) \, = \, d(g)$, 
\item[(iii)] (sub-additivity) for any $\displaystyle (g,h) \in G^{(2)}$,   
$d(gh) \, \leq \, d(g) + d(h)$,  
\end{enumerate}
\end{enumerate} 
\label{dnormedgroupoid}
\end{definition}

If $Ob(G)$ is a singleton then $G$ is just a group and the previous definition corresponds exactly to the definition \ref{dnormedgroup} of a normed group. As in the case of normed groups,  normed groupoids induce metric spaces too. 

\begin{proposition}
Let $(G,d)$ be  a normed groupoid. For any   $x \in Ob(G)$ the pair $\displaystyle (\alpha^{-1}(x), d_{x})$ is a metric space, where for any $g, h \in G$ with $\alpha(g) = \alpha(h) = x$ we define $\displaystyle d_{x}(g,h) \, = \, d( g h^{-1})$. 

Moreover, any  normed groupoid is a  disjoint union of metric spaces 
\begin{equation}
G \, = \, \bigcup_{x \in Ob(G)} \alpha^{-1}(x) \quad , 
\label{disu}
\end{equation}
such that  for any $u \in G$ the "right translation"  
$$\displaystyle R_{u} : \alpha^{-1}\left(\omega(u)\right) \rightarrow \alpha^{-1}\left(\alpha(u)\right) \quad , \quad R_{u} (g) \, = \, gu$$
is an isometry, that is for any  $\displaystyle g, h \in \alpha^{-1}\left(\omega(u)\right)$  
$$d_{\omega(u)} (g,h) \, = \, d_{\alpha(u)} ( R_{u}(g) , R_{u}(h)) \quad .$$
\label{pgroupoid}
\end{proposition}

\paragraph{Proof.}
If $\alpha(g) = \alpha(h)= x$ then $\displaystyle (g, h^{-1}) \in G^{(2)}$, therefore  
$\displaystyle d_{x}(g,h)$ is well defined. The proof of the first part is straightforward, the properties (i), (ii), (iii) of the groupoid norm $\rho$ transform respectively into (i), (ii), (iii) properties of the distance  $\displaystyle d_{x}$. 

For the second part of the proposition remark that  $\displaystyle R_{u}$ is well defined and moreover 
$\displaystyle R_{u}(g) \left( R_{u}(h)\right)^{-1} \, = \, g h^{-1}$. Then
$$d_{\alpha(u)} ( R_{u}(g) , R_{u}(h))  \, = \, d\left( R_{u}(g) \left( R_{u}(h)\right)^{-1}\right) \, = \, $$ 
$$\, = \, d(g h^{-1}) \, = \, d_{\omega(u)} (g,h)$$ 
and the proof is done. \quad $\square$

Conversely, any metric spaces is identified with a normed groupoid. 

\begin{proposition}
Let $(X,d)$ be a metric space and consider the  "trivial groupoid" $G = X \times X$,  with 
multiplication $\displaystyle (x,y) (y,z) \, = \, (x,z)$
and inverse $\displaystyle (x,y)^{-1} = (y,x)$. Then $(G,d)$ is a normed groupoid and moreover 
any component of the decomposition (\ref{disu}) of $G$  is isometric with $(X,d)$. 

Moreover if $G = X \times X$ is the trivial groupoid associated to the set $X$ and $d$ is a norm on $G$ then $(X,d)$ is a metric space.
\label{enoughgen}
\end{proposition}

\paragraph{Proof.}
We have $\alpha(x,y) = (y, y)$ and $\omega(x,y) = (x,x)$, therefore  the set of objects of the trivial groupoid is  
$Ob(G) \, = \, \left\{ (x,x) \mbox{ : } x \in X\right\}$. This set  can be identified with $X$ by the bijection $(x,x) \mapsto x$.  Moreover, for any $x  \in X$ we have $\displaystyle \alpha^{-1}((x,x)) \, = \, X \times \left\{ x\right\}$. 

The distance  $d: X\times X \rightarrow [0,+\infty)$ is a groupoid norm, seen as   
$d: G \rightarrow [0, +\infty)$. Indeed (i) ($d(x,y) = 0$ if and only if $(x,y) \in Ob(G)$) is equivalent with $d(x,y) = 0$ if and only if $x = y$. The symmetry condition (ii) is just the symmetry of the distance $d(x,y) = d(y,x)$. Finally the sub-additivity of $d$ as a groupoid norm is equivalent with the triangle inequality.  In conclusion $(X \times X,d)$ is a normed groupoid if and only if $(X,d)$ is a metric space. 

For any $x  \in X$, let us compute  the distance $\displaystyle d_{(x,x)}$, which is the distance on the space 
$\displaystyle \alpha^{-1}((x,x)) $. We have  
$$d_{(x,x)}( (u,x) , (v, x) ) \, = \, d((u,x) (v,x)^{-1}) \, = \, d((u,x)(x,v)) \, = \, d(u,v)$$
therefore the metric space $\displaystyle (\alpha^{-1}((x,x)), d_{(x,x)}) $ is isometric with 
$(X,d)$ by the isometry $(u,x) \mapsto u$, for any $u \in X$. \quad $\square$

In conclusion normed groups make good examples of metric spaces.

\subsection{Gromov-Hausdorff distance}

We shall denote by $CMS$ the set of isometry classes of compact metric spaces. This set is endowed with the Gromov distance and with the topology is induced by this distance. 

The  Gromov-Haudorff distance shall be introduced by  way of (cartographic like) maps. 
Although this definition is well known, the cartographic analogy  was explained first time in detail in \cite{maps}, which we follow here.

\begin{definition}
Let $\rho \subset X \times Y$ be a relation which represents a map of 
$dom \ \rho \, \subset (X,d)$ into $ im \ \rho \, \subset (Y,D)$. To this map 
are associated three quantities: accuracy, precision and resolution. 

 The accuracy of the map $\rho$ is defined by: 
\begin{equation}
acc(\rho) \, = \, \sup \left\{ \mid D(y_{1}, y_{2}) - d(x_{1},x_{2}) \mid \mbox{
: } (x_{1},y_{1}) \in \rho \, , \, (x_{2},y_{2}) \in \rho \right\}
\label{acc1}
\end{equation}
 The resolution of $\rho$ at  $y \in im \ \rho$ is 
\begin{equation}
res(\rho)(y) \, = \, \sup \left\{ d(x_{1},x_{2}) \mbox{
: } (x_{1},y) \in \rho \, , \, (x_{2},y) \in \rho \right\}
\label{resy1}
\end{equation}
and the resolution of $\rho$ is given by: 
\begin{equation}
res(\rho) \, = \, \sup \left\{ res(\rho)(y) \mbox{
: } y \in \, im \ \rho \right\}
\label{res1}
\end{equation}
 The precision of $\rho$ at $x \in dom \ \rho$ is 
\begin{equation}
prec(\rho)(x) \, = \, \sup \left\{ D(y_{1},y_{2}) \mbox{
: } (x,y_{1}) \in \rho \, , \, (x,y_{2}) \in \rho \right\}
\label{precx1}
\end{equation}
and the precision of $\rho$ is given by: 
\begin{equation}
prec(\rho) \, = \, \sup \left\{ prec(\rho)(x) \mbox{
: } x \in \, dom \ \rho \right\}
\label{prec1}
\end{equation}
\label{defacc}
\end{definition}

We may need to perform also a "cartographic generalization", starting from   a relation $\rho$, with domain $M = \, dom(\rho)$  which is $\varepsilon$-dense in $(X,d)$. 

\begin{definition}
A subset $M \subset X$ of a metric space $(X,d)$ is $\varepsilon$-dense in 
$X$ if for any $u \in X$ there is $x \in M$ such that $d(x,u) \leq \varepsilon$. 

Let $\rho \subset X \times Y$ be a relation such that $dom \ \rho$ is 
$\varepsilon$-dense in $(X,d)$ and $im \ \rho$  is $\mu$-dense in 
 $(Y,D)$. We define then $\displaystyle \bar{\rho} \subset X \times Y$ by:
 $\displaystyle (x,y) \in \bar{\rho}$ if there is $\displaystyle (x',y') \in
 \rho$ such that $\displaystyle d(x,x')\leq \varepsilon$ and $\displaystyle 
 D(y,y') \leq \mu$.
 \label{defgencart1}
\end{definition}

If $\rho$ is a relation as described in definition \ref{defgencart1} then 
accuracy $acc(\rho)$, $\varepsilon$ and $\mu$ control the precision 
$prec(\rho)$ and resolution $res(\rho)$. Moreover, the accuracy, precision and 
resolution of $\displaystyle \bar{\rho}$ are controlled by those of 
$\rho$ and $\varepsilon$, $\mu$, as well. 

\begin{proposition}
Let $\rho$ and $\displaystyle \bar{\rho}$ be as described in definition
\ref{defgencart1}. Then: 
\begin{enumerate}
\item[(a)] $\displaystyle res(\rho) \, \leq \, acc(\rho)$, 
\item[(b)] $\displaystyle prec(\rho) \, \leq \, acc(\rho)$, 
\item[(c)] $\displaystyle res(\rho) + 2 \varepsilon \leq \, res(\bar{\rho}) \leq
\, acc(\rho) + 2(\varepsilon + \mu)$,
\item[(d)] $\displaystyle prec(\rho) + 2 \mu \leq \, prec(\bar{\rho}) \leq
\, acc(\rho) + 2(\varepsilon + \mu)$,
\item[(e)] $\displaystyle \mid acc(\bar{\rho}) - \, acc(\rho) \mid \leq
2(\varepsilon + \mu)$. 
\end{enumerate}
\label{propacc1}
\end{proposition}

\paragraph{Proof.} Remark that (a), (b) are immediate consequences of definition
\ref{defacc} and that 
(c) and (d) must have identical proofs, just by switching $\varepsilon$ with 
$\mu$ and $X$ with $Y$ respectively. I shall prove therefore (c) and (e). 

For proving (c), consider $y \in Y$. By definition of 
$\displaystyle \bar{\rho}$ we write 
$$\left\{ x \in X \mbox{ : } (x,y) \in \bar{\rho} \right\} \, = \, \bigcup_{(x',y')
\in \rho , y' \in \bar{B}(y,\mu)} \bar{B}(x',\varepsilon)$$
Therefore we get 
$$res(\bar{\rho})(y) \, \geq \, 2 \varepsilon + \sup \left\{ res(\rho)(y') 
\mbox{ : } y' \in \,
im(\rho) \cap \bar{B}(y,\mu) \right\} $$
By taking the supremum over all $y \in Y$ we obtain the inequality 
$$res(\rho) + 2 \varepsilon \leq \, res(\bar{\rho})$$
For the other inequality, let us consider $\displaystyle (x_{1},y), (x_{2},y)
\in \bar{\rho}$ and $\displaystyle (x_{1}', y_{1}'), (x_{2}', y_{2}') \in \rho$
such that $\displaystyle d(x_{1},x_{1}') \leq \varepsilon, d(x_{2},x_{2}') 
\leq \varepsilon, D(y_{1}',y) \leq \mu,  D(y_{2}',y) \leq \mu$. Then: 
$$d(x_{1},x_{2}) \leq 2 \varepsilon + d(x_{1}',x_{2}') \leq 2 \varepsilon + \, 
acc(\rho) + d(y_{1}',y_{2}') \leq 2 (\varepsilon + \mu) + \, 
acc(\rho)$$
Take now a supremum and arrive to the desired inequality. 

For the proof of (e)  let us consider for $i = 1,2$   
$\displaystyle (x_{i},y_{i}) \in \bar{\rho}, (x_{i}',y_{i}') \in \rho$ such 
that $\displaystyle d(x_{i}, x_{i}') \leq \varepsilon, D(y_{i},y_{i}') \leq
\mu$.   It is then  enough to take absolute values and transform 
the following equality  
$$d(x_{1},x_{2}) - D(y_{1},y_{2}) = d(x_{1},x_{2}) - d(x_{1}',x_{2}') + 
d(x_{1}',x_{2}') - D(y_{1}',y_{2}') + $$ 
$$+ D(y_{1}',y_{2}') - D(y_{1},y_{2})$$ 
into well chosen, but straightforward, inequalities. \quad $\square$

The Gromov-Hausdorff distance is simply  the optimal lower bound for the accuracy of maps of $(X,d)$ into $(Y,D)$.

\begin{definition}
Let $\displaystyle (X, d)$, $(Y,D)$,  be a pair of metric spaces and $\mu > 0$. 
We shall say that $\mu$ is admissible if  there is a relation 
$\displaystyle \rho \subset X \times Y$ such that  
$\displaystyle dom \ \rho = X$, $\displaystyle im \ \rho = Y$, and $acc(\rho) \leq \mu$.
The Gromov-Hausdorff distance  between $\displaystyle (X,d)$ and $\displaystyle  
(Y,D)$  is   the infimum of admissible numbers $\mu$. 
\label{defgh}
\end{definition}

This is a true distance on the set of  isometry classes of pointed compact metric spaces.

\subsection{Length in metric spaces}
\label{slen}

In a metric space, the distance function associates a number to a pair of points. In length metric spaces we have a length functional defined over Lipschitz curves. This functional is defined over a space of curves, therefore is a more sophisticated object. Length dilation structures are the correspondent of dilation structures for length metric spaces. It is surprising that, as far as I know, before \cite{buligachar} there were no previous efforts to describe the behaviour of the length functional restricted to smaller and smaller regions around a point in a length metric space. 

The following definitions and results are standard, see for example the first chapter of  \cite{amb}.

\begin{definition}
The  (upper) dilation of a map $\displaystyle f: (X, d) \rightarrow (Y,D)$,  
in a point $u \in Y$ is 
$$ Lip(f)(u) = \limsup_{\varepsilon \rightarrow 0} \  
\sup  \left\{ 
\frac{D(f(v), f(w))}{d(v,w)} \ : \ v \not = w \ , \ v,w \in B(u,\varepsilon)
 \right\}$$
\end{definition}
Clearly, in the  case of a   derivable function 
$f: \mathbb{R} \rightarrow \mathbb{R}^{n}$ the upper dilation is  
$\displaystyle Lip(f)(t)  =  \|\dot{f}(t)\|$. 

A function  $f:(X, d) \rightarrow (Y, D)$  is Lipschitz if there is a positive 
 constant $C$ such that for any $x,y \in X$ we have 
 $\displaystyle D(f(x),f(y)) \leq C \, d(x,y)$. The number $Lip(f)$ is the smallest 
 such positive constant. Then  for any $x \in X$ we have the obvious relation  
$\displaystyle Lip(f)(x) \ \leq \ Lip(f)$.

A curve is a continuous function $c: [a,b] \rightarrow X$. The image of a curve is 
called "path". Geometrically speaking, length measures paths, not curves, therefore the length  functional, if defined over a class of curves, should be independent on the reparameterization of the path (image of the curve). 

\begin{definition}
Let $(X,d)$ be a metric space. There are  several ways to define  a notion of length.   The length of a  curve with $L^{1}$ upper dilation $c: [a,b] \rightarrow X$ is 
$\displaystyle L(f) = \int_{a}^{b} Lip(c)(t) \mbox{ d}t$. The  variation of any curve  $c: [a,b] \rightarrow X$ is   
$$Var(c)  = \sup \left\{ \sum_{i=0}^{n} d(c(t_{i}), 
c(t_{i+1})) \ \mbox{ : }  a = t_{0} <  t_{1} < ... < t_{n} < t_{n+1} =  b
\right\}$$ 
The length of the path $A = c([a,b])$ is  the
one-dimensional Hausdorff measure of the path, i.e. 
$$l(A) \ = \ \lim_{\delta \rightarrow 0}  
 \inf \left\{ \sum_{i \in I} diam \ E_{i}  \mbox{ : } diam \ E_{i} 
< \delta \ , \ \ A \subset \bigcup_{i \in I} E_{i} \right\} $$
\label{deflenght}
\end{definition}

The definitions are not equivalent, but  for any Lipschitz curve $c: [a,b] \rightarrow X$, we have 
$$\displaystyle L(c) \ = \ Var(c) \ \geq \ \mathcal{H}^{1}(c([a,b]))$$ 
Moreover, if  $c$ is  injective (i.e. a simple curve) then  $\displaystyle \mathcal{H}^{1}(c([a,b])) \  =  \ Var(f)$. 

It is important to know the fact (which will be used repeatedly) that any Lipschitz curve  $c$ admits a reparametrisation $c'$ such that $Lip(c')(t) = 1$ for almost any $t \in [a,b]$.

We associate a length functional to a metric space. 

\begin{definition}
We shall denote by $l_{d}$ the length functional induced by the distance $d$, 
defined only on the family of Lipschitz curves.   If the metric space $(X,d)$ is connected by Lipschitz curves, then the length induces a new distance $d_{l}$, given by: 
$$d_{l}(x,y) \  = \ \inf \ \left\{ l_{d}(c([a,b])) \mbox{ : } 
c: [a,b] \rightarrow X \ \mbox{ Lipschitz } , \right.$$
$$\left. \ c(a)=x \ , \ c(b) = y \right\}$$

A length metric space is a metric space $(X,d)$, connected by Lipschitz curves, such that $d  = d_{l}$. 
\label{dpath}
\end{definition}

Lipschitz curves in complete length metric spaces are absolutely continuous. 
Indeed, here is the definition of an absolutely continuous curve 
(definition 1.1.1, chapter 1,   \cite{amb}). 
 
 \begin{definition}
 Let $(X,d)$ be a complete metric space. A curve $c:(a,b)\rightarrow X$ is absolutely 
 continuous if there exists $m\in L^{1}((a,b))$ such that for any $a<s\leq t<b$ we have 
 $$d(c(s),c(t)) \leq \int_{s}^{t} m(r) \mbox{ d}r   .$$
 Such a function $m$ is called a  upper gradient of the curve $c$. 
 \label{defac}
 \end{definition}
 
 For a Lipschitz curve $c:[a,b]\rightarrow X$ in a 
 complete length metric space such a function 
 $m\in L^{1}((a,b))$  is the upper dilation  $Lip(c)$. 
More can be said about the expression of the upper dilation. We need first to introduce the notion of metric derivative of a Lipschitz curve. 

\begin{definition}
A curve $c:(a,b)\rightarrow X$ is  metrically derivable in $t\in(a,b)$ if the limit 
$$md(c)(t) = \lim_{s\rightarrow t} \frac{d(c(s),c(t))}{\mid s-t \mid}$$
exists and it is finite. In this case $md(c)(t)$ is called the  metric
derivative of $c$ in $t$. 
\label{defmd}
\end{definition}

For the proof of the following theorem see \cite{amb}, theorem 1.1.2, chapter 1. 

\begin{theorem}
Let $(X,d)$ be a complete metric space and $c:(a,b)\rightarrow X$ be an absolutely continuous curve. 
Then $c$ is metrically  derivable for $\mathcal{L}^{1}$-a.e. $t\in(a,b)$. Moreover the function $md(c)$ belongs to $L^{1}((a,b))$ and it is minimal in the following sense: $md(c)(t)\leq m(t)$  for  $\mathcal{L}^{1}$-a.e. $t\in(a,b)$, for each upper gradient $m$ of the curve $c$. 
\label{tupper}
\end{theorem}

\subsection{Metric profiles. Metric tangent space}

To any locally compact metric  space we  associate a metric profile \cite{buliga3,buliga4}. 
This metric profile is a way of organizing the information given by the distance function in order to get an understanding of the local behaviour of the distance around a point of the space. We need to consider local compactness in order to have compact small balls in the next definition. 

Let us denote by $CMS'$ the set of isometry classes of pointed compact metric spaces. An element of $CMS'$ is denoted like $[X,d,x]$ and is the equivalence class of a compact metric space $(X,d)$ with a specified point $x \in X$, with respect to the following equivalence relation: two pointed compact metric spaces $(X,d,x)$ and $(Y,D,y)$  are equivalent if there is a surjective isometry $f: (X,d) \rightarrow (Y,D)$  such that $f(x) = y$. 

The set $CMS'$ is endowed with the GH distance for pointed metric spaces. In order to define this distance we have to slightly modify definition \ref{defgh}, by restricting the class of maps (relations) of $(X,d)$ into $(Y,D)$ -- $\displaystyle \rho \subset X \times Y$ such that  
$\displaystyle dom \ \rho = X$, $\displaystyle im \ \rho = Y$ -- to those which satisfy also 
$(x,y) \in \rho$. 

We can define now metric profiles. 

\begin{definition}
\label{dmprof}
The metric profile associated to the locally metric space $(M,d)$ is  the assignment 
(for small enough $\varepsilon > 0$) 
$$(\varepsilon > 0 , \ x \in M) \ \mapsto \  \mathbb{P}^{m}(\varepsilon, x) = \left[\bar{B}(x,1), 
\frac{1}{\varepsilon} d, x\right] \in CMS' $$
\end{definition}

The metric profile of the space at a point is therefore  a curve in another metric space, namely $CMS'$. with a Gromov-Hausdorff distance. It is not any curve, but one which has certain properties which can be expressed with the help of the GH distance.  Indeed, for any $\varepsilon, b >0$,  sufficiently small, we have 
$$
\mathbb{P}^{m}(\varepsilon b, x) =  \mathbb{P}^{m}_{d_{b}}(\varepsilon,x)
$$
where  $d_{b} = (1/b) d$ and 
$\mathbb{P}^{m}_{d_{b}}(\varepsilon,x) = [\bar{B}(x,1),\frac{1}{\varepsilon}d_{b}, x]$.

These curves give interesting local and infinitesimal information about the metric space. For example, what kind of metric space has constant metric profile with respect to one of its points?

\begin{definition}
\label{defmetcone}
A metric cone $(X,d, x)$ is a locally compact metric space $(X,d)$, with a marked point $x \in X$ such 
that for any $a,b \in (0,1]$ we have 
\[\displaystyle \mathbb{P}^{m}(a,x)  =  \mathbb{P}^{m}(b,x)\] 
\end{definition}

Metric cones are self-similar, in the sense that they have dilations.

\begin{definition}
Let $(X,d, x)$ be a metric cone. For any $\varepsilon \in (0,1]$  a dilation is a function $\displaystyle \delta^{x}_{\varepsilon}: \bar{B}(x,1) \rightarrow \bar{B}(x,\varepsilon)$ such that 
\begin{itemize}
\item[-] $\displaystyle \delta^{x}_{\varepsilon}(x) = x$, 
\item[-] for any $u,v \in X$ we have 
$$d\left(\delta^{x}_{\varepsilon}(u), \delta^{x}_{\varepsilon}(v)\right) =\varepsilon \, d(u,v) $$
\end{itemize}
\end{definition}

The existence of dilations for metric cones comes from the definition \ref{defmetcone}. 
Indeed, dilations are just  isometries from $\displaystyle (\bar{B}(x,1), d, x)$ to $ (\bar{B}, \frac{1}{a}d, x)$.

\begin{definition}
\label{defmetspace}
A (locally compact) metric space $(M,d)$ admits a (metric) tangent space in $x \in M$ if the associated metric profile  
$\varepsilon \mapsto \mathbb{P}^{m}(\varepsilon, x)$  admits a prolongation by continuity in 
$\varepsilon = 0$, i.e if the following limit exists: 
\begin{equation}
\label{limmetspace}
[T_{x}M,d^{x}, x]  =  \lim_{\varepsilon \rightarrow 0} \mathbb{P}^{m}(\varepsilon, x)
\end{equation}
\end{definition}

The connection between metric cones, tangent spaces and metric profiles in the abstract sense is made by the following proposition. 

\begin{proposition}
\label{propmetcone}
Metric tangent spaces are  metric cones. 
\end{proposition}

\paragraph{Proof.} 
A tangent space  $[V,d_{v}, v]$ exists if and only if we have the limit from the relation (\ref{limmetspace}), that is of there exists a prolongation by continuity to $\varepsilon = 0$  of the  metric profile $\displaystyle \mathbb{P}^{m}(\cdot , x)$.  For any $a \in (0,1]$ we have
$\displaystyle \left[\bar{B}(x,1), \frac{1}{a}d^{x}, x\right] = \lim_{\varepsilon \rightarrow 0} \mathbb{P}^{m}(a \varepsilon, x)$, therefore we have 
$$\left[\bar{B}(x,1), \frac{1}{a}d^{x}, x\right]  = [T_{x}M,d^{x}, x] $$ 
which proves the thesis. \quad $\square$

We may also define abstract metric profiles. The previously defined metric profiles are abstract metric profiles, but we shall see further (related to the Mitchell theorem in sub-riemannian geometry, for example) that abstract metric profiles are useful too.  

\begin{definition}
\label{dprofile}
An abstract  metric profile is a curve $\mathbb{P}:[0,a] \rightarrow CMS$ such that
\begin{enumerate}
\item[(a)] it is continuous at $0$,
\item[(b)]for any $b \in [0,a]$ and  $\varepsilon \in (0,1]$ we have 
\[d_{GH} (\mathbb{P}(\varepsilon b), \mathbb{P}^{m}_{d_{b}}(\varepsilon,x_{b}))  \  = \ O(\varepsilon)\]
\end{enumerate}
The function $\mathcal{O}(\varepsilon)$ may change with  $b$.
We used the notations
\[
\mathbb{P}(b) = [\bar{B}(x,1) ,d_{b}, x_{b}] \quad \mbox{  and } \quad 
\mathbb{P}^{m}_{d_{b}}(\varepsilon,x) = \left[\bar{B}(x,1),\frac{1}{\varepsilon}d_{b}, x_{b}\right] 
\]
\end{definition}

\subsection{Curvdimension and curvature}
\label{securv}

In the case of a riemannian manifold $(X,g)$, with smooth enough (typically $\displaystyle 
\mathcal{C}^{1}$) metric $g$, the tangent metric spaces exist for any point of the manifold. 
Ideed, the tangent metric space $\displaystyle [T_{x}X,d^{x}, x]$ is the isometry class of a small neighbourhood of the origin of the tangent space (in differential geometric sense) $\displaystyle T_{x}X$, with $\displaystyle d^{x}$ being the euclidean distange induced by the norm given by $\displaystyle g_{x}$.  Moreover, we have the following description of the sectional curvature. 

\begin{proposition}
Let $(X,d)$ be a $\displaystyle \mathcal{C}^{4}$ smooth riemannian manifold with $d$ the length distance induced by the riemannian metric $g$. Suppose that for a point $x \in X$ the sectional curvature is bounded in the sense that for any 
linearly independent $\displaystyle u, v \in T_{x}X$ we have $\displaystyle \mid K_{x}(u,v) \mid \leq C$. Then for any sufficiently small $\varepsilon > 0$ we have 
\begin{equation}
\frac{1}{\varepsilon^{2}} d_{GH}(\mathbb{P}^{m}(\varepsilon,x), [T_{x}X,d^{x}, x]) \, \leq \, \frac{1}{3} C + \mathcal{O}(\varepsilon)
\label{sectional}
\end{equation}
\label{pcurvature}
\end{proposition}

\paragraph{Proof.}
This is well known, in another form. Indeed, for small enough $\varepsilon$, consider the geodesic exponential map which associates to any $\displaystyle u \in W \subset T_{x}X$ (in a neighbourhood $W$ of the origin which is independent of $\varepsilon$) the point 
$\exp_{x} \varepsilon u$.  Define now the distance 
$$d_{\varepsilon}^{x} (u, v) \, = \, \frac{1}{\varepsilon} d (\exp_{x} \varepsilon u , 
\exp_{x} \varepsilon v)$$
We can choose the neighbourhood $W$ to be the unit ball with respect to the distance $\displaystyle d^{x}$ in order to get  the following estimate: 
$$\sup \left\{ \mid d_{\varepsilon}^{x} (u, v) - d^{x}(u,v)\mid \mbox{ : } 
u,v \in W \right\} \geq d_{GH}(\mathbb{P}^{m}(\varepsilon,x), [W,d^{x}, 0])$$ 
where $\displaystyle d^{x}(u,v) = \|u-v\|_{x}$. 
In the given regularity settings, we shall use the following expansion of $\displaystyle 
d^{x}_{\varepsilon}$: if $u,v$ are linearly independent then 
\begin{equation}
d^{x}_{\varepsilon}(u,v) \, = \, d^{x}(u,v) - \frac{1}{6} \varepsilon^{2} K_{x}(u,v) 
\frac{\|u\|_{x}^{2} \|v\|_{x}^{2} - \langle u,v \rangle_{x}^{2}}{d^{x}(u,v)} + \varepsilon^{2} \mathcal{O}(\varepsilon)
\label{mainriem}
\end{equation}
where $K$ is the sectional curvature of the metric $g$. (If $u,v$ are linearly dependent then 
$\displaystyle  d^{x}_{\varepsilon}(u,v)  =  d^{x}(u,v)$.) From here we easily obtain that 
$$\frac{1}{\varepsilon^{2}} d_{GH}(\mathbb{P}^{m}(\varepsilon,x), [T_{x}X,d^{x}, x]) \, \leq \, \frac{1}{3} \sup\left\{ \mid K_{x}(u,v) \mid \mbox{ : } u,v \in W \mbox{ lin. indep.}\right\} + \mathcal{O}(\varepsilon)$$ 
which ends the proof. \quad $\square$

This proposition makes us define the "curvdimension" and "curvature" of an (abstract) metric profile. 

\begin{definition}
Let $\displaystyle \mathbb{P}$ be an abstract metric profile. The curvdimension of this 
abstract metric profile is 
\begin{equation}
curvdim \, \mathbb{P} \, = \, \sup \left\{ \alpha > 0 \mbox{ : } \lim_{\varepsilon \rightarrow 0} \frac{1}{\varepsilon^{\alpha}} d_{GH}(\mathbb{P}(\varepsilon), \mathbb{P}(0)) = 0 \right\} 
\label{dcurvdim}
\end{equation}
and, in the case that the curvdimension equals $\beta > 0$ then the $\beta$-curvature of $\mathbb{P}$ is the number $M > 0$ such that 
\begin{equation}
\lim_{\varepsilon \rightarrow 0} \log_{\varepsilon} \left( \frac{1}{M}  d_{GH}(\mathbb{P}(\varepsilon), \mathbb{P}(0)) \right) = \beta 
\label{dcurv}
\end{equation}
In case $\mathbb{P}$ is the metric profile of a point $x \in X$ in a metric space $(X,d)$ then the curvdimension at $x$ and curvature at $x$ are the curvimension, respectively curvature, of the metric profile of $x$. 
\label{defcurv}
\end{definition}

It follows that non-flat riemannian (smooth enough) spaces have curvdimension 2. Also, any metric cone has curvdimension equal to $0$ (meaning "all metric cones are flat objects"). In particular, finite dimensional normed vector spaces are flat (as they should be).

\paragraph{Question.} What is the curvdimension of a sub-riemannian space? Carnot groups endowed with left invariant Carnot-Carath\'eodory distances (or any left-invariant distance coming from a norm on the Carnot group seen as a conical group) have curvdimension equal to zero, i.e. they are "flat". In the non-flat case, that is when the metric profiles are not constant, what is then the curvdimension at a generic point of a sub-riemannian space? In the paper \cite{buliga3}, theorem 10.1, then in section 8, \cite{buliga4}, we presented evidence for the fact that metric contact 3 dimensional manifolds have curvdimension smaller than 2. Is it, in this case, equal to 1?

\section{Groups with dilations}

We shall see that for a dilation structure (or "dilatation structure", or "metric space with dilations") the metric tangent spaces   have the structure of a normed local group with dilations. The notion has been introduced in published \cite{buliga2}, \cite{buligadil1}; it appears first time in \cite{buligasrlie1}, which we follow here.

 We shall work with local groups and local actions instead of the usual global ones. We shall use "uniform local group" for a local group endowed with its canonical  uniform structure.

Let $\Gamma$ be a topological commutative group, endowed with a continuous morphism 
$\mid \cdot \mid : \Gamma \rightarrow (0,+\infty)$. For example $\Gamma$ could be $(0,+\infty)$ with the operation of multiplication of positive real numbers and the said morphism could be the identity. Or $\Gamma$ could be the set of complex numbers different from $0$, with the operation of multiplication of complex numbers and morphism taken to be the modulus function. 
Also, $\Gamma$ could be the set of integers with the operation of addition and the morphism could be the exponential function. Many other possibilities exist (like a product between a finite commutative group with one of the examples given before). 

It is useful further to just think that $\Gamma$ is like in the first example, because in these notes we are not going to use the structure of $\Gamma$ in order to put more geometrical objects on the metric space (like we do, for example, in the paper \cite{buligabraided}).

The elements of $\Gamma$ will be denoted with small greek letters, like $\varepsilon, \mu, ...$. By covention, whenever we write "$\varepsilon \rightarrow 0$", we really mean "$\mid \varepsilon \mid \rightarrow 0$". Also "$\displaystyle \mathcal{O}(\varepsilon)$" means  "$\displaystyle \mathcal{O}(\mid\varepsilon\mid)$", and so on.

\begin{definition}
$(G,\delta)$ is a local  group with dilations if  $G$ is a local group and  $\displaystyle \delta : \Gamma \rightarrow \mathcal{C}(G,G)$ is a local  action of $\Gamma$  on $G$, such that
\begin{enumerate}
\item[H0.] the following limit $\displaystyle \lim_{\varepsilon \rightarrow 0} 
\delta_{\varepsilon} x  =  e$  is uniform with respect to $x$ in a compact neighbourhood of the identity element  $e$.
\item[H1.] the limit $ \displaystyle \beta(x,y)  =  \lim_{\varepsilon \rightarrow 0} \delta_{\varepsilon}^{-1} \left((\delta_{\varepsilon}x) (\delta_{\varepsilon}y ) \right)$
is is uniform with respect to $(x,y)$ in a compact neighbourhood of $(e,e)$.
\item[H2.] the limit $\displaystyle  \lim_{\varepsilon \rightarrow 0} \delta_{\varepsilon}^{-1}
\left( ( \delta_{\varepsilon}x)^{-1}\right)  =  x^{-1}$ is uniform with respect to $x$ in a compact neighbourhood of the identity element  $e$.
\end{enumerate}
\label{defgwd}
\end{definition}

\begin{definition} A normed local group with dilations $(G,  \| \cdot \|, \delta)$ is a 
local group with dilations  $(G, \delta)$ endowed with a continuous norm,   
 $\displaystyle \|\cdot \| : G \rightarrow [0,+\infty)$  which satisfies: 
 \begin{enumerate}
 \item[(a)] there is a function  $\displaystyle \|\cdot \| : U \subset G \rightarrow [0,+\infty)$ defined on a neighbourhood $U$ of $e$, such that the limit $\displaystyle \lim_{\varepsilon \rightarrow 0} \frac{1}{\mid \varepsilon \mid} \| \delta_{\varepsilon} x \| = \| x\|^{N}$  is uniform with respect to $x$ in compact set, 
 \item[(b)] if $\displaystyle \| x\|^{N} = 0$ then $x=e$.
  \end{enumerate}
In a normed local group with dilations we consider the  left invariant (locally defined) distance given by 
\begin{equation}
d(x,y) = \| x^{-1}y\| \quad . 
\label{dnormed}
\end{equation}
and  dilations based in any point $x \in G$ by 
 \begin{equation}
 \delta^{x}_{\varepsilon} u = x \delta_{\varepsilon} ( x^{-1}u)  . 
 \label{dilat}
 \end{equation}
  \label{dnco}
  \end{definition}

\subsection{Conical groups}

\begin{definition}
A normed conical group $(N, \| \cdot \|, \delta)$ is a normed  group with dilations  such that  for any $\varepsilon \in \Gamma$: (a)  the dilation 
 $\delta_{\varepsilon}$ is a group morphism  and  (b) the norm is homogeneous, that is 
  $\displaystyle  \| \delta_{\varepsilon} x \| = \mid \varepsilon\mid \,  \| x \|$. 
\end{definition}

Then, a  conical group appears as the  infinitesimal version of a group with 
dilations (\cite{buligadil1} proposition 2). For the proof see the more general theorem 
\ref{tgene}.

\begin{proposition}
Let $(G, \| \cdot \|, \delta$ be a normed local group with dilations. Then $\displaystyle (G, \| \cdot \|^{N},  \delta)$ is a local normed conical group, with operation 
$\beta$,  dilations $\delta$ and homogeneous norm $\displaystyle \| \cdot \|^{N}$.
\label{here3.4}
\end{proposition}

\subsection{Carnot groups}
\label{carnotgroups}

Carnot groups appear in sub-riemannian geometry 
 as models of tangent spaces,   \cite{bell}, \cite{gromovsr}, \cite{pansu}. In particular such groups can be endowed with a structure of a sub-riemannian manifold. Here we are interested in the fact that they are particular examples of conical groups.

\begin{definition}
A Carnot group is a pair $\displaystyle (N, V_{1})$ formed by a real connected simply connected Lie group $N$,   with  a distinguished subspace  
$V_{1}$ of  the Lie algebra $Lie(N)$, such that  
$$Lie(N) \ = \ \sum_{i=1}^{m} V_{i} \ , \ \ V_{i+1} \ = \ [V_{1},V_{i}]$$
The number $m$ is called the step of the group. T number $\displaystyle Q \ = \ \sum_{i=1}^{m} i \ dim V_{i}$ is called the homogeneous dimension of the group. 
\label{dccgroup}
\end{definition}

Because the group $N$ is nilpotent and simply connected, it follows that the (Lie group) 
exponential mapping is a diffeomorphism. It is customary to identify then  the group with the algebra.  We obtain a set $N$ equal to some $\displaystyle \mathbb{R}^{n}$,  endowed with a a Lie
algebra structure (that is a real vector space and a Lie bracket) and a Lie group structure (that is a Lie group operation, denoted multiplicatively,  defined for any pair of elements of $N$, with the $0$ element of the vector space $N$ as neutral element for the group operation). 

The  Baker-Campbell-Hausdorff formula connects the Lie bracket and the group operation. 
Indeed, the algebra being nilpotent, it follows that the group operation
is polynomial, because the Baker-Campbell-Hausdorff formula contains only a finite 
number of terms. Thus, the group operation is expressed as a function of the Lie bracket operation.  Moreover,  Lie algebra endomorphisms are group endomorphisms and th converse is also true.

The most simple case is when the Lie bracket is a constant function equal to $0$ and $\displaystyle  V_{1} = N = \mathbb{R}^{n}$. In this case the group operation is the vector space addition (even if defined multiplicatively). The step of this group is equal to 1 and the homogeneous dimension equals $n$. 

Let us take $\Gamma = (0, +\infty)$, $\mid \cdot \mid : \Gamma \rightarrow (0, \infty)$ the dentity function  and let us define  for any 
$\varepsilon > 0$, the dilation: 
$$ x \ = \ \sum_{i=1}^{m} x_{i} \ \mapsto \ \delta_{\varepsilon} x \
= \ \sum_{i=1}^{m} \varepsilon^{i} x_{i}$$
Any such dilation is a group morphism and a Lie algebra morphism.

Let us choose  an  euclidean norm $\| \cdot \|$ on $\displaystyle V_{1}$.  We shall endow the group $N$ with a norm coming from a Carnot-Carth\'eodory distance  (general definition in section \ref{secsr}). Remark that for any $\displaystyle x \in V_{1}$ and any $\varepsilon > 0$ we have $\displaystyle \| \delta_{\varepsilon} x \| = \varepsilon \| x \|$.

Indeed, by definition the space  $V_{1}$  generates  $N$ (as a Lie algebra), therefore  any element $x \in N$ can be written as a product of elements from $V_{1}$.  A controlled way to do this is  described in the following slight reformulation of Lemma 1.40, Folland, Stein \cite{folstein}).  

\begin{lemma}
Let $N$ be a Carnot group and $X_{1}, ..., X_{p}$ an orthonormal basis 
for $V_{1}$. Then there is  a natural number $M$ and a function $g: \left\{ 1,...,M \right\} 
\rightarrow \left\{ 1,...,p\right\}$ such that any 
$x \in N$ can be written as: 
\begin{equation}
x \ = \ \prod_{i = 1}^{M} \exp(t_{i}X_{g(i)})
\label{fp2.4}
\end{equation}
Moreover, if $x$ is sufficiently close (in Euclidean norm) to
$0$ then each $t_{i}$ can be chosen such that $\mid t_{i}\mid \leq C 
\| x \|^{1/m}$
\label{p2.4}
\end{lemma}

From these data we may construct a norm on the Carnot group $N$, by the intermediary of a Carnot-Carath\'eodory (CC for short) distance. Here we give an algebraic definition of this distance.

 The (Carnot-Carath\'eodory) norm  on the Carnot group is defined as  
$$\| x \| \ = \ \inf \left\{ \sum_{i \in I} \| x_{i} \|  \ \mbox{ :  all finite sets } I \mbox { and all decompositions }  \quad x = \prod_{i \in I}  x_{i} \mbox{ where all } x_{i} \in V_{1} \right\}$$ 
The CC norm is then  finite (by lemma \ref{p2.4}) for any two $x \in N$ and it is also continuous. All in all $(N, \| \cdot \|, \delta)$ is a (global) normed conical group.

\subsection{Contractible groups}

\begin{definition}
A contractible group is a pair $(G,\alpha)$, where $G$ is a  
topological group with neutral element denoted by $e$, and $\alpha \in Aut(G)$ 
is an automorphism of $G$ such that: 
\begin{enumerate}
\item[-] $\alpha$ is continuous, with continuous inverse, 
\item[-] for any $x \in G$ we have the limit $\displaystyle 
\lim_{n \rightarrow \infty} \alpha^{n}(x) = e$. 
\end{enumerate}
\end{definition}

For a contractible group  $(G,\alpha)$, the automorphism $\alpha$ is compactly contractive (Lemma 1.4 (iv) \cite{siebert}), that is: for each compact set 
$K \subset G$ and open set $U \subset G$, with $e \in U$, there is 
$\displaystyle N(K,U) \in \mathbb{N}$ such that for any $x \in K$ and $n \in
\mathbb{N}$, $n \geq N(K,U)$, we have $\displaystyle \alpha^{n}(x) \in U$. 

If $G$ is locally compact then $\alpha$ compactly contractive is equivalent with:  each identity neighbourhood of $G$ contains an $\alpha$-invariant neighbourhood. Further on we shall assume without mentioning that all groups are locally compact.

Any conical group is a contractible  group. 
Indeed, it suffices to associate to a conical group $(G,\delta)$ the 
contractible group $\displaystyle (G,\delta_{\varepsilon})$, for a fixed 
$\varepsilon \in \Gamma$ with $\nu(\varepsilon) < 1$.

Conversely, to any contractible group $(G,\alpha)$ we may  associate the conical group $(G,\delta)$, with 
$\displaystyle \Gamma = \left\{ \frac{1}{2^{n}} \mbox{ : } n \in \mathbb{N}
\right\}$ and for any $n \in \mathbb{N}$ and $x \in G$ 
$$\displaystyle \delta_{\frac{1}{2^{n}}} x \ = \ \alpha^{n}(x) \quad . $$
(Finally, a local conical
group has only locally the structure of a contractible group.)

The structure of contractible groups is known to some detail, due to results by Siebert
\cite{siebert}, Wang \cite{wang},  Gl\"{o}ckner and Willis \cite{glockwill}, 
Gl\"{o}ckner \cite{glockner} and others (see  references in the mentioned papers). 

Related to contractible groups, here is the 
definition of a contractive automorphism group \cite{siebert}, definition 5.1. 

\begin{definition}
Let $G$ be a locally compact group. An automorphism group on $G$ is a family 
$\displaystyle T= \left( \tau_{t}\right)_{t >0}$ in $Aut(G)$, such that 
$\displaystyle \tau_{t} \, \tau_{s} = \tau_{ts}$ for all $t,s > 0$. 

The contraction  group of $T$ is defined by 
$$C(T) \ = \ \left\{ x \in G \mbox{ : } \lim_{t \rightarrow 0} \tau_{t}(x) = e
\right\} \quad .$$
The automorphism group $T$ is contractive if $C(T) = G$. 
\end{definition}

It is obvious that a contractive automorphism  group $T$ induces on $G$ a 
structure of conical group. Conversely, any conical group with $\Gamma = 
(0,+\infty)$ has an associated contractive automorphism group (the group of 
dilations based at the neutral element). 

Siebert, proposition 5.4 \cite{siebert}, gives a very useful description of a class of contractible groups. 

\begin{proposition}
For a locally compact group $G$ the following assertions are equivalent: 
\begin{enumerate}
\item[(i)] $G$ admits a contractive automorphism group;
\item[(ii)] $G$ is a simply connected Lie group whose Lie algebra admits a 
positive graduation.
\end{enumerate}
\label{psiebert}
\end{proposition}

These groups are almost Carnot groups. Indeed, what is missing is the fact that the first elemet of the graduation generates the Lie algebra.

\section{Dilation structures}
\label{dilst}

In this paper I use the denomination "dilation structure", or "metric space with dilations", compared with older papers, where the name "dilatation structure" was used. 

We shall use here a slightly particular version of dilation structures. 
For the general definition of a dilation structure see \cite{buligadil1} (the general definition applies for dilation structures over ultrametric spaces as well). 

\subsection{Normed groupoids with dilations}

\paragraph{Notions of convergence.}
We need a topology on a normed groupoid $(G,d)$, induced by the norm. 

%For this reason, we
%shall always denote the limit of a net $\displaystyle (s_{\varepsilon})$, 
%$\displaystyle s_{\varepsilon} \in \mathbb{R}$ for any $\varepsilon \in I$, 
% by $\displaystyle \lim_{\varepsilon \rightarrow 0} s_{\varepsilon}$ instead of 
% $\displaystyle \lim_{\varepsilon \in I} s_{\varepsilon}$.

\begin{definition}
 A net of arrows  $\displaystyle (a_{\varepsilon}) $ 
 simply converges to the arrow $a \in G$ (we write $\displaystyle a_{\varepsilon}
\rightarrow a$) if: 
\begin{enumerate}
\item[(i)] for any $\displaystyle \varepsilon \in I$ there are elements 
$\displaystyle g_{\varepsilon}, h_{\varepsilon} \in G$ such that $\displaystyle h_{\varepsilon} a_{\varepsilon} g_{\varepsilon} = 
a$, 
\item[(ii)] we have $\displaystyle \lim_{\varepsilon \in I} d(g_{\varepsilon}) =
0$ and $\displaystyle \lim_{\varepsilon \in I} d(h_{\varepsilon}) = 0$. 
\end{enumerate}

 A net of arrows  $\displaystyle (a_{\varepsilon})$ left-converges 
 to the arrow $a \in G$ (we write 
$\displaystyle a_{\varepsilon} \lrightarrow a$) if for all  $i \in I$ we have 
$\displaystyle (a_{\varepsilon}^{-1}, a) \in G^{(2)}$ and 
moreover $\displaystyle \lim_{\varepsilon \in I} d(a_{\varepsilon}^{-1} a) = 0$. 

 A net of arrows  $\displaystyle (a_{\varepsilon})$ right-converges 
 to the arrow $a \in G$ (we write 
$\displaystyle a_{\varepsilon} \rrightarrow a$) if for all  $i \in I$ we have 
 $\displaystyle (a_{\varepsilon}, a^{-1}) \in G^{(2)}$ and 
moreover $\displaystyle \lim_{\varepsilon \in I} d(a_{\varepsilon} a^{-1}) = 0$. 
\label{defconv}
\end{definition}

It is clear that if $\displaystyle a_{\varepsilon} \lrightarrow a$ or 
$\displaystyle a_{\varepsilon} \rrightarrow a$ then 
$\displaystyle a_{\varepsilon} \rightarrow a$. 

Right-convergence of $\displaystyle a_{\varepsilon}$ to $a$ is just convergence 
of $\displaystyle a_{\varepsilon}$ to $a$ in the distance $\displaystyle
d_{\alpha(a)}$, that is $\displaystyle \lim_{\varepsilon \in I}
d_{\alpha(a)}(a_{\varepsilon}, a)= 0$. 

Left-convergence of $\displaystyle a_{\varepsilon}$ to $a$ is just convergence 
of $\displaystyle a_{\varepsilon}^{-1}$ to $a^{-1}$ in the distance 
$\displaystyle
d_{\omega(a)}$, that is $\displaystyle \lim_{\varepsilon \in I}
d_{\omega(a)}(a_{\varepsilon}^{-1}, a^{-1})= 0$.

\begin{proposition}
\label{propconv}
Let $(G,d)($ be a normed groupoid.  
\begin{enumerate}
\item[(i)] If $\displaystyle a_{\varepsilon} \lrightarrow a$ and $\displaystyle a_{\varepsilon}
\lrightarrow b$ then $a=b$. If  $\displaystyle a_{\varepsilon} \rrightarrow a$ and $\displaystyle a_{\varepsilon}
\rrightarrow b$ then $a=b$. 
\item[(ii)] The following are equivalent: 
\begin{enumerate}
\item[1.] $G$ is a Hausdorff topological groupoid with respect to the topology induced by 
the simple convergence, 
\item[2.] $d$ is a separable norm, 
\item[3.] for any net $\displaystyle (a_{\varepsilon}) $, if $\displaystyle a_{\varepsilon} 
\rightarrow a$ and $\displaystyle a_{\varepsilon}
\rightarrow b$ then $a=b$.
\item[4.] for any net $\displaystyle (a_{\varepsilon}) $, if $\displaystyle a_{\varepsilon} 
\rrightarrow a$ and $\displaystyle a_{\varepsilon}
\lrightarrow b$ then $a=b$.
\end{enumerate}
\end{enumerate}
\end{proposition}

\paragraph{Proof.}
(i) We prove only the first part of the conclusion.  We can write $\displaystyle b^{-1} a = b^{-1} a_{\varepsilon} 
a_{\varepsilon}^{-1} a$, therefore 
$$d(b^{-1}a) \leq d(b^{-1} a_{\varepsilon}) + d(a_{\varepsilon}^{-1} a)$$
The right hand side of this inequality is arbitrarily small, so $\displaystyle 
d(b^{-1}a) =0$, which implies $a=b$. 

(ii) Remark that the structure maps are continuous with
respect to the topology induced by the simple convergence. We need only to prove
the uniqueness of limits. 

 3. $\Rightarrow$ 4. is trivial. In order to prove that 
4.$\Rightarrow$ 3., consider an arbitrary  net  
$\displaystyle (a_{\varepsilon}) $ such that  $\displaystyle a_{\varepsilon} 
\rightarrow a$ and $\displaystyle a_{\varepsilon}
\rightarrow b$. This means that there exist nets $\displaystyle 
(g_{\varepsilon}) , (g'_{\varepsilon}) , (h_{\varepsilon}) , 
(h'_{\varepsilon}) $ such that 
$\displaystyle h_{\varepsilon} a_{\varepsilon} g_{\varepsilon} = a$, 
$\displaystyle h'_{\varepsilon} a_{\varepsilon} g'_{\varepsilon} = b$ 
and $\displaystyle \lim_{i \in I} \left( d(g_{\varepsilon}) + 
d(g'_{\varepsilon}) + d(h_{\varepsilon}) + d(h'_{\varepsilon}) \right) = 0$. 
Let $\displaystyle g"_{\varepsilon} = g_{\varepsilon}^{-1} g'_{\varepsilon}$ 
and $\displaystyle h"_{\varepsilon} = h'_{\varepsilon} h_{\varepsilon}^{-1}$. 
We have then $\displaystyle b = h"_{\varepsilon} a g"_{\varepsilon}$ and 
$\displaystyle  \lim_{i \in I} \left( d(g"_{\varepsilon}) + d(h"_{\varepsilon})
\right) = 0$. Then $\displaystyle h"_{\varepsilon} a \lrightarrow b$ and 
$\displaystyle h"_{\varepsilon}a \rrightarrow a$. We deduce that $a=b$. 

1.$\Leftrightarrow$ 3. is trivial. So is 3. $\Rightarrow$ 2. 
We finish the proof by showing that 2. $\Rightarrow$ 3. 
By a reasoning made previously, it is enough to prove that: 
if $ b = h_{\varepsilon} a g_{\varepsilon}$ and 
$\displaystyle \lim_{i \in I} \left( d(g_{\varepsilon}) + d(h_{\varepsilon}) 
\right) = 0$ then $a=b$. Because $d$ is separable it follows that 
$\alpha(a) = \alpha(b)$ and $\omega(a) = \omega(b)$. We have then 
$\displaystyle a^{-1} b = a^{-1} h_{\varepsilon} a g_{\varepsilon}$, therefore 
$$d(a^{-1} b) \leq d(a^{-1} h_{\varepsilon} a) + d(g_{\varepsilon})$$
 The  norm $d$ induces a left invariant  distance on the vertex group of all 
 arrows $g$ such that $\alpha(g) = \omega(g) = \alpha(a)$. This distance is 
 obviously continuous with respect to the simple convergence  in the group. 
 The net $\displaystyle a^{-1} h_{\varepsilon} a$ simply converges to 
 $\alpha(a)$ by the continuity of the multiplication (indeed, 
 $\displaystyle h_{\varepsilon}$ simply converges to $\alpha(a)$). 
 Therefore  $\displaystyle \lim_{i \in I} d(a^{-1} h_{\varepsilon} a) = 0$.
It follows that $\displaystyle d(a^{-1}b)$ is arbitrarily small, therefore $a=b$. 
\quad $\quad \square$

By adapting the definition of a normed group with dilations to a normed groupoid with dilations, we get the following structure.

\begin{definition}
A  normed groupoid $(G,d, \delta)$ with dilations is a separated normed groupoid $(G,d)$ endowed with  a map assigning to any 
$\varepsilon \in \Gamma$ a transformation $\delta_{\varepsilon}: \, 
dom(\varepsilon) \rightarrow \, im(\varepsilon)$ which satisfies the following:
\begin{enumerate}
\item[A1.] For any 
$\varepsilon \in \Gamma$ $\alpha \delta_{\varepsilon} = \alpha$. Moreover $\displaystyle
\varepsilon \in \Gamma \mapsto \delta_{\varepsilon}$ is an action of $\Gamma$ on
$G$, that is for any 
$\varepsilon, \mu \in \Gamma$ we have $\displaystyle \delta_{\varepsilon} 
\delta_{\mu} = \delta_{\varepsilon \mu}$, $\displaystyle 
\left(\delta_{\varepsilon} \right)^{-1} = \delta_{\varepsilon^{-1}}$ and 
$\delta_{e} = \, id$. 
\item[A2.] For any $x \in Ob(G)$ and any $\varepsilon \in \Gamma$ we have 
$\displaystyle \delta_{\varepsilon}(x) = x$. Moreover the transformation 
$\displaystyle \delta_{\varepsilon}$ contracts $dom(\varepsilon)$ to $X = Ob(G)$ uniformly on bounded sets, which means that 
the net $\displaystyle d \, \delta_{\varepsilon}$ converges to the constant
function $0$, uniformly on bounded sets.
\item[A3.] There is a function 
$\displaystyle \bar{d}: U \rightarrow \mathbb{R}$ 
which is the limit  
$$ \lim_{\varepsilon \rightarrow 0} \frac{1}{\mid\varepsilon\mid} \, 
d \, \delta_{\varepsilon} (g)  \, = \, \bar{d}(g)$$
uniformly on bounded sets in the sense of 
definition \ref{defconv}. Moreover, if $\bar{d}(g) = 0$ then $g \in Ob(G)$.
\item[A4.] the net $\displaystyle dif_{\varepsilon}$ converges uniformly on
 bounded sets to a function $\bar{dif}$. 
\end{enumerate}
The domains and codomains of a dilation of $(G,d)$ satisfy the following 
Axiom A0: 
\begin{enumerate}
\item[(i)] for any $\varepsilon \in \Gamma$ $\displaystyle Ob(G) = X \subset 
\, dom(\varepsilon)$ and $\displaystyle dom(\varepsilon) = 
dom(\varepsilon)^{-1}$, 
\item[(ii)]  for any bounded set $K \subset Ob(G)$ there are  $1<A<B$  such that 
for any $\varepsilon \in \Gamma$, $\mid \varepsilon\mid \leq 1$: 
$$d^{-1}(\mid \varepsilon\mid) \, \cap \, \alpha^{-1}(K) \, \subset \,  
\delta_{\varepsilon} \, \left( d^{-1}(A) \cap \alpha^{-1}(K) \right)  
\, \subset  \, dom(\varepsilon^{-1}) \, \cap \, \alpha^{-1}(K)  \, \subset $$
\begin{equation}
\subset \,  
\delta_{\varepsilon} \, \left( d^{-1}(B) \cap \alpha^{-1}(K) \right)  \, \subset
\, \delta_{\varepsilon} \, \left( dom(\varepsilon) \cap \alpha^{-1}(K)\right)
\label{eax0}
\end{equation}
\item[(iii)]  for any bounded set $K \subset Ob(G)$ there are $R> 0$ and 
$\displaystyle \varepsilon_{0} \in (0,1]$ such that 
for any $\varepsilon \in \Gamma$, $\displaystyle \mid \varepsilon\mid \leq
\varepsilon_{0}$ and any $\displaystyle g,h \in \mid d^{-1}(R) \cap \,
\alpha^{-1}(K)$  we have:
\begin{equation}
dif(\delta_{\varepsilon} g, \delta_{\varepsilon} h) \, \in \,
dom(\varepsilon^{-1})
\label{edomdif}
\end{equation}
\end{enumerate}

\label{ddefor}
\end{definition}

\subsection{Dilation structures, definition}

By proposition \ref{enoughgen}, any metric space $(X,d)$ may be seen as the normed groupoid $(X \times X, d)$. Let us see what happens if we endow this trivial groupoid with dilations, according to definition \ref{ddefor}. For any $\varepsilon \in \Gamma$ we have a  dilation 
$$\displaystyle \delta_{\varepsilon} : dom(\varepsilon) \subset X^{2} \rightarrow im(\varepsilon) \subset X^{2}$$ 
which satisfies a number of axioms. Let us take them one by one. 

\paragraph{A1. for trivial groupoids.} For any 
$\varepsilon \in \Gamma$ $\displaystyle \alpha \delta_{\varepsilon} = \alpha$ is equivalent to 
the existence of a locally defined function $\displaystyle \delta^{x}_{\varepsilon}$, for any $x \in X$, such that 
 $$\delta_{\varepsilon}(y,x) = (\delta^{x}_{\varepsilon} y , x)$$  
The domain of definition of $\displaystyle \delta^{x}_{\varepsilon}$ is $dom(\varepsilon) \cap 
\left\{ (y,x) \mbox{ : } y \in X \right\}$ (similarily for the image). 
The fact that  $\displaystyle
\varepsilon \in \Gamma \mapsto \delta_{\varepsilon}$ is an action of $\Gamma$ on
$X^{2}$,  translates into: for any $\varepsilon, \mu \in \Gamma$ and $x \in X$ we have 
$\displaystyle \delta_{\varepsilon}^{x} \delta_{\mu}^{x} = \delta_{\varepsilon \mu}^{x}$, $\displaystyle \left(\delta_{\varepsilon}^{x} \right)^{-1} = \delta_{\varepsilon^{-1}}^{x}$ and 
$\delta_{e}^{x} = \, id$. 

\paragraph{A2. for trivial groupoids.} The objects of the trivial groupoid $\displaystyle X^{2}$ are of the form $(x,x)$ with $x \in X$. We use what we already know from A1 to deduce that the axiom A2 says:  for any $x \in X$ and any 
$\varepsilon \in \Gamma$ we have $\displaystyle \delta_{\varepsilon}^{x} x = x$. Moreover the transformation $\displaystyle (y,x) \mapsto (\delta_{\varepsilon}^{x} y, x)$ contracts the domain $dom(\varepsilon)$ to $\left\{ (x,x) \mbox{ : } x \in X\right\}$,  uniformly on  sets 
$\displaystyle A \subset X^{2}$ which are "bounded" in the sense: there is a $M> 0$ such that for any $(x,y) \in A$ we have $d(x,y) \leq M$. This  means that 
the net of functions $\displaystyle (x,y) \mapsto d (\delta_{\varepsilon}^{x} y, x)$ converges to the constant function $0$, uniformly with respect to $(x,y) \in A$, where $A$ is bounded in the sense explained before.

\paragraph{A3. for trivial groupoids.} A simple computation shows that a pair $\displaystyle (g,h) \in G \times_{\alpha} G$ has the form $g = (u,x)$, $h = (v,x)$. Then 
$$dif(\delta_{\varepsilon} g, \delta_{\varepsilon} h) = \left( \delta_{\varepsilon^{-1}}^{\delta^{x}_{\varepsilon} v} \, \delta_{\varepsilon}^{x} u  \right)$$
The axiom A3 says that for any $x \in X$ there is a function $\displaystyle d^{x}$, locally defined on pair of points $(u,v) \in X \times X$, such that 
$$ \lim_{\varepsilon \rightarrow 0} \frac{1}{\mid\varepsilon\mid} \, 
d \left( \delta_{\varepsilon^{-1}}^{\delta^{x}_{\varepsilon} v} \, \delta_{\varepsilon}^{x} u  \right) \, = \, d_{x}(u,v)$$
uniformly with respect to $d(x,u)$, $d(x,v)$. Moreover, if   $\displaystyle d_{x}(u,v) = 0$ then $u = v$. 

\paragraph{A4. for trivial groupoids.} Using A2 for trivial groupoids, this axiom says that 
the net $\displaystyle \delta_{\varepsilon^{-1}}^{\delta^{x}_{\varepsilon} v} \, \delta_{\varepsilon}^{x} u $ converges to a function $\displaystyle \Delta^{x}(v,u)$ uniformly with respect to $d(x,u)$, $d(x,v)$. 

The axiom A0 can also be detailed. All in all we see that trivial normed groupoids with dilations corespond to strong dilation structures. defined next.

\begin{definition}
Let $(X,d)$ be a complete metric space such that for any $x  \in X$ the 
closed ball $\bar{B}(x,3)$ is compact. A dilation structure $(X,d, \delta)$ 
over $(X,d)$ is the assignment to any $x \in X$  and $\varepsilon \in (0,+\infty)$ 
of a invertible homeomorphism, defined as: if 
$\displaystyle   \varepsilon \in (0, 1]$ then  $\displaystyle 
 \delta^{x}_{\varepsilon} : U(x)
\rightarrow V_{\varepsilon}(x)$, else 
$\displaystyle  \delta^{x}_{\varepsilon} : 
W_{\varepsilon}(x) \rightarrow U(x)$,  such that the following axioms are satisfied: 
\begin{enumerate}
\item[{\bf A0.}]  For any $x \in X$ the sets $ \displaystyle U(x), V_{\varepsilon}(x), 
W_{\varepsilon}(x)$ are open neighbourhoods of $x$.  There are  numbers  $1<A<B$ such that for any $x \in X$  and any 
$\varepsilon \in (0,1)$ we have 
  the following string of inclusions:
$$ B_{d}(x, \varepsilon) \subset \delta^{x}_{\varepsilon}  B_{d}(x, A) 
\subset V_{\varepsilon}(x) \subset 
W_{\varepsilon^{-1}}(x) \subset \delta_{\varepsilon}^{x}  B_{d}(x, B) $$
Moreover for any compact set $K \subset X$ there are $R=R(K) > 0$ and 
$\displaystyle \varepsilon_{0}= \varepsilon(K) \in (0,1)$  such that  
for all $\displaystyle u,v \in \bar{B}_{d}(x,R)$ and all 
$\displaystyle \varepsilon  \in (0,\varepsilon_{0})$,  we have 
$$\delta_{\varepsilon}^{x} v \in W_{\varepsilon^{-1}}( \delta^{x}_{\varepsilon}u) \ .$$

\item[{\bf A1.}]  We  have 
$\displaystyle  \delta^{x}_{\varepsilon} x = x $ for any point $x$. 
We also have $\displaystyle \delta^{x}_{1} = id$ for any $x \in X$. 
Let us define the topological space
$$ dom \, \delta = \left\{ (\varepsilon, x, y) \in (0,+\infty) \times X 
\times X \mbox{ :  if } \varepsilon \leq 1 \mbox{ then } y 
\in U(x) \,
\, , 
\right.$$ 
$$\left. \mbox{  else } y \in W_{\varepsilon}(x) \right\} $$ 
with the topology inherited from $(0,+\infty) \times X \times X$ endowed with
 the product topology. Consider also 
$\displaystyle Cl(dom \, \delta)$, 
the closure of 
$dom \, \delta$ in $\displaystyle [0,+\infty) \times X \times X$. The function $\displaystyle \delta : dom \, \delta 
\rightarrow  X$ defined by $\displaystyle \delta (\varepsilon,  x, y)  = 
\delta^{x}_{\varepsilon} y$ is continuous. Moreover, it can be continuously 
extended to the set $\displaystyle Cl(dom \, \delta)$ and we have 
$$\lim_{\varepsilon\rightarrow 0} \delta_{\varepsilon}^{x} y \, = \, x  $$

\item[{\bf A2.}] For any  $x, \in X$, $\displaystyle \varepsilon, \mu \in (0,+\infty)$
 and $\displaystyle u \in U(x)$   we have the equality: 
$$ \delta_{\varepsilon}^{x} \delta_{\mu}^{x} u  = \delta_{\varepsilon \mu}^{x} u $$ 
whenever one of the sides are well defined.

\item[{\bf A3.}]  For any $x$ there is a distance  function $\displaystyle (u,v) \mapsto d^{x}(u,v)$, defined for any $u,v$ in the closed ball (in distance d) $\displaystyle 
\bar{B}(x,A)$, such that 
$$\lim_{\varepsilon \rightarrow 0} \quad \sup  \left\{  \mid 
\frac{1}{\varepsilon} d(\delta^{x}_{\varepsilon} u, \delta^{x}_{\varepsilon} v) \ - \ d^{x}(u,v) \mid \mbox{ :  } u,v \in \bar{B}_{d}(x,A)\right\} \ =  \ 0$$
uniformly with respect to $x$ in compact set. 

\end{enumerate}

The dilation structure is strong if it satisfies the following supplementary condition: 

\begin{enumerate}
\item[{\bf A4.}] Let us define 
$\displaystyle \Delta^{x}_{\varepsilon}(u,v) =
\delta_{\varepsilon^{-1}}^{\delta^{x}_{\varepsilon} u} \delta^{x}_{\varepsilon} v$. 
Then we have the limit 
$$\lim_{\varepsilon \rightarrow 0}  \Delta^{x}_{\varepsilon}(u,v) =  \Delta^{x}(u, v)  $$
uniformly with respect to $x, u, v$ in compact set. 
\end{enumerate}
\label{defweakstrong}
\end{definition}
 
We shall use many times from now the words "sufficiently close". This deserves 
a definition. 

\begin{definition}
Let  $(X,d, \delta)$ be a strong dilation structure. We say that a property 
$$\displaystyle \mathcal{P}(x_{1},x_{2},
x_{3}, ...)$$ is true  for $\displaystyle x_{1}, x_{2}, x_{3},
...$ {\bf sufficiently close} if for any compact, non empty set $K \subset X$, there
is a positive constant $C(K)> 0$ such that $\displaystyle \mathcal{P}(x_{1},x_{2},
x_{3}, ...)$ is true for any $\displaystyle x_{1},x_{2},
x_{3}, ... \in K$ with $\displaystyle d(x_{i}, x_{j}) \leq C(K)$.
\end{definition}

\section{Examples of dilation structures}

\subsection{Snowflakes, nonstandard dilations in the plane}

\paragraph{Snowflake construction.} This is an adaptation of a standard construction for metric spaces with dilations: if $(X,d,\delta)$ is a dilation structure then $(X,d_{a}, \delta(a))$ is also a dilation structure, for any $a \in (0,1]$, where
$$d_{a}(x,y) \ = \ \left( d(x,y) \right)^{a} \ \ , \ \ \delta(a)_{\varepsilon}^{x}\ = \ \delta^{x}_{\varepsilon^{\frac{1}{a}}} \ .$$

In particular, if $X =  \mathbb{R}^{n}$ then  for any $a \in (0,1]$ we may take the distance and dilations 
$$ d_{a}(x,y) \ = \ \| x-y \|^{\alpha} \ \ \ ,  \ \ \delta^{x}_{\varepsilon} y \ =  \ x + \varepsilon^{\frac{1}{a}} (y - x) \ .$$

\paragraph{Nonstandard dilations.} In the plane $\displaystyle X = \mathbb{R}^{2}$,  endowed with the euclidean distance, we may consider another one-parameter group of linear transformations instead of the familiar homotheties. Indeed,  for any complex number  $z= 1+ i \theta$ let us  define  the dilations
$$\delta_{\varepsilon} x = \varepsilon^{z} x  \ .$$
Then  $(X,\delta, +, d)$ is a conical group, therefore the dilations
$$\delta^{x}_{\varepsilon} y = x + \delta_{\varepsilon} (y-x)  \ $$
together with the euclidean distance, form a dilation structure. 

Two such dilation structures, constructed respectively by using the numbers 
$1+ i \theta$ and $1+ i \theta'$,  are equivalent (see the section \ref{secequi} for the definition of equivalent dilation structures)) if and only if $\theta = \theta'$. 

Such dilation structures give examples of metric spaces with dilations which don't have the Radon-Nikodym property, see section \ref{radon}.

\subsection{Normed groups with dilations}

The following result is theorem 15 \cite{buligadil1}. 

\begin{theorem}
Let $(G, \delta, \| \cdot \|)$ be  a locally compact  normed local group with dilations. Then $(G, d, \delta)$ is 
a dilation structure, where $\delta$ are the dilations defined by (\ref{dilat}) and the distance $d$ is induced by the norm as in (\ref{dnormed}). 
\label{tgrd}
\end{theorem}

\paragraph{Proof.}
The axiom A0 is straightforward from definition \ref{defgwd},  axiom H0,  and because the dilation structure is left invariant, in the sense that the transport by left translations in 
$G$ preserves the dilations $\delta$. We also trivially have axioms A1 and A2 satisfied. 

For the axiom A3  remark that $\displaystyle d(\delta_{\varepsilon}^{x} u , \delta_{\varepsilon}^{x} v) = d(x \delta_{\varepsilon}(x^{-1}u),  x \delta_{\varepsilon}(x^{-1}u)) = d( \delta_{\varepsilon}(x^{-1}u) ,  \delta_{\varepsilon}(x^{-1}v))$. 
Let us denote $\displaystyle U= x^{-1}u$, $\displaystyle V= x^{-1}v$ and for $\varepsilon > 0$ let 
$$\displaystyle \beta_{\varepsilon}(u, v) =  \delta_{\varepsilon}^{-1}
\left((\delta_{\varepsilon}u) (\delta_{\varepsilon}v ) \right) .$$
 We have then: 
$$\frac{1}{\varepsilon} d(\delta_{\varepsilon}^{x} u , \delta_{\varepsilon}^{x} v) = 
\frac{1}{\varepsilon} \|\delta_{\varepsilon} \beta_{\varepsilon} \left(  \delta_{\varepsilon}^{-1}
\left( ( \delta_{\varepsilon}V)^{-1}\right) , U\right)\| \quad . $$

Define the function 
$$d^{x}(u,v) = \| \beta( V^{-1} , U) \|^{N} . $$
From definition \ref{defgwd} axioms H1, H2, and from definition \ref{dnco} (d), we obtain that axiom A3 
is satisfied. 

For the axiom A4 we have to compute: 
$$\Delta^{x}(u,v) = \delta_{\varepsilon^{-1}}^{\delta^{x}_{\varepsilon} u} \delta^{x}_{\varepsilon} v = 
\left(\delta^{x}_{\varepsilon} u \right) \left( \delta_{\varepsilon}\right)^{-1} \left( \left(\delta^{x}_{\varepsilon} u \right)^{-1} \left(\delta^{x}_{\varepsilon} v \right)\right)  = $$
$$= \left(x\delta_{\varepsilon} U \right)  \beta_{\varepsilon}\left(  \delta_{\varepsilon}^{-1}
\left( ( \delta_{\varepsilon}V)^{-1}\right) , U\right) \rightarrow x \beta\left( V^{-1}, U \right)  $$
as $\varepsilon \rightarrow 0$. Therefore axiom A4 is satisfied.\quad  \quad $\square$

\subsection{Riemannian manifolds} 
\label{riemann}

The following interesting   quotation from Gromov book \cite{gromovbook}, pages 85-86, motivates some of the  ideas underlying  dilation structures, 
especially in the very particular case of a riemannian manifold: 

``{\bf 3.15. 
Proposition:} {\it Let $(V, g)$ be a Riemannian manifold with $g$ continuous.
For each $v \in V$ the spaces $(V, \lambda d , v)$ Lipschitz converge as  $\lambda \rightarrow \infty$ to the 
tangent space $(T_{v}V, 0)$ with its Euclidean metric $g_{v}$.} 

$\mathbf{Proof_{+}:}$ {\it Start with a $C^{1}$ map $\mathbb{(R}^{n}, 0) \rightarrow  (V, v)$ whose differential is 
isometric at 0. The $\lambda$-scalings of this provide almost isometries between large balls in $\mathbb{R}^{n}$
 and those in $\lambda V$ for  $\lambda \rightarrow \infty$.} 
 {\bf Remark:} {\it In fact we can define Riemannian manifolds as 
 locally compact path metric spaces that satisfy the conclusion of Proposition 3.15.}``

The problem of domains and codomains left aside, any chart of a  Riemannian manifold induces locally a dilation 
structure on the manifold. Indeed, take $(M,d)$ to be a $n$-dimensional  Riemannian manifold with $d$ the distance on 
$M$ induced by the Riemannian structure. Consider a diffeomorphism $\phi$ of an open set $U \subset M$ onto 
$\displaystyle V \subset \mathbb{R}^{n}$  and transport the dilations from $V$ to $U$ (equivalently, transport the distance $d$ from 
$U$ to $V$). There is only one thing to check in order to see that we got a dilation structure: the axiom A3, expressing 
the compatibility of the distance $d$ with the dilations. But this is just a metric way to express the distance 
on the tangent space of $M$ at $x$ as a limit of rescaled distances (see  Gromov Proposition 3.15, \cite{gromovbook},
p. 85-86). Denoting by $\displaystyle g_{x}$ the metric tensor at $x \in U$, we have:  
$$\left[ d^{x}(u,v) \right]^{2}  = $$ 
$$= g_{x}\left( \frac{d}{d \, \varepsilon}_{|_{\varepsilon = 0}}\phi^{-1}\left(\phi(x) + \varepsilon (\phi(u) -
\phi(x))\right) ,  \frac{d}{d \, \varepsilon}_{|_{\varepsilon = 0}}\phi^{-1}\left(\phi(x) + \varepsilon (\phi(v) -
\phi(x))\right) \right)$$

A different example of a dilation structure on a riemannian manifold comes from the setting of proposition \ref{pcurvature}, section \ref{securv}. 

Let $M$ be a smooth enough riemannian manifold   and $\exp$ be the geodesic exponential. 
To any point $x \in M$ and any vector $\displaystyle v \in T_{x} M$ the point   $\displaystyle \exp_{x}(v) \in M$ is located on the geodesic passing thru $x$ 
and tangent to $v$; if we parameterize this geodesic with respect to length, such that the tangent at $x$ is parallel and 
has the same direction as $v$, then $\displaystyle \exp_{x}(v) \in M$ has the coordinate equal with the length of $v$ with 
respect to the norm on $\displaystyle T_{x} M$.  We define implicitly  the dilation based at $x$, of coefficient 
$\varepsilon > 0$ by the relation: 
 $$\delta^{x}_{\varepsilon} \exp_{x} (u) \, = \, \exp_{x} \left( \varepsilon u \right) \quad . $$

\begin{proposition}
The above example is a strong dilation structure. 
\label{priemann}
\end{proposition}

\paragraph{Proof.} This field of dilations satisfies trivially A0, A1, A2, only A3 and A4 are left to be checked.  

Proposition \ref{pcurvature} provides a proof for A3. For the proof of A4 see the section 
\ref{srsmooth}, where normal and adapted frames are defined for sub-riemannian manifolds. In particular the same construction works for riemannian manifolds, where one can attach normal frames to geodesic coordinate systems. \quad $\square$

\section{Length dilation structures}
\label{seclds}

Consider $(X,d)$ a complete, locally compact metric space, and a triple 
 $(X,d,\delta)$  which satisfies  A0, A1, A2. Denote by $\displaystyle Lip([0,1],X,d)$ the space of $d$-Lipschitz curves $c:[0,1] \rightarrow X$. Let also $\displaystyle 
l_{d}$ denote the length functional associated to the distance $d$.

\begin{definition}
For any $\varepsilon \in (0,1)$ and $x \in X$ we define the length functional at scale $\varepsilon$, relative to $x$, to be 
$$l_{\varepsilon}(x,c) \ = \ l^{x}_{\varepsilon}(c) \ = \ \frac{1}{\varepsilon}
\, 
l_{d}(\delta^{x}_{\varepsilon} c)  $$
The domain of definition of the functional 
$\displaystyle l_{\varepsilon}$ is the space: 
 $$\mathcal{L}_{\varepsilon}(X,d,  \delta) \ = \ \left\{(x ,c) \in X  
 \times \mathcal{C}([0,1],X) 
\mbox{ : } c: [0,1] \in U(x) \, \, , \,  \right.$$ 
$$\left. \delta^{x}_{\varepsilon}c \mbox{ is } 
d-Lip \mbox{ and } Lip(\delta^{x}_{\varepsilon}c)\,  \leq \, 2  \, 
l_{d}(\delta^{x}_{\varepsilon}c) 
\right\} $$

\label{thespaceleps}
\end{definition}

The last condition from the definition of $\displaystyle 
\mathcal{L}_{\varepsilon}(X,d,  \delta)$ is a selection of parameterization 
of the path $c([0,1])$. Indeed, by the reparameterization theorem, if 
$\displaystyle \delta^{x}_{\varepsilon} c :[0,1] \rightarrow (X,d)$ is a 
$d$-Lipschitz curve of length $\displaystyle L = l_{d}(\delta^{x}_{\varepsilon}c)$ 
then $\displaystyle \delta^{x}_{\varepsilon}c([0,1])$ can be reparameterized by length, that is there exists a 
increasing  function $\phi:[0,L] \rightarrow [0,1]$ such that $\displaystyle 
c'= \delta^{x}_{\varepsilon} c\circ \phi$ is 
a $d$-Lipschitz curve with $Lip(c') \leq 1$. But we can use a second affine
reparameterization which sends $[0,L]$ back to $[0,1]$ and we get a Lipschitz
curve $c"$ with $c"([0,1]) = c'([0,1])$ and $\displaystyle Lip(c") \leq 2 
l_{d}(c)$.

In the definition of dilation structures we  use uniform convergence of distances. Here we need a notion of convergence for length functionals.  This is the  Gamma-convergence, a notion used many time in calculus of variations. A good reference is the book  \cite{dalmaso}. 

\begin{definition}
Let $Z$ be a metric space with distance function $D$ and $\displaystyle \left(
l_{\varepsilon}\right)_{\varepsilon > 0}$ be a family of functionals $\displaystyle 
l_{\varepsilon}: Z_{\varepsilon} \subset Z \rightarrow [0,+\infty]$. Then 
$\displaystyle l_{\varepsilon}$ Gamma-converges to the functional 
$\displaystyle l: Z_{0} \subset Z \rightarrow [0,+\infty]$ if: 
\begin{enumerate}
\item[(a)] ({\bf liminf inequality}) for any function $\displaystyle \varepsilon \in
(0,\infty)  \mapsto 
x_{\varepsilon} \in Z_{\varepsilon}$ such that $\displaystyle \lim_{\varepsilon
\rightarrow 0} x_{\varepsilon} \, = \, x_{0} \in Z_{0}$ we have 
$$l(x_{0}) \, \leq \, \liminf_{\varepsilon \rightarrow 0}
l_{\varepsilon}(x_{\varepsilon})$$
\item[(b)] ({\bf existence of a recovery sequence}) For any $\displaystyle x_{0} \in Z_{0}$ 
and for any sequence $\displaystyle \left( \varepsilon_{n} \right)_{n \in \mathbb{N}}$
such that $\displaystyle \lim_{n \rightarrow \infty} \varepsilon_{n} \, = \, 0$ there
is a sequence $\displaystyle \left( x_{n} \right)_{n \in \mathbb{N}}$ with
$\displaystyle x_{n} \in Z_{\varepsilon_{n}}$ for any $n \in \mathbb{N}$, such that 
$$l(x_{0}) \, = \, \lim_{n \rightarrow \infty}
l_{\varepsilon_{n}}(x_{n})$$
\end{enumerate}
\label{defgammamet}
\end{definition}

For our needs we shall take  $Z$ to be the space
 $X \times \mathcal{C}([0,1],X)$ endowed with the distance 
 $$ D((x,c), (x', c')) \ = \ \max\left\{ d(x, x') \, , \, \sup \left\{ d(c(t),
 c'(t)) \mbox{ : } t \in [0,1] \right\} \right\} $$

Let  $\displaystyle \mathcal{L}(X,d,\delta)$be the class of all 
$\displaystyle (x ,c) \in X  \times \mathcal{C}([0,1],X)$ which appear as  limits  
$\displaystyle (x_{n}, c_{n}) \rightarrow (x,c)$, with $\displaystyle 
(x_{n}, c_{n}) \in \mathcal{L}_{\varepsilon_{n}}(X,d,  \delta)$, 
the family $\displaystyle (c_{n})_{n}$ is $d$-equicontinuous and $\displaystyle 
\varepsilon_{n} \rightarrow 0$ as $n \rightarrow \infty$.

\begin{definition}
A triple $(X,d,\delta)$ is a  length dilation structure if $(X,d)$ is a
complete, locally compact metric space such that A0, A1, 
A2,  are satisfied, together with  the following axioms: 
\begin{enumerate}
\item[\bf{A3L.}] there is a functional $\displaystyle l : \mathcal{L}(X,d,  \delta) \rightarrow
[0,+\infty]$ such that for any $\varepsilon_{n} \rightarrow 0$ as $n \rightarrow
\infty$ the sequence of functionals 
$\displaystyle l_{\varepsilon_{n}}$ Gamma-converges to the functional $l$. 
\item[\bf{A4+}] Let us define 
$\displaystyle \Delta^{x}_{\varepsilon}(u,v) =
\delta_{\varepsilon^{-1}}^{\delta^{x}_{\varepsilon} u} \delta^{x}_{\varepsilon}
v$ and $\displaystyle \Sigma^{x}_{\varepsilon}(u,v) = 
\delta_{\varepsilon^{-1}}^{x} \delta_{\varepsilon}^{\delta^{x}_{\varepsilon} u}
v$. 
Then we have the limits 
$$\lim_{\varepsilon \rightarrow 0}  \Delta^{x}_{\varepsilon}(u,v) =  \Delta^{x}(u, v)  $$
$$\lim_{\varepsilon \rightarrow 0}  \Sigma^{x}_{\varepsilon}(u,v) =  \Sigma^{x}(u, v)  $$
uniformly with respect to $x, u, v$ in compact set. 
\end{enumerate} 
\label{deflds}
\end{definition}

\begin{remark}
For strong dilation structures the axioms A0 - A4 imply A4+, cf. corollary 9  \cite{buligadil1}. The
transformations $\displaystyle \Sigma^{x}_{\varepsilon}(u, \cdot)$ have the
interpretation of approximate left translations in the tangent space of $(X,d)$
at $x$.  
\end{remark}

For any $\varepsilon \in (0,1)$ and any $x \in X$ the length functional 
$\displaystyle l^{x}_{\varepsilon}$ induces a distance on $U(x)$: 
$$ \mathring{d}^{x}_{\varepsilon}(u,v) \ = \ \inf\left\{ l^{x}_{\varepsilon}(c)
\mbox{ : } (x,c) \in \mathcal{L}_{\varepsilon}(X,d,  \delta) \, , \, c(0) = u \,
, \, c(1) = v \right\} $$
In the same way the length functional $l$ from A3L induces a distance
$\displaystyle \mathring{d}^{x}$ on $U(x)$.

Gamma-convergence implies that 
\begin{equation}
\mathring{d}^{x}(u,v) \, \geq \, \limsup_{\varepsilon \rightarrow 0} \mathring{d}^{x}_{\varepsilon}(u,v)
\label{dsup}
\end{equation}
but We don't believe that, at this level of generality, we could have equality without supplementary hypotheses. This means that, probably,  there exist length dilation structures which are not  strong dilation structures.

\section{Properties of dilation structures}
\label{sprop}

\subsection{Metric profiles associated with dilation structures}

In this subsection we shall look at dilation structures from the metric point of view, by  using Gromov-Hausdorff distance and metric profiles.

Let us denote by $(\delta, \varepsilon)$ the distance on 
$$\bar{B}_{d^{x}}(x,1) \ = \ \left\{ y \in X \mbox{ : } d^{x}(x,y) \leq 1 \right\}$$ given by
$$(\delta, \varepsilon)(u,v) \ = \ \frac{1}{\varepsilon} d(\delta^{x}_{\varepsilon} u , \delta^{x}_{\varepsilon} v) \ .$$

The following theorem is a generalization of the   Mitchell \cite{mit} theorem 1, concerning  sub-riemannian geometry.

\begin{theorem}
Let $(X,d,\delta)$ be a dilation structure. 

For all $u,v \in X$ such that $\displaystyle d(x,u)\leq 1$ and $\displaystyle d(x,v) \leq 1$  and all $\mu \in (0,A)$ we have: 
$$d^{x}(u,v) \ = \ \frac{1}{\mu} d^{x}(\delta_{\mu}^{x} u , \delta^{x}_{\mu} v) \ .$$
Therefore $\displaystyle (U(x), d^{x}, x)$ is a metric cone. 

Moreover, if the dilation structure is strong  then  the curve $\displaystyle \varepsilon> 0 \mapsto \mathbb{P}^{x}(\varepsilon) \ = \ [\bar{B}_{d^{x}}(x,1), (\delta, \varepsilon), x]$ is an abstract  metric profile.

Finally,  we have  the following limit: 
$$\lim_{\varepsilon \rightarrow 0} \ \frac{1}{\varepsilon} \sup \left\{  \mid d(u,v) - d^{x}(u,v) \mid \mbox{ : } d(x,u) \leq \varepsilon \ , \ d(x,v) \leq \varepsilon \right\} \ = \ 0 \ .$$
therefore if A4 holds then  $(X,d)$ admits a metric tangent space 
in $x$. 
\label{thcone}
\end{theorem}

\paragraph{Proof.}
For fixed  $\mu \in (0,1)$ and variable $\varepsilon \in (0,1)$  we have: 
$$\mid \frac{1}{\mu} \frac{1}{\varepsilon} d(\delta_{\varepsilon}^{x}\delta^{x}_{\mu} u, \delta_{\varepsilon}^{x}\delta^{x}_{\mu} v) \ - \ d^{x}(u,v) \mid \ = \  \mid \frac{1}{\varepsilon \mu} d(\delta_{\varepsilon \mu}^{x}u, \delta_{\varepsilon \mu}^{x} v) \ - \ d^{x}(u,v) \mid \ .$$
We pass pass to the limit with $\varepsilon \rightarrow 0$ and we obtain the cone property of the distance $\displaystyle d^{x}$. 

Let us  prove that $\mathbb{P}^{x}$ is an abstract metric profile. For this we have to compare two pointed metric spaces, namely 
$\displaystyle  \left((\delta^{x}, \varepsilon \mu) , \bar{B}_{d^{x}}(x,1), x\right)$ and 
$\displaystyle  \left(\frac{1}{\mu}(\delta^{x}, \varepsilon), \bar{B}_{\frac{1}{\mu}(\delta^{x}, \varepsilon)}(x,1), x \right)$. Let $u \in X$ such that 
$$\frac{1}{\mu}(\delta^{x}, \varepsilon)(x, u) \leq 1 \ .$$
From the axioms of dilation structures and the cone property we obtain the estimate:
$$\frac{1}{\varepsilon} d^{x}( x, \delta^{x}_{\varepsilon} u) \ \leq \ (\mathcal{O}(\varepsilon) + 1) \mu$$
therefore 
$\displaystyle d^{x}(x,u)  \  \leq \ (\mathcal{O}(\varepsilon) + 1) \mu$. 
It follows that for any $\displaystyle u \in \bar{B}_{\frac{1}{\mu}(\delta^{x}, \varepsilon)}(x,1)$ we can choose a point $\displaystyle w(u) \in \bar{B}_{d^{x}}(x,1)$ such that
$$\frac{1}{\mu} d^{x}(u, \delta^{x}_{\mu} w(u)) \ = \ \mathcal{O}(\varepsilon) \ .$$
Then, by using twice A3, we obtain
$$\mid \frac{1}{\mu}(\delta^{x}, \varepsilon) (u_{1}, u_{2}) \ - \ (\delta^{x}, \varepsilon \mu) (w(u_{1}), w(u_{2}) ) \mid \ = $$ 
$$= \ \mid \frac{1}{\varepsilon \mu} d( \delta^{x}_{\varepsilon} u_{1}, \delta^{x}_{\varepsilon} u_{2}) \ - \ \frac{1}{\varepsilon \mu} d(\delta^{x}_{\varepsilon}\delta^{x}_{ \mu} w(u_{1}), \delta^{x}_{\varepsilon}\delta^{x}_{ \mu} w(u_{2}))\mid \ \leq \ $$
$$\leq  \ \frac{1}{\mu} \mathcal{O}(\varepsilon) \ + \ 
\frac{1}{\mu} \mid d^{x} (u_{1}, u_{2}) \ - \ d^{x}(\delta^{x}_{\mu} w(u_{1}) , \delta^{x}_{\mu} w(u_{2})) \mid \ \leq $$  
$$\leq  \ \frac{1}{\mu} \mathcal{O}(\varepsilon) \ + \ \frac{1}{\mu} 
d^{x}(u_{1}, \delta^{x}_{\mu} w(u_{1})) \  + \ \frac{1}{\mu} 
d^{x}(u_{1}, \delta^{x}_{\mu} w(u_{2})) \ \leq  $$ 
$$ \  \leq \   \frac{1}{\mu} \mathcal{O}(\varepsilon) + \mathcal{O}(\varepsilon) \  .$$
This shows that the property (b) of an abstract  metric profile is satisfied. For the property (a) of an abstract metric profile we do the following.  By A0,  for $\varepsilon \in (0,1)$ and $\displaystyle u,v \in  \bar{B}_{d}(x,\varepsilon)$ there exist $\displaystyle U,V  \in \bar{B}_{d}(x,A)$ such that $$u = \delta^{x}_{\varepsilon} U , v = \delta^{x}_{\varepsilon} V . $$
By the cone property we have 
$$\frac{1}{\varepsilon} \mid d(u,v) - d^{x}(u,v) \mid = \mid \frac{1}{\varepsilon} d(\delta^{x}_{\varepsilon} U,\delta^{x}_{\varepsilon} V) - d^{x}(U,V) \mid . $$
By A2 we have 
$$  \mid \frac{1}{\varepsilon} d(\delta^{x}_{\varepsilon} U,\delta^{x}_{\varepsilon} V) - d^{x}(U,V) \mid \leq \mathcal{O}(\varepsilon) .   \quad \square$$

\subsection{The tangent bundle of a dilation structure}
\label{induced}

The following proposition contains the main relations between the approximate difference,
sum and inverse operations. In \cite{buligadil1} I explained these relations as appearing from the equivalent formalism using binary decorated trees.

\begin{proposition}
Let $\displaystyle (X,\circ_{\varepsilon})_{\varepsilon \in \Gamma}$ be
a $\Gamma$-irq. Then we have the relations: 
\begin{enumerate}
\item[(a)] $\displaystyle \Delta^{x}_{\varepsilon}(u,
\Sigma^{x}_{\varepsilon}(u,v)) \, = \, v$  (difference is the inverse of sum) 
\item[(b)] $\displaystyle \Sigma^{x}_{\varepsilon}(u,
\Delta^{x_{\varepsilon}}(u,v)) \, = \, v$  (sum is the inverse of difference) 
\item[(c)]  $\displaystyle \Delta^{x}_{\varepsilon}(u, v) \, = \, \Sigma^{x \circ_{\varepsilon} u}_{\varepsilon} 
(inv_{\varepsilon}^{x} u , v)$ (difference  approximately equals  the sum of the inverse) 
\item[(d)]  $\displaystyle inv_{\varepsilon}^{x\circ u} \, inv_{\varepsilon}^{x}
\, u  \, = \, u $ (inverse operation is approximatively an involution) 
\item[(e)] $\displaystyle \Sigma^{x}_{\varepsilon}(u, \Sigma^{x\circ_{\varepsilon}u}_{\varepsilon}
(v , w)) \, = \, \Sigma^{x}_{\varepsilon}(\Sigma^{x}_{\varepsilon}(u,v), w) $ (approximate associativity of the sum)  
\item[(f)] $\displaystyle  inv^{x}_{\varepsilon} \, u \, = \,  \Delta^{x}_{\varepsilon}( u , x)$
\item[(g)]  $\displaystyle  \Sigma^{x}_{\varepsilon} (x, u) \, = \,  u $
(neutral element at right).
\end{enumerate}
\label{pplay}
\end{proposition}

The next theorem is the generalization of proposition \ref{here3.4}.  It is the main structure theorem fot the tangent bundle associated to a dilation structure, see theorems 7, 8,
10 in \cite{buligadil1}. 

\begin{theorem}
Let $(X,d,\delta)$ be a strong dilation structure. Then for any $x \in X$ 
$\displaystyle (U(x), \Sigma^{x}, \delta^{x})$ is a conical group. Moreover, left translations of this group are 
$\displaystyle d^{x}$ isometries. 
\label{tgene}
\end{theorem}

\paragraph{Proof.} 
(I.) Let us define  the "infinitesimal translations" 
$$\displaystyle L^{x}_{u}(v) =  \lim_{\varepsilon \rightarrow 0}  \Delta^{x}_{\varepsilon}(u,v)$$ and prove that they are  $\displaystyle d^{x}$ isometries. 

From theorem \ref{thcone} we get the limit, as $\varepsilon \rightarrow 0$: 
\begin{equation}
\sup \left\{ \frac{1}{\varepsilon} \mid d(u, v) \ - \ d^{x}(u,v) \mid \mbox{ : } d(x,u) \leq \frac{3}{2}\varepsilon, \ 
d(x,v) \leq \frac{3}{2}\varepsilon \right\} \rightarrow 0
\label{estiminf}
\end{equation}

 For any  $\varepsilon > 0$ sufficiently small the points 
$\displaystyle x, \delta^{x}_{\varepsilon}u, \delta^{x}_{\varepsilon} v , \delta^{x}_{\varepsilon} w$ are close one to another. Indeed, we have 
$\displaystyle d(\delta^{\varepsilon}_{x} u, \delta^{\varepsilon}_{x} v) = \varepsilon( d^{x}(u,v) + \mathcal{O}(\varepsilon))$. Therefore, if we choose $u,v,w$ such that $\displaystyle d^{x}(u,v)< 1, d^{x}(u,w)< 1$, then there is $\eta>0$ such that 
for all $\varepsilon\in (0,\eta)$ we have 
$$d(\delta^{\varepsilon}_{x} u, \delta^{\varepsilon}_{x} v) \leq \frac{3}{2}\varepsilon \quad , \quad  d(\delta^{\varepsilon}_{x} u, 
\delta^{\varepsilon}_{x} v) \leq \frac{3}{2}\varepsilon \quad . $$

We use  (\eqref{estiminf})  for the basepoint $\displaystyle \delta^{x}_{\varepsilon} u$  to get, as $\varepsilon \rightarrow 0$
$$\frac{1}{\varepsilon} \mid d(\delta^{x}_{\varepsilon} v , \delta^{x}_{\varepsilon} w) \ - \ 
d^{\delta^{x}_{\varepsilon} u}(\delta^{x}_{\varepsilon} v , \delta^{x}_{\varepsilon} w) \mid \rightarrow 0$$
From the cone property of the distance $\displaystyle d^{\delta^{x}_{\varepsilon} u}$  
\begin{equation}
\mid \frac{1}{\varepsilon} d(\delta^{x}_{\varepsilon} v , \delta^{x}_{\varepsilon} w) \ - \ 
d^{\delta^{x}_{\varepsilon} u}\left( \delta_{\varepsilon^{-1}}^{\delta^{x}_{\varepsilon} u} \delta^{x}_{\varepsilon} v ,  \delta_{\varepsilon^{-1}}^{\delta^{x}_{\varepsilon} u} \delta^{x}_{\varepsilon} w 
\right) \mid \rightarrow 0
\label{esti2}
\end{equation}
as $\varepsilon \rightarrow 0$. By the axioms A1, A3, the function $\displaystyle (x,u,v) \mapsto d^{x}(u,v)$ is uniformly continuous on compact sets, being an uniform limit of continuous functions.  We prove the fact that   the "infinitesimal translations are are  $\displaystyle d^{x}$ isometries 
by passing  to the limit in the LHS of (\eqref{esti2})  and by using 
this uniform continuity.   

(II.)  If for any $x$ the distance $\displaystyle d^{x}$ is non degenerate then there exists $C>0$ such that: 
  for any $x$ and $u$ with $d(x,u) \leq C$ there exists a $\displaystyle d^{x}$ 
isometry $\displaystyle \Sigma^{x}(u, \cdot)$ obtained as the limit:
$$ \lim_{\varepsilon \rightarrow 0} \Sigma^{x}_{\varepsilon}(u,v) = \Sigma^{x}(u, v)$$ 
uniformly with respect to $x, u, v$ in compact set. 
 
Indeed, from the step (I.)  we know that $\displaystyle \Delta^{x}(u,\cdot)$ is  
a $\displaystyle d^{x}$ isometry. If $\displaystyle d^{x}$ is non degenerate then  $\displaystyle \Delta^{x}(u,\cdot)$ is invertible. Let  $\displaystyle \Sigma^{x}(u,\cdot)$ be the inverse. 

From proposition \ref{pplay} we know that $\displaystyle \Sigma^{x}_{\varepsilon}(u, \cdot)$ is the inverse of $\displaystyle \Delta^{x}_{\varepsilon}(u, \cdot)$. Therefore 
$$d^{x}( \Sigma^{x}_{\varepsilon}(u, w), \Sigma^{x}(u, w) ) = d^{x}( \Delta^{x}(u,  \Sigma^{x}_{\varepsilon}(u, w)) , w) = $$
$$ = d^{x}( \Delta^{x}(u,  \Sigma^{x}_{\varepsilon}(u, w)) ,  \Delta^{x}_{\varepsilon}(u,  \Sigma^{x}_{\varepsilon}(u, w)) . $$
From the uniformity of convergence in step (I.) and the uniformity assumptions in axioms of dilation structures, the conclusion follows. 

(III.) We start by proving that $\displaystyle (U(x), \Sigma^{x})$ is a local uniform group. The uniformities are induced by the distance $d$. 

Indeed,  according to  proposition \ref{pplay}, we can pass to the limit with $\varepsilon \rightarrow 0$ and 
define: 
$$inv^{x}(u) = \lim_{\varepsilon \rightarrow 0} \Delta^{x}_{\varepsilon}(u,x) = \Delta^{x}(u,x) . $$
From relation (d) , proposition (\ref{pplay}) we get (after passing to the limit with $\varepsilon \rightarrow 0$) 
$$inv^{x}(inv^{x}(u)) = u .$$
We shall see that $\displaystyle inv^{x}(u)$ is the inverse of $u$. 
Relation (c), proposition  (\ref{pplay}) gives: 
\begin{equation}
\Delta^{x}(u,v) = \Sigma^{x} (inv^{x}(u), v)
\label{po}
\end{equation}
therefore  relations (a), (b) from proposition \ref{pplay} give 
\begin{equation}
\Sigma^{x}(inv^{x}(u), \Sigma^{x}(u,v)) = v , 
\label{pb}
\end{equation}
\begin{equation}
\Sigma^{x}(u, \Sigma^{x}(u,v)) = v .
\label{pa}
\end{equation}
Relation (e) from proposition \ref{pplay} gives 
\begin{equation}
\Sigma^{x}(u, \Sigma^{x}(v,w)) = \Sigma^{x}(\Sigma^{x}(u,v), w) 
\label{pc}
\end{equation}
which shows that $\displaystyle \Sigma^{x}$ is an associative operation. From (\ref{pa}), (\ref{pb}) we 
obtain that for any $u,v$ 
\begin{equation}
\Sigma^{x}(\Sigma^{x}(inv^{x}(u), u),v) = v , 
\label{pd}
\end{equation}
\begin{equation}
\Sigma^{x}(\Sigma^{x}(u, inv^{x}(u)), v) = v .
\label{pe}
\end{equation}
Remark that for any $x$, $v$ and $\varepsilon \in (0,1)$ we have $\displaystyle 
\Sigma^{x}(x,v) =  v$. 
    Therefore $x$ is a neutral element at left for the operation $\displaystyle \Sigma^{x}$. 
From the definition of $\displaystyle inv^{x}$, relation (\ref{po}) and the fact that $\displaystyle inv^{x}$ 
is equal to its inverse, we get that $x$ is an inverse at right too: for any $x$, $v$ we have
$$\Sigma^{x}(v,x) = v . $$
Replace now $v$ by $x$ in relations (\ref{pd}), (\ref{pe}) and prove that indeed $\displaystyle inv^{x}(u)$ 
is  the inverse of $u$.

We also have to prove that $\displaystyle (U(x), \Sigma^{x})$ admits $\displaystyle \delta^{x}$ as dilations.In this reasoning we need the axiom A2 in strong form. 

Namely we have to prove that for any $\mu \in (0,1)$ we have 
$$\delta^{x}_{\mu} \Sigma^{x}(u,v) = \Sigma^{x}(\delta^{x}_{\mu} u, \delta^{x}_{\mu} v) . $$
For this is sufficient to notice that 
$$\delta^{x}_{\mu} \Delta^{x}_{\varepsilon \mu} (u,v)  = \Delta^{x}_{\varepsilon} (\delta^{x}_{\mu} u, \delta^{x}_{\mu} v) $$
and pass to the limit as $\varepsilon \rightarrow 0$.   \quad  $\square$

\begin{definition}
The (local) conical group $\displaystyle (U(x), \Sigma^{x}, \delta^{x})$ is the tangent space 
of $(X,d, \delta)$ at $x$ (in the sense of dilation structures). We  denote it by  
$\displaystyle T_{x} (X, d, \delta) =  (U(x), \Sigma^{x}, \delta^{x})$, or by $\displaystyle T_{x} X$ if $(d,\delta)$ are clear from the context. 
\label{deftdil}
\end{definition}

\paragraph{Compatibility of topologies.} The axiom A3 implies that for any $x \in X$ the function $\displaystyle d^{x}$ is continuous, therefore  open sets with respect to $\displaystyle d^{x}$ are open with respect to $d$.

 If  $(X,d)$ is separable and $\displaystyle d^{x}$ is non degenerate then the uniformities induced by  $d$ and  $\displaystyle d^{x}$ are the same. Indeed, let 
 $\displaystyle \left\{u_{n} \mbox{ : } n \in \mathbb{N}\right\}$ 
 be a dense set in $U(x)$, with $\displaystyle x_{0}=x$. 
 We can embed $\displaystyle (U(x),  (\delta^{x}, \varepsilon))$  isometrically in the separable Banach space 
 $\displaystyle l^{\infty}$, for any $\varepsilon \in (0,1)$, by the function 
 $$\phi_{\varepsilon}(u) = \left( \frac{1}{\varepsilon} d(\delta^{x}_{\varepsilon}u,  \delta^{x}_{\varepsilon}x_{n}) - \frac{1}{\varepsilon} d(\delta^{x}_{\varepsilon}x,  \delta^{x}_{\varepsilon}x_{n})\right)_{n}  . $$
 A reformulation of the first part  of theorem \ref{thcone} is that on compact sets $\displaystyle \phi_{\varepsilon}$ uniformly converges to the isometric embedding of $\displaystyle (U(x), d^{x})$ 
 $$\phi(u) = \left(  d^{x}(u,  x_{n}) - d^{x}(x, x_{n})\right)_{n}  . $$
Remark that the uniformity induced by $(\delta,\varepsilon)$ is the same as the uniformity 
induced by $d$, and that it is the same induced from the uniformity on 
$\displaystyle l^{\infty}$ by 
the embedding $\displaystyle \phi_{\varepsilon}$. We proved that the uniformities induced by 
 $d$ and  $\displaystyle d^{x}$ are the same.

From previous considerations we deduce the following characterization of  tangent
spaces associated to a dilation structure. 

\begin{corollary}
Let $(X,d,\delta)$ be a strong dilation structure with group $\Gamma = (0,+\infty)$. 
Then for any $x \in X$ the local group   $\displaystyle (U(x), \Sigma^{x})$ is locally a simply connected Lie group  whose Lie algebra admits a positive graduation (a homogeneous group).
 \label{cortang}
\end{corollary}

\paragraph{Proof.}
Indeed, from previous considerations,   $\displaystyle (U(x), \Sigma^{x})$ is a locally compact
group  which admits 
$\displaystyle \delta^{x}$ as a contractive automorphism group (from theorem \ref{tgene}). Instead of Siebert proposition \ref{psiebert}, we need now a version for local groups. Fortunately, theorem 1.1 \cite{recent} states that a locally compact, locally connected, contractible (with  Siebert' wording) group is locally isomorphic to a contractive Lie group.   \quad $\square$

 Straightforward modifications in the proof of the previous 
theorem allow us to extend some results to  length dilation structures.

\begin{theorem}
Let $(X,d,\delta)$ be a  a length dilation
structure. Then: 
 \begin{enumerate}
 \item[(a)]  $\displaystyle \Sigma^{x}$ is a local group operation on $U(x)$,
 with $x$ as neutral element and $\displaystyle \, inv^{x}$ as the inverse element
 function; for any $\varepsilon \in (0,1]$ the 
 dilation $\displaystyle \delta^{x}_{\varepsilon}$ is an automorphism with 
 respect to the group operation; 
  \item[(b)]the  length functional $\displaystyle l^{x} = l(x, \cdot)$ is invariant with
  respect to left translations $\displaystyle \Sigma^{x}(y, \cdot)$, $y \in
  U(x)$;  moreover, for any $\mu \in (0,1]$ the equality
 $$l(x, \delta^{x}_{\mu} c) \, =  \, \mu \, l (x, c)$$
\end{enumerate}
\label{tgenep}
\end{theorem}

\paragraph{Proof.}
Notice that the axiom A4+ is all that we need in order to
transform  the proof of theorem 10 \cite{buligadil1} into a proof of this theorem.
Indeed, for this we need the existence of the limits from A4+ and the algebraic
relations from theorem 11 \cite{buligadil1} which are true only from A0, A1, A2.

If $\displaystyle (\delta^{x}_{\varepsilon}y ,c) \in
\mathcal{L}_{\varepsilon}(X,d,\delta)$ then $\displaystyle 
(x, \Sigma^{x}_{\varepsilon}(y, \cdot) c) \in
\mathcal{L}_{\varepsilon}(X,d,\delta)$ and moreover 
$$l_{\varepsilon}(\delta^{x}_{\varepsilon}y ,c) \, = \, l_{\varepsilon} 
(x, \Sigma^{x}_{\varepsilon}(y, \cdot) c)$$
Indeed, this is true because of the equality: 
$$\delta^{\delta^{x}_{\varepsilon} y}  c \, = \, \delta^{x}_{\varepsilon} 
\Sigma^{x}_{\varepsilon}(y, \cdot) c$$
By passing to the limit with $\varepsilon \rightarrow 0$ and using A3L and A4+ we
get 
$$l(x, c) \, = \, l(x, \Sigma^{x}(y, \cdot) c)$$

For any $\varepsilon, \mu > 0$ (and sufficiently small) 
$\displaystyle (x,c) \in \mathcal{L}_{\varepsilon \mu}(X,d,\delta)$ is equivalent with 
$\displaystyle (x, \delta^{x}_{\mu} c) \in
\mathcal{L}_{\varepsilon}(X,d,\delta)$ 
and moreover: 
$$l_{\varepsilon}(x, \delta^{x}_{\mu} c) \, = \, \frac{1}{\varepsilon} \, 
l_{d}(\delta^{x}_{\varepsilon \mu} c) \, = \, \mu \, l_{\varepsilon \mu} (x, c)$$
We pass to the limit with $\varepsilon \rightarrow 0$ and we get the desired
equality. \quad $\square$

\subsection{Differentiability with respect to a pair of dilation structures}

For any pair of strong dilation structures or length dilation structures there is an associated  notion  of differentiability (section 7.2 \cite{buligadil1}). 
First we need the definition of a morphism of conical groups. 

\begin{definition}
 Let $(N,\delta)$ and $(M,\bar{\delta})$ be two  conical groups. A function $f:N\rightarrow M$ is a conical group morphism if $f$ is a group morphism and for any $\varepsilon>0$ and $u\in N$ we have 
 $\displaystyle f(\delta_{\varepsilon} u) = \bar{\delta}_{\varepsilon} f(u)$. 
\label{defmorph}
\end{definition}

The definition of the derivative, or differential,  with respect to dilations structures follows. In the case of a pair of Carnot groups this is just the definition of the Pansu derivative introduced in \cite{pansu}.

 \begin{definition}
 Let $(X, d, \delta)$ and $(Y, \overline{d}, \overline{\delta})$ be two 
 strong dilation structures  or length dilation structures and $f: U \subset X \rightarrow Y$ be a continuous function defined on an open subset of $X$. The function $f$ is differentiable in $x \in U$ if there exists a 
 conical group morphism  $\displaystyle D \, f(x):T_{x}X\rightarrow T_{f(x)}Y$, defined on a neighbourhood of $x$ with values in  a neighbourhood  of $f(x)$ such that 
\begin{equation}
\lim_{\varepsilon \rightarrow 0} \sup \left\{  \frac{1}{\varepsilon} \overline{d} \left( f\left( \delta^{x}_{\varepsilon} u\right) ,  \overline{\delta}^{f(x)}_{\varepsilon} D \, f(x)  (u) \right) \mbox{ : } d(x,u) \leq \varepsilon \right\}Ê  = 0 , 
\label{edefdif}
\end{equation}
The morphism $\displaystyle D \, f(x) $ is called the derivative, or differential,  of $f$ at $x$.

\label{defdiffer}
\end{definition}

\subsection{Equivalent dilation structures}
\label{secequi}

In the following we adopt a notion of (local) equivalence of dilation structures. 

\begin{definition}
Two strong dilation structures $(X, \delta , d)$ and $(X,
\overline{\delta} , \overline{d})$   are (locally) equivalent if 
\begin{enumerate}
\item[(a)] the identity  map $\displaystyle id: (X, d) \rightarrow (X, \overline{d})$ is bilipschitz, uniformly on compact sets, that is for any compact set $K \subset X$ there  are numbers $R=R(K) >0$ and $c=c(K), C= C(K)> 0$such that for any $x \in K$ the restriction of the identity on the ball $\displaystyle B_{d}(x,R)$ is bilipschitz, with Lipschitz constant smaller  than $C$ and Lipschitz constant of the inverse smaller than by $\displaystyle \frac{1}{c}$,  
\item[(b)]  for any $x \in X$ there are functions $\displaystyle P^{x}, Q^{x}$ (defined for $u \in X$ sufficiently close to $x$) such that  
\begin{equation}
\lim_{\varepsilon \rightarrow 0} \frac{1}{\varepsilon} \overline{d} \left( \delta^{x}_{\varepsilon} u ,  \overline{\delta}^{x}_{\varepsilon} Q^{x} (u) \right)  = 0 , 
\label{dequiva}
\end{equation}
\begin{equation}
 \lim_{\varepsilon \rightarrow 0} \frac{1}{\varepsilon} d \left( \overline{\delta}^{x}_{\varepsilon} u ,  
 \delta^{x}_{\varepsilon} P^{x} (u) \right)  = 0 , 
\label{dequivb}
\end{equation}
uniformly with respect to $x$, $u$ in compact sets. 
\end{enumerate}
\label{dilequi}
\end{definition}

We shall keep the  word "local" further, only if needed.  

\begin{proposition}
 $(X, d, \delta)$ and $(X, , \overline{d}, \overline{\delta})$  are equivalent  if and 
only if 
\begin{enumerate}
\item[(a)] the identity  map $\displaystyle id: (X, d) \rightarrow (X,
\overline{d})$ is locally  bilipschitz, 
\item[(b)]  for any $x \in X$ there are conical group morphisms: 
 $$\displaystyle P^{x}: T_{x}(X, \overline{\delta} , \overline{d}) 
 \rightarrow T_{x} (X, \delta , d) \mbox{ and } \displaystyle  Q^{x}: T_{x} (X, \delta , d) \rightarrow 
 T_{x}(X, \overline{\delta} , \overline{d})$$
  such that the following limits exist  
\begin{equation}
\lim_{\varepsilon \rightarrow 0}  \left(\overline{\delta}^{x}_{\varepsilon}\right)^{-1}  \delta^{x}_{\varepsilon} (u) = Q^{x}(u) , 
\label{dequivap}
\end{equation}
\begin{equation}
 \lim_{\varepsilon \rightarrow 0}  \left(\delta^{x}_{\varepsilon}\right)^{-1}  \overline{\delta}^{x}_{\varepsilon} (u) = P^{x}(u) , 
\label{dequivbp}
\end{equation}
and are uniform with respect to $x$, $u$ in compact sets. 
\end{enumerate}
\label{pdilequi}
\end{proposition}

 The next theorem shows a link between the tangent bundles of equivalent dilation structures. 
 
 \begin{theorem} 
 Let $(X, d, \delta)$ and $(X, \overline{d}, \overline{\delta})$  be  equivalent
 strong  dilation structures.  Then for any $x \in X$ and 
 any $u,v \in X$ sufficiently close to $x$ we have:
 \begin{equation}
 \overline{\Sigma}^{x}(u,v) = Q^{x} \left( \Sigma^{x} \left( P^{x}(u) , P^{x}(v) \right)\right) . 
 \label{isoequiv}
 \end{equation}
 The two tangent bundles  are therefore isomorphic in a natural sense. 
 \label{tisoequiv}
 \end{theorem}

\subsection{Distribution of a dilation structure}

 Let $(X,d,\delta)$ be a strong dilation structure or a length dilation structure. We have then 
 a notion of differentiability for curves in $X$, seen as continuous functions from (a open interval  in) $\mathbb{R}$, with the usual dilation structure, to $X$ with the dilation structure $(X,d,\delta)$.  
 Further we want to see what  differentiability in the sense of definition \ref{defdiffer} means for curves. 
In proposition \ref{pintrinsicd} we shall arrive  to a  notion of a distribution in a dilation structure, with the geometrical meaning of a cone of all possible derivatives of curves passing through a point. 
 
\begin{definition}
In a normed conical group $N$   we shall denote by $D(N)$ the set of all 
$u\in N$ with the property that $\displaystyle 
\varepsilon \in ((0,\infty),+) \mapsto \delta_{\varepsilon} u \in N$ is a morphism of
groups.
\label{defdisn}
\end{definition}
$D(N)$ is always non empty, because it contains the neutral element of $N$. 
$D(N)$ is also a cone, with dilations $\displaystyle   \delta_{\varepsilon}$, 
and a closed set.

 \begin{proposition}
 Let $(X,d,\delta)$ be a strong dilation structure or a length dilation structure and  
 let $\displaystyle c : [a, b] \rightarrow (X,d)$ be a continuous curve. For any 
 $x \in X$ and any $\displaystyle y \in T_{x}(X,d,\delta)$ we denote by 
 $$inv^{x}(y) \, = \, \Delta^{x}(y, x)$$ 
 the inverse of $y$ with respect to the group operation in  $\displaystyle  T_{x}(X,d,\delta)$. 
 Then the following are equivalent: 
 \begin{enumerate}
 \item[(a)] $c$ is derivable in $t \in (a,b)$ with respect to the dilation structure $(X,d,\delta)$; 
 \item[(b)] there exists 
  $\displaystyle \dot{c}(t) \in D(T_{c(t}(X,d,\delta))$ such that 
$$\frac{1}{\varepsilon} d(c(t+\varepsilon) , \delta_{\varepsilon}^{c(t)} \dot{c}(t)) \rightarrow 0 $$
$$\frac{1}{\varepsilon} d(c(t-\varepsilon) , \delta_{\varepsilon}^{c(t)} 
inv^{c(t)}( \dot{c}(t))) \rightarrow 0  $$
\end{enumerate}
\label{pintrinsicd}
\end{proposition}

\paragraph{Proof.}  
 It is  straightforward that a conical group morphism $f: \mathbb{R} \rightarrow N$ is defined by its value $f(1)\in N$. Indeed, for any $a>0$ we have $\displaystyle f(a) = \delta_{a} f(1)$ and for any 
 $a<0$ we have $\displaystyle f(a) = \delta_{a} f(1)^{-1}$. From the morphism property we also 
 deduce that  
 $$\delta v = \left\{ \delta_{a} v \mbox{ : } a>0 , v=f(1) \mbox{ or } v=f(1)^{-1} \right\}$$
 is a one parameter group and that for all $\alpha, \beta >0$ we have 
 $\displaystyle \delta_{\alpha+\beta} u = \delta_{\alpha}u \,  \delta_{\beta}u$. We have
 therefore a bijection between conical group morphisms $f: \mathbb{R} \rightarrow
 (N,\delta)$  and elements of $D(N)$. 
 
 The curve  
 $\displaystyle c : [a,b] \rightarrow (X,d)$ is derivable in $t \in (a,b)$ if 
 and only if there is a morphism of normed conical groups 
 $\displaystyle f: \mathbb{R} \rightarrow T_{c(t}(X,d,\delta)$ such that 
 for any $a \in \mathbb{R}$ we have 
 $$\lim_{\varepsilon \rightarrow 0} \frac{1}{\varepsilon} \, d(c(t+\varepsilon a), 
 \delta^{c(t)}_{\varepsilon} f(a)) \ = \ 0$$
 Take  $\displaystyle \dot{c}(t) = f(1)$. Then $\displaystyle \dot{c}(t) \in 
 D(T_{c(t}(X,d,\delta))$. For any $a > 0$ we have $\displaystyle f(a) =
 \delta^{c(t)}_{a} \dot{c}(t)$; otherwise if $a < 0$ we have 
 $\displaystyle f(a) =
 \delta^{c(t)}_{a} \, inv^{c(t)} \,\dot{c}(t)$. This implies the equivalence 
 stated on the proposition.
 \quad $\square$

\section{Supplementary properties of dilation structures}

At this level of generality, dilation structures come in many flavors. Further we shall propose two supplementary properties which may be satisfied by a dilation structure: the Radon-Nikodym property and the property of being tempered. It will turn out that sub-riemannian spaces may be endowed with dilation structures having the RNP (therefore true in particular for riemannian spaces), but genuinely sub-riemannian spaces don't have dilation structures which are tempered, in contradistinction to the riemannian spaces. 

\subsection{The Radon-Nikodym property}
\label{radon}

 \begin{definition} 
A strong dilation structure or a length dilation structure
 $(X,d, \delta)$ has the  Radon-Nikodym property (or rectifiability
 property, or RNP) if  any  
Lipschitz curve $\displaystyle c : [a,b] \rightarrow (X,d)$ is derivable 
almost everywhere. 
\label{defrn}
 \end{definition}

\paragraph{Three examples.}

The first example is obvious. For  $\displaystyle (X,d)  =  ( \mathbb{V}, d)$, a real, finite dimensional,
normed vector space, with distance $d$ induced by the norm, the (usual) 
 dilations $\displaystyle \delta^{x}_{\varepsilon}$ are given by:  
$$ \delta_{\varepsilon}^{x} y \ = \ x + \varepsilon (y-x) $$
Dilations are defined everywhere.  Axioms 0,1,2 are obviously 
true. Concerning the axiom A3, remark that for any $\varepsilon > 0$, $x,u,v \in X$ we 
have $\displaystyle \frac{1}{\varepsilon} d(\delta^{x}_{\varepsilon} u , 
\delta^{x}_{\varepsilon} v ) \ = \ d(u,v)$. It follows that $\displaystyle d^{x} = d$ for any $x \in X$. Concerning the axiom A4, remark that  for any $\varepsilon > 0$ and any $x,u,v \in X$ we have
$$\Delta_{\varepsilon}^{x}(u,v) \ = \ 
x + \varepsilon  (u-x) + \frac{1}{\varepsilon} \left( x+ \varepsilon(v-x) - x - \varepsilon(u-x) \right) \ = \  x + \varepsilon  (u-x) + v - u$$ 
therefore this quantity converges to $\displaystyle x + v - u \ = \ x + (v - x) - (u - x)$. 
as $\varepsilon \rightarrow 0$. 
For this dilation structure, the RNP as in definition \ref{defrn} is just the usual Radon-Nikodym property.

Further is an example of a dilation structure which does not have the
Radon-Nikodym property.  Take $\displaystyle X = \mathbb{R}^{2}$ with the euclidean 
distance $\displaystyle d$. For any $z = 1+ i \theta \in \mathbb{C}$, with $\theta \in \mathbb{R}$,  we define dilations 
$$\delta_{\varepsilon} x = \varepsilon^{z} x  \ .$$
Then  $\displaystyle (\mathbb{R}^{2},d, \delta)$ 
is a dilation structure, with  dilations $\displaystyle \delta^{x}_{\varepsilon} y = x + \delta_{\varepsilon} (y-x)$.

Two such dilation structures (constructed with the help of complex numbers 
$1+ i \theta$ and $1+ i \theta'$) are equivalent if and only if $\theta = \theta'$.  
Moreover, if $\theta \not = 0$ then there are no non trivial 
Lipschitz curves in $X$ which are differentiable almost everywhere. It means
that such a dilation structure does not have the Radon-Nikodym property. 

More than this, such a dilation structure does not satisfy the obviously reformulated Rademacher theorem (which states that a Lipschitz function is derivable -- in the sense of dilation structures -- almost everywhere with respect to the $2$-Hausdorff measure). Indeed  any holomorphic and Lipschitz function from $X$ to $X$ (holomorphic in the usual sense on $X = \mathbb{R}^{2} = \mathbb{C}$) is differentiable almost everywhere (classically and, equivalently, in the sense of dilation structures), but there are Lipschitz functions from $X$ to $X$ which are not differentiable almost everywhere: it suffices to take a  $\displaystyle \mathcal{C}^{\infty}$ function from  $\displaystyle \mathbb{R}^{2}$ to $\displaystyle \mathbb{R}^{2}$ which is not holomorphic.

The last example concerns riemannian manifolds endowed with dilation structures as in proposition \ref{priemann}. 

\begin{proposition}
The dilation structure of a riemannian manifold constructed from the geodesic spray, as in proposition \ref{priemann}, has the RNP from definition \ref{defrn}.
\end{proposition}

\paragraph{Proof.} Indeed, locally this dilation structure is trivially equivalent to the first dilation structure which was constructed in section \ref{riemann}. That dilation structure is just the one of $\displaystyle \mathbb{R}^{n}$ with the usual dilations and an euclidean distance. The RNP in this case is just the usual RNP, therefore, by corollary \ref{corequi} (from the section \ref{subsdist}), we get the result. \quad $\square$

\subsection{Radon-Nikodym property, representation of length,  distributions}
\label{subsdist}

 \begin{theorem}
Let $(X,d,\delta)$ be a strong dilation structure with the Radon-Nikodym property, 
over a complete length metric space $(X,d)$. Then   for any $x, y \in X$ we have 
$$d(x,y) \ = \ \inf \left\{ \int_{a}^{b} d^{c(t)}(c(t),\dot{c}(t)) \mbox{ d}t  \mbox{ :
} c:[a,b]\rightarrow X \mbox{ Lipschitz }, \right. $$
$$\left.  c(a) = x , c(b) = y \right\}  $$
\label{fleng}
\end{theorem}

\paragraph{Proof.}
By theorem \ref{tupper},  for almost every $t\in(a,b)$ 
the upper dilation of $c$ in $t$ can be expressed as the limit
$$Lip(c)(t) = \lim_{s\rightarrow t} \frac{d(c(s),c(t))}{\mid s-t \mid} $$
For a  dilation structure with the RNP,  for almost every $t \in [a,b]$ there is 
$\displaystyle \dot{c}(t) \in D(T_{c(t)} X)$ such that 
$$\displaystyle \frac{1}{\varepsilon} d(c(t+\varepsilon) , \delta_{\varepsilon}^{c(t)} \dot{c}(t)) \rightarrow 0$$ 
It follows that  for almost every $t \in [a,b]$ we have
$$\displaystyle Lip(c)(t) = \lim_{\varepsilon\rightarrow 0} 
\frac{1}{\varepsilon} d(c(t+\varepsilon),c(t)) = d^{c(t)}(c(t),\dot{c}(t))$$ 
which imples the representation of length. \quad $\square$

As a  consequence, the distance $d$ is uniquely determined by the 
distribution in the sense of dilation structures and  the norm on it.
  
 \begin{corollary}
Let $(X,d,\delta)$ and $(X,\bar{d},\bar{\delta})$ be two strong dilation structures 
with the Radon-Nikodym property , which are also complete length metric spaces, 
such that for any $x \in X$ we have  $\displaystyle D(T_{x}(X,d,\delta)) =
D(T_{x}(X,\bar{d}, \bar{\delta}))$ and $\displaystyle d^{x}(x,u ) = \bar{d}^{x}(x, u)$ for any 
$\displaystyle u \in  D(T_{x}(X,d,\delta))$. Then $d = \bar{d}$.
 \end{corollary}

Another consequence is that the RNP is transported by the equivalence of dilation structures. More precisely we have the following. 
 
\begin{corollary}
Let  $(X, d, \delta)$ and $(X, \overline{d}, \overline{\delta})$  be  
equivalent strong dilation structures. Then for any $x \in X$ we have 
$$\displaystyle Q^{x} (D(T_{x}(X, \delta , d))) \ = \ D(T_{x}(X, \overline{\delta} ,
\overline{d}))  $$

If  $(X, d, \delta)$ has the Radon-Nikodym property , then 
$(X, \overline{d}, \overline{\delta})$ has the same property. 

Suppose that $(X, d, \delta)$ and $(X, \overline{d}, \overline{\delta})$  are 
complete length spaces with the Radon-Nikodym property . If the functions 
$\displaystyle P^{x}, Q^{x}$ from definition \ref{dilequi} (b) are isometries, 
then $\displaystyle d = \overline{d}$. 
\label{corequi}
\end{corollary}

 \subsection{Tempered dilation structures}
\label{stemp}

The notion of a tempered dilation structure extends  
the results of Venturini \cite{venturini} and Buttazzo, De Pascale and 
Fragala \cite{buttazzo1} (propositions 2.3, 2.6 
and a part of theorem 3.1)  to dilation structures. 

The following definition associates  a class of distances 
$\mathcal{D}(X,\bar{d}, \bar{\delta})$  to a strong dilation structure 
$(X, \bar{d}, \bar{\delta})$. This is a  which generalization of  the class of distances 
$\mathcal{D}(X)$ from \cite{buttazzo1}, definition 2.1. 

\begin{definition}
To a  strong dilation structure $(X, \bar{d}, \bar{\delta})$   
 we associte  the class $\mathcal{D}(X,
\bar{d}, \bar{\delta})$ of all length distance functions $d$ on $X$ such that  for any $\varepsilon > 0$ and any $x, u, v$ sufficiently close the are constants $0 < c < C$ with the property  
\begin{equation}
c \, \bar{d}^{x}(u,v) \, \leq \, \frac{1}{\varepsilon} \,
d(\bar{\delta}^{x}_{\varepsilon} u , \bar{\delta}^{x}_{\varepsilon} v ) \, \leq 
\, C \, \bar{d}^{x}(u,v) 
\label{new2.3}
\end{equation}
The dilation structure $(X, \bar{d}, \bar{\delta})$ is tempered if 
 $\bar{d} \in \mathcal{D}(X, \bar{d}, \bar{\delta})$. 

On $\mathcal{D}(X, \bar{d}, \bar{\delta})$ we put the topology of uniform
convergence (induced by distance $\bar{d}$) on compact subsets of $X \times
X$. 
\label{dtempered}
\end{definition}

To any distance $d \in \mathcal{D}(X, \bar{d}, \bar{\delta})$ we associate
the function: 
$$\phi_{d}(x, u) \, = \, \limsup_{\varepsilon \rightarrow 0}
\frac{1}{\varepsilon} \, d(x, \delta^{x}_{\varepsilon} u ) $$
defined for any $x, u \in X$ sufficiently close. We have therefore 
\begin{equation}
c \, \bar{d}^{x}(x, u) \, \leq \, \phi_{d}(x,u) \, \leq \, C \, \bar{d}^{x}(x,u)
\label{new2.6}
\end{equation}

Notice that if $d \in \mathcal{D}(X, \bar{d}, \bar{\delta})$ then for any 
$x, u, v$ sufficiently close we have the relations 
$$- \bar{d}(x,u) \, O(\bar{d}(x,u)) \, + \, c \,  \bar{d}^{x}(u,v) \, \leq $$
$$ \leq \, d(u,v) \, \leq \,   \, C \, \bar{d}^{x}(u,v) \, + \, 
 \bar{d}(x,u) \, O(\bar{d}(x,u))$$

\paragraph{Important remark.} If $c: [0,1] \rightarrow X$ is a $d$-Lipschitz curve and 
$d \in \mathcal{D}(X, \bar{d}, \bar{\delta})$ then we may decompose it 
in a finite family of curves $\displaystyle c_{1}, ... , c_{n}$ (with $n$ depending on $c$) 
such that there are $\displaystyle x_{1}, ... , x_{n} \in X$ with 
$\displaystyle c_{k} $ is $\displaystyle \bar{d}^{x_{k}}$-Lipschitz. Indeed, the image of the 
curve $c([0,1])$ is compact, therefore we may cover it with a finite number of
balls $\displaystyle B(c(t_{k}), \rho_{k}, \bar{d}^{c(t_{k})})$ and apply 
(\ref{new2.3}). If moreover $(X, \bar{d}, \bar{\delta})$ is tempered then 
it follows that $c: [0,1] \rightarrow X$  $d$-Lipschitz curve is equivalent
with $c$ $\bar{d}$-Lipschitz curve. 

By using the same arguments as in the proof of theorem \ref{fleng}, we get the
following extension of proposition 2.4 \cite{buttazzo1}. 

\begin{proposition}
If $(X, \bar{d}, \bar{\delta})$ is tempered, with the Radon-Nikodym 
property, and $d \in \mathcal{D}(X, \bar{d}, \bar{\delta})$ then 
$$d(x,y) \ = \ \inf \left\{ \int_{a}^{b} \phi_{d}(c(t),\dot{c}(t)) \mbox{ d}t  \mbox{ :
} c:[a,b]\rightarrow X \mbox{ $\bar{d}$-Lipschitz }, \right. $$
$$\left.  c(a) = x , c(b) = y \right\}  $$
\label{new2.4}
\end{proposition}

The next theorem is a generalization of a part of theorem 3.1
\cite{buttazzo1}.  

\begin{theorem}
Let $(X, \bar{d}, \bar{\delta})$ be a strong dilation structure which is 
tempered, with the Radon-Nikodym 
property, and $\displaystyle d_{n} \in \mathcal{D}(X, \bar{d}, \bar{\delta})$ 
 a sequence  of distances converging to $d \in \mathcal{D}(X, \bar{d},
 \bar{\delta})$. Denote by $\displaystyle L_{n}, L$ the length functional induced
 by the distance $\displaystyle d_{n}$, respectively by $d$. 
 Then $\displaystyle L_{n}$ $\Gamma$-converges to $L$. 
 \label{new3.1}
 \end{theorem}
 
 \paragraph{Proof.}
We have to prove the liminf inequality and the existence of a recovery sequence, i.e.  parts (a), (b) respectively,  of the definition \ref{defgammamet} of $\Gamma$-convergence of length functionals.  The proof is inspired by the one of implication (i) $\Rightarrow$ (iii) from theorem  3.1 \cite{buttazzo1}, p. 252-253,  we only need to 
replace everywhere expressions like  $\mid x - y\mid$ by $\bar{d}(x,y)$ and 
use proposition \ref{new2.4}, relations (\ref{new2.6}) and (\ref{new2.3})
instead of respectively proposition 2.4 and relations (2.6) and (2.3)
\cite{buttazzo1}.

\paragraph{Proof of (a).} Let us take any sequence of curves $\displaystyle (c_{n})_{n}$, with 
$\displaystyle d_{n}$-Lipschitz curve $\displaystyle c_{n}: [0,1] \rightarrow X$ for every $n$. We suppose that $\displaystyle c_{n}$ converges uniformly to the curve $c$. We want to prove that  
$$L(c) \leq \liminf_{n \rightarrow \infty} L_{n}(c_{n})$$
For any $\eta > 0$ there is a number $m = m(\eta)$ and a  division $\displaystyle \Delta_{\eta} = \left\{ t_{0} = 0, ... , t_{m} = 1 \right\}$ of $[0,1]$ such that 
$$L(c) - \eta \leq \sum_{i = 0}^{m-1} d(c(t_{i}), c(t_{i+1}))$$
From the uniform convergence of $\displaystyle c_{n}$ to $c$ it follows that there is a compact set $K \subset X \times X$ and a number $N = N(K)$ such that for any $n \geq N$ and any $s,t \in [0,1]$ we have $\displaystyle (c_{n}(s), c_{n}(t)) \in K$. 

For any $n \geq N(K)$, 
$$\mid d_{n}(c_{n}(t_{i}), c_{n}(t_{i+1})) - d (c(t_{i}), c(t_{i+1})) \mid \leq 
\mid d_{n}(c_{n}(t_{i}), c_{n}(t_{i+1})) - d (c_{n}(t_{i}), c_{n}(t_{i+1})) \mid + $$ 
$$ + \mid d(c_{n}(t_{i}), c_{n}(t_{i+1})) - d (c(t_{i}), c(t_{i+1})) \leq \sup_{K} \mid d_{n} - d \mid + 2 \sup_{[0,1]} d(c_{n}, c) \leq \varepsilon_{n}$$
where $\displaystyle \varepsilon_{n} = \sup_{K} \mid d_{n} - d \mid + 2 C \sup_{[0,1]} \bar{d}d(c_{n}, c)$ and the number $C>0$ comes from $d \in \mathcal{D}(X, \bar{d}, \bar{\delta})$, see the comments after definition \ref{dtempered}. From this inequality we obtain: 
$$L(c) \leq \eta + L_{n}(c_{n}) + \varepsilon_{n} m(\eta)$$
We pass to the limit with $n \rightarrow \infty$ and we use $\displaystyle \varepsilon_{n} \rightarrow 0$ (from the convergence of $\displaystyle c_{n}$ to $c$) and we get: 
$$L(c) \leq \eta + \liminf_{n \rightarrow \infty} L_{n}(c_{n})$$ 
which gives the desired liminf inequality by the fact that $\eta$ is arbitrary.

\paragraph{Proof of (b).} We have the curve $c$ which is $d$-Lipschitz (therefore $\bar{d}$-Lipschitz) and we want to construct a recovery sequence of curves $\displaystyle c_{n}$. 

Let us consider an increasing sequence $k(n)$ of natural numbers, for the moment without any supplementary assumption. For each $n$ we division $[0,1]$ into $k(n)$ intervals of equal length, denote by $\displaystyle t_{i}^{n}$ the elements of the division, $i = 0, ..., k(n)$, and we define the curve $c_{n}$ to be  one such that $c_{n}(t_{i}^{n}) = c(t_{i}^{n})$ for all $i$ and, denoting by $\displaystyle c^{i}_{n}$ the restriction of $\displaystyle c_{n}$ to the interval $\displaystyle [t_{i}^{n}, t_{i+1}^{n}]$,  such that $\displaystyle c^{i}_{n}$ to be an almost geodesic with respect to $\displaystyle d_{n}$, that is 
\begin{equation}
L_{n}(c^{i}_{n}) \leq d_{n}(c(t_{i}^{n}), c(t_{i+1}^{n})) + \frac{1}{2^{k(n)}}
\label{nhere}
\end{equation}
We want to prove that $\displaystyle c_{n}$ is a recovery sequence, for a well chosen sequence 
$k(n)$. 

It is not restricting the generality to suppose that the length $L(c)$ is equal to one. 
Let $K \subset X \times X$ be a compact set which contains the set $\left\{ (c(t),y) 
\mbox{ : }  d(c(t), y) \leq 1 \right\}$. We choose then the sequence $k(n)$ to be one with the property 
$$\lim_{n \rightarrow \infty} k(n) \sup_{K} \mid d_{n} - d \mid = 0$$
Then, there is a number $N(K)$ such that for any $ n \geq N(K)$  and any For any $t \in [0,1]$ and any $s,t \in [0,1]$ we have $\displaystyle (c_{n}(s), c_{n}(t)) \in K$. 

For any $n \geq N(K)$ and any $t \in [0,1]$ we denote by $\displaystyle [t^{-}_{n}, t^{+}_{n}$ the interval of the division of $[0,1]$ into $k(n)$ intervals of equal length, where $t$ belongs. Let also 
$\displaystyle c_{n}^{t}$ be the restriction of $\displaystyle c_{n}$ to the interval 
$\displaystyle [t^{-}_{n}, t^{+}_{n}$. 

 Then 
$\displaystyle \sup_{[0,1]} \bar{d}(c_{n}(t), c(t)) \leq A_{n} + B_{n}$, where 
$\displaystyle A_{n} = \sup_{[0,1]} \bar{d}(c(t^{+}_{n}), c(t))$ and 
$\displaystyle B_{n} = \sup_{[0,1]} \bar{d}(c_{n}(t^{+}_{n}), c_{n}(t))$. Trivially $\displaystyle A_{n} \rightarrow 0$ as $n \rightarrow \infty$. The term $\displaystyle B_{n}$ may be estimaded as follows. By proposition \ref{new2.4}, relations (\ref{new2.6}),  (\ref{new2.3}) and (\ref{nhere}), there are numbers $0< c < C$ such that
$$c \, \bar{d}(c_{n}(t^{+}_{n}), c_{n}(t)) \leq L_{n}(c^{t}_{n}) + \frac{1}{2^{k(n)}} \leq C \,  \bar{d}( c(t^{-}_{n}) , c(t^{+}_{n})) + \frac{1}{2^{k(n)}}$$
We have then 
$$L(c) \leq \sum_{i=0}^{k(n) -1} d(c(t_{i}^{n}), c(t_{i+1}^{n})) = 
\sum_{i=0}^{k(n) -1} d_{n}(c_{n}(t_{i}^{n}), c_{n}(t_{i+1}^{n})) + $$ 
$$+ \sum_{i=0}^{k(n) -1} \left(  d(c(t_{i}^{n}), c(t_{i+1}^{n})) - d_{n}(c(t_{i}^{n}), c(t_{i+1}^{n})) \right) \geq L_{n}(c_{n}) - \frac{k(n)}{2^{k(n)}} -   k(n) \sup_{K} \mid d - d_{n} \mid $$ 
We apply $\displaystyle \limsup_{n \rightarrow \infty}$ to his inequality and we prove the claim. \quad $\square$

As a corollary we obtain   a large class of examples of length dilation structures.

\begin{corollary}
If  $(X, \bar{d}, \bar{\delta})$ is a strong dilation structure which is 
tempered and it has  the Radon-Nikodym property then it is a length dilation structure. 
\label{cortemp}
\end{corollary}

\paragraph{Proof.}
Indeed, from the hypothesis we deduce that $\displaystyle
\bar{\delta}^{x}_{\varepsilon} \bar{d} \, \in \, 
 \mathcal{D}(X, \bar{d}, \bar{\delta})$. For any sequence 
 $\displaystyle \varepsilon_{n} \rightarrow 0$ we thus obtain a sequence of 
 distances $\displaystyle d_{n} \, = \, \bar{\delta}^{x}_{\varepsilon_{n}} 
 \bar{d}$ converging to $\bar{d}^{x}$. We apply now theorem \ref{new3.1} 
 and we get the result. 
\quad $\square$

We arive therefore to the following characterization of riemannian manifolds.

\begin{theorem}
The dilation structure associated to  a riemannian manifold, as in proposition \ref{priemann}, is tempered and is a length dilation structure. 

Conversely, if  $(X, d, \delta)$ is a strong dilation structure which is 
tempered,  it has  the Radon-Nikodym property and moreover for any $x \in X$ the tangent space is a commutative local group, then any open, with compact closure subset of $X$ can be endowed with a $\displaystyle C^{1}$ riemannian structure which gives a distance $d'$ which is bilipschitz equivalent with $d$. 
\label{tchariem}
\end{theorem}

\paragraph{Proof.}  We already know that this dilation structure has the RNP. It is also tempered, because of the estimate (\ref{mainriem}) from the proof of the proposition 
\ref{pcurvature}. By corollary \ref{cortemp} it is a length dilation structure. 

For the converse assertion, remark that the only local conical groups, which are locally compact and admit a one parameter group of dilations (that is the abelian group $\Gamma$ is 
$(0,+\infty)$ are (open neighbourhoods of $0$ in) $\displaystyle \mathbb{R}^{n}$. From the fact that the dilation structure is tempered, it follows that locally $(X,d)$ is bilipschitz with 
$\displaystyle \mathbb{R}^{n}$ endowed with the norm constructed as (\ref{new2.6}). But any 
norm on $\displaystyle \mathbb{R}^{n}$ is bilipschitz with an euclidean norm. 
\quad $\square$

\section{Dilation structures on sub-riemannian manifolds}
\label{srsmooth}

In \cite{buligasr}, followed in this section,   we proved that we can associate dilation structures to regular sub-Riemannian manifolds. This result is the source of inspiration of the notion of a coherent projection, section \ref{cohp}.

In this section we use differential geometric tools, mainly the normal frames, definition \ref{defnormal}. This has been proved by Bella\"{\i}che  \cite{bell}, starting with theorem 4.15 and ending in the first half of section 7.3 (page 62). We shall not give an exposition of this long proof, although a streamlined version of it would be very useful.

From the existence of normal frames  we shall get the existence  of certain  dilation structures regular sub-riemannian manifolds, theorem \ref{structhm}.  
 From this, according to the general theory of dilation structures, via Siebert type results alsoo, follow all  classical results concerning the structure of the tangent space to a point of a regular sub-riemannian manifold. 

In particular, this is evidence for the fact that the classical differential calculus is needed   only in the part concerning the existence of normal frames and after this stage an intrinsic way of reasoning is possible.

Let us compare with, maybe, the closest results, which are those  of Vodopyanov \cite{vodopis}, \cite{vodokar}. In both approaches c the tangent space to a point is defined only locally, as a neighbourhood of the point, in the manifold. The difference is that Vodopyanov proves the existence of the (locally defined) operation on the tangent space  under very weak regularity assumptions on the sub-riemannian  manifold, by using the differential structure of the underlying manifold.  In distinction, we prove in \cite{buligadil1}, in  an abstract setting,   that the very existence of a dilation structure induces a locally defined operation.

\subsection{Sub-riemannian manifolds}
\label{secsr}

$M$ is a connected,  $n$ dimensional,  real manifold. A ( differential geometric)  
distribution on $M$ is a smooth  subbundle $D$ of $M$. Such a distribution associates 
to any point $x \in M$ a vector space $D_{x} \subset T_{x}M$.
The dimension of the distribution $D$ at point $x \in M$ is $\displaystyle m(x) = \, dim \, D_{x}$. The function $x \in M \mapsto m(x)$ is locally constant, because of the distribution is smooth. We  shall suppose  that the dimension of the distribution is
globally constant and we  denote it by $m$. In general  $m \leq n$. The typical case we are interested in is  $m < n$.

A horizontal curve $c:[a,b] \rightarrow M$ is a curve which is almost everywhere derivable and for almost any $t \in [a,b]$ we have $\displaystyle \dot{c}(t) \in D_{c(t)}$.
The class of horizontal curves will be denoted by $Hor(M,D)$.

\begin{definition}
The distribution $D$ is completely non-integrable if $M$ is locally connected
by horizontal  curves, that is curves in $Hor(M,D)$.
\label{defcnint}
\end{definition}
The Chow condition (C) \cite{chow} is sufficient  for the distribution $D$ to be completely non-integrable.

\begin{theorem} (Chow) Let us suppose there is a positive integer number $k$ (called the rank of the distribution $D$) such that for any $x \in X$ there is a topological  open ball
$U(x) \subset M$ with $x \in U(x)$ such that there are smooth vector fields
$\displaystyle X_{1}, ..., X_{m}$ in $U(x)$ with the property:

(C) the vector fields $\displaystyle X_{1}, ..., X_{m}$ span $\displaystyle
D_{x}$ and these vector fields together with  their iterated
brackets of order at most $k$ span the tangent space $\displaystyle T_{y}M$
at every point $y \in U(x)$.

Then the distribution $D$ is completely non-integrable in the sense of
definition \ref{defcnint}.
\label{tchow}
\end{theorem}

\begin{definition}
A sub-riemannian (SR) manifold is a triple $(M,D, g)$, where $M$ is a
connected manifold, $D$ is a completely non-integrable distribution on $M$, and $g$ is a metric (Euclidean inner-product) on the distribution (or horizontal bundle)  $D$.
\label{defsr}
\end{definition}

With the notations from condition (C), let us  define on $U = U(x)$ a filtration of bundles
as follows. First we define  the class of horizontal vector fields on $U$
$$\mathcal{X}^{1}(U(x),D) \ = \ \left\{ X \in \mathcal{X}^{\infty}(U) \mbox{ : }
\forall y \in U(x) \ , \ X(y) \in D_{y} \right\}$$
Next, we define inductively for all positive integers $j$ the following vector fields: 
$$ \mathcal{X}^{j+1} (U(x),D) \ = \ \mathcal{X}^{j}(U(x),D) \, + \,
 [ \mathcal{X}^{1}(U(x),D),
\mathcal{X}^{j}(U(x),D)]$$
where $[ \cdot , \cdot ]$ denotes the bracket of vector fields. 
All in all we obtain  the  filtration $\displaystyle \mathcal{X}^{j}(U(x),D) \subset
\mathcal{X}^{j+1} (U(x),D)$. By evaluation at  at $y \in U(x)$, we get a filtration 
$$V^{j}(y,U(x),D) \ = \ \left\{ X(y) \mbox{ : } X \in \mathcal{X}^{j}(U(x),D)\right\}$$
According to Chow condition (C),  there is a positive integer $k$ such that
 for all $y \in U(x)$ we have 
$$D_{y} =  V^{1}(y,U(x),D) \subset V^{2}(y,U(x),D) \subset ... \subset
V^{k}(y, U(x),D) = T_{y}M $$
Consequently, to the sub-riemannian manifold is associated the string of
numbers:
$$ \nu_{1}(y) = \dim V^{1}(y, U(x),D) < \nu_{2}(y) = \dim V^{2}(y, U(x),D)
< ... < n = \dim M$$
Generally $k$, $\displaystyle \nu_{j}(y)$ may vary from a point to another. The number $k$ is called the step of the distribution at $y $.

\begin{definition}
The distribution $D$, which satisfies the Chow condition (C),  is regular if  $\displaystyle \nu_{j}(y)$ are constant on the manifold $M$.

The sub-riemannian manifold $M,D,g)$ is regular if $D$ is regular and for any
$x \in M$ there is a topological ball $U(x) \subset M$ with $x \in U(M)$ and
 an orthonormal (with respect to the metric $g$)  family of
  smooth vector fields $\displaystyle \left\{ X_{1}, ..., X_{m}\right\}$ in
  $U(x)$ which satisfy the condition (C).
\label{dreg}
\end{definition}

The lenght of  a horizontal curve is obtained from the metric $g$ by 
$$l(c) \ = \ \int_{a}^{b} \left(g_{c(t)} (\dot{c}(t),
\dot{c}(t))\right)^{\frac{1}{2}} \mbox{ d}t$$

\begin{definition}
The Carnot-Carath\'eodory distance (or CC distance) associated to the sub-riemannian manifold is the distance induced by the length $l$ of horizontal curves:
$$d(x,y) \ = \ \inf \left\{ l(c) \mbox{ : } c \in Hor(M,D) \
, \ c(a) = x \ , \  c(b) = y \right\} $$
\end{definition}

The Chow condition ensures the existence of a horizontal path linking any two sufficiently
closed points, therefore the CC distance is  locally finite. The distance
depends only on the distribution $D$ and metric $g$, and not on the choice of
vector fields $\displaystyle X_{1}, ..., X_{m}$ satisfying the condition (C).
The space $(M,d)$ is locally compact and complete, and the topology
induced by the distance $d$ is the same as the topology of the manifold $M$.
(These important details may be recovered from reading carefully the proofs of Chow theorem given by  Bella\"{\i}che \cite{bell} or Gromov \cite{gromovsr}.)

\subsection{Sub-riemannian dilation structures associated to normal frames}

In the following we suppose that $M$ is a regular sub-riemannian manifold. The dimension of $M$ as a differential manifold is denoted by $n$, the
 step of the regular sub-riemannian manifold $(M,D,g)$ is denoted by $k$, the
dimension of the distribution is $m$, and there are numbers $\displaystyle
\nu_{j}$, $j = 1, ..., k$ such that for any $x \in M$ we have
$\displaystyle dim \, V^{j}(x) = \nu_{j}$. The Carnot-Carath\'eodory distance
is denoted by $d$.Further, we stay in a small open neighbourhood of an arbitrary, but fixed point $\displaystyle x_{0} \in M$. We shall no longer mention the
dependence of various objects on $\displaystyle x_{0}$, on the neighbourhood
$\displaystyle U(x_{0})$, or the distribution $D$.

\begin{definition}
An adapted frame $\displaystyle \left\{ X_{1}, ... , X_{n} \right\}$ is an ordered 
collection of smooth vector fields, constructed according to the following recipe. 

The first $m$ vector fields  $X_{1}, ... , X_{m}$  satisfy
the condition (C). Therefore, for any point $x$ the vectors  $\displaystyle
X_{1}(x), ... , X_{m}(x)$ form a basis for $\displaystyle D_{x}$.  

We associate to any word $\displaystyle a_{1} .... a_{q}$ with letters in the alphabet $1, ... ,m$ the vector field equal to the   multi-bracket$\displaystyle [X_{a_{1}}, [ ... , X_{a_{q}}] ... ]$. We add,  in the lexicographic order, $n-m$ elements to the set 
$\displaystyle \left\{ X_{1}, ... , X_{m} \right\}$ until we get a collection
$\displaystyle \left\{ X_{1}, ... , X_{n} \right\}$ such that:
for any $j = 1, ..., k$ and for any point $x$ the set 
$\displaystyle \left\{X_{1}(x), ..., X_{\nu_{j}}(x) \right\}$ is a basis for
$\displaystyle V^{j}(x)$.
\label{defadapt}
\end{definition}

Let $\displaystyle \left\{ X_{1}, ... , X_{n} \right\}$ be an adapted frame.
For any $j = 1, ..., n$  the degree $\displaystyle deg \, X_{j}$ of the vector
field $\displaystyle X_{j}$  is defined
as the only positive integer $p$ such that for any point $x$ we have
$$X_{j}(x) \in V^{p}_{x} \setminus V^{p-1}(x)$$

In  definition below, the key  are the uniform convergence assumptions. This is in line 
with Gromov suggestions in the last section of Bella\"{\i}che \cite{bell}. 

\begin{definition}
A normal  frame is an adapted frame $\displaystyle \left\{ X_{1}, ... , X_{n} \right\}$ which satisfies the following supplementary properties:
\begin{enumerate}
\item[(a)] the limit 
$$\lim_{\varepsilon \rightarrow 0_{+}} \frac{1}{\varepsilon} \, d\left(
\exp \left( \sum_{1}^{n}\varepsilon^{deg\, X_{i}} a_{i} X_{i} \right)(y), y \right) \ = \ A(y, a) \in
(0,+\infty)$$
exists and is uniform with respect to $y$ in compact sets and $\displaystyle a=(a_{1}, ...,
a_{n}) \in W$, with $\displaystyle W \subset \mathbb{R}^{n}$ compact
neighbourhood of $\displaystyle 0 \in \mathbb{R}^{n}$,
\item[(b)] for any compact set $K\subset M$ with diameter (with respect to the
 CC distance $d$) sufficiently small,  and for any $i = 1, ..., n$ there
are functions 
$$ P_{i}(\cdot, \cdot, \cdot): U_{K} \times U_{K} \times K \rightarrow  \mathbb{R}$$
 with $\displaystyle U_{K} \subset \mathbb{R}^{n}$ a sufficiently small compact neighbourhood of 
$\displaystyle 0 \in \mathbb{R}^{n}$ such that for any   $x \in K$ and any $\displaystyle 
a,b \in U_{K}$ we have
$$\exp \left( \sum_{1}^{n} a_{i} X_{i} \right) (x) \ = \
\exp \left( \sum_{1}^{n} P_{i}(a, b, y) X_{i} \right) \circ \exp \left( \sum_{1}^{n}  b_{i} X_{i} \right) (x) $$
and such that the following limit exists
$$\lim_{\varepsilon \rightarrow 0_{+}}
\varepsilon^{-deg \, X_{i}} P_{i}(\varepsilon^{deg \, X_{j}} a_{j}, \varepsilon^{deg \, X_{k}} b_{k}, x)   \in
\mathbb{R}$$
and it is uniform with respect to $x  \in K$ and $\displaystyle a, b \in U_{K}$.
\end{enumerate}
\label{defnormal}
\end{definition}

The condition (a) definition \ref{defnormal} is a part of the conclusion of Gromov approximation theorem, namely when one point coincides with the center of nilpotentization; also condition (b) is equivalent with a statement of Gromov concerning the convergence of rescaled vector fields to their nilpotentization (an informed reader must at least follow in all details the papers Bella\"{\i}che \cite{bell} and Gromov \cite{gromovsr}, where differential calculus in the classical sense is heavily used). 

In the case of a Lie group $G$ endowed with a left invariant distribution, normal frames are very easy to visualize.  The left-invariant distribution is completely non-integrable if and only it is generated by the left  translation of a vector subspace $D$ of the algebra
$\mathfrak{g} = T_{e}G$ which  generates the whole Lie algebra 
$\mathfrak{g}$. Let us take $\displaystyle \left\{ X_{1}, ..., X_{m}\right\}$ to be  a collection of $m = \, dim \, D$ left-invariant independent vector fields.  Let us define with their help an adapted frame, as explained in definition \ref{defadapt}. This frame is in fact normal.

\begin{definition}
Let $(M,d,g)$ be a  regular sub-riemannian
manifold and let 
$\displaystyle \left\{ X_{1}, ..., X_{n} \right\}$ be a normal frame (locally defined, but for simplicity here we neglect this detail). To this pair we  associate a triple $(M,d, \delta)$, where: $d$ is  the
Carnot-Carath\'eodory distance, and for any point $x \in M$
and any $\varepsilon \in (0,+\infty)$ (sufficiently small if necessary),
the dilation $\displaystyle \delta^{x}_{\varepsilon}$ is defined by  the expression 
$$\delta^{x}_{\varepsilon} \left(\exp\left( \sum_{i=1}^{n} a_{i} X_{i} \right)(x)\right) \  = \
\exp\left( \sum_{i=1}^{n} a_{i} \varepsilon^{deg X_{i}}  X_{i} \right)(x)$$
\label{dsrdil}
\end{definition}

We shall prove that $(M,d, \delta)$ is  a dilation structure. This
allows us to get the main results concerning the infinitesimal geometry of a
regular sub-riemannian manifold.

\begin{theorem} 
 Let $(M,D, g)$ be a regular sub-riemannian manifold and $U \subset M$ an open set which admits a normal frame. Then $(U, d,  \delta)$ is a strong dilation structure. 
\label{structhm}
\end{theorem}

\paragraph{Proof.} 
We only need to prove that A3 and A4 hold. The fact that A3 is true is a result similar to Gromov local approximation theorem \cite{gromovsr}, p. 135, or to Bella\"{\i}che theorem 7.32 \cite{bell}.

\paragraph{A3 is satisfied.} We prove that   the limit 
$\displaystyle \lim_{\varepsilon \rightarrow 0} \frac{1}{\varepsilon} \, d\left(
\delta^{x}_{\varepsilon} u, \delta^{x}_{\varepsilon} v \right) \ = \
d^{x}(u,v)$ exists and it uniform with respect to $x,u,v$ sufficiently close.  We represent $u$ and $v$ with respect to the normal frame: 
$$u \ = \ \exp \left( \sum_{1}^{n} u_{i} X_{i} \right)(x) \quad \quad , \quad
v \ = \ \exp \left( \sum_{1}^{n} v_{i} X_{i} \right)(x)$$
Let us denote by $\displaystyle u_{\varepsilon} = \delta_{\varepsilon}^{x} u = \exp \left( \sum_{1}^{n}
\varepsilon^{deg \, X_{i}} u_{i} X_{i} \right)(x)$. Then 
$$\frac{1}{\varepsilon} \, d\left(
\delta^{x}_{\varepsilon} u, \delta^{x}_{\varepsilon} v \right) \ = \
\frac{1}{\varepsilon} \, d\left(
\delta^{x}_{\varepsilon} \exp \left( \sum_{1}^{n} u_{i} X_{i} \right)(x),
\delta^{x}_{\varepsilon} \exp \left( \sum_{1}^{n} v_{i} X_{i} \right)(x) \right)
\ =$$
$$= \ \frac{1}{\varepsilon} \, d\left(
 \exp \left( \sum_{1}^{n} \varepsilon^{deg \, X_{i}} u_{i} X_{i} \right)(x),
\exp \left( \sum_{1}^{n} \varepsilon^{deg \, X_{i}}  v_{i} X_{i} \right)(x) \right)
$$
Let us  make the  notation: for any $i=1, ..., n$
$ \displaystyle a_{i}^{\varepsilon} \ = \ \varepsilon^{-deg \,
X_{i}} \, P_{i}(\varepsilon^{deg \, X_{j}} v_{j}, \varepsilon^{deg \, X_{k}} u_{k}, x)$. 
By the second part of property (b), definition \ref{defnormal}, the vector
$\displaystyle a^{\varepsilon} \in \mathbb{R}^{n}$ converges to a finite value $\displaystyle a^{0} \in 
\mathbb{R}^{n}$, as $\varepsilon \rightarrow 0$,  uniformly with respect
to $x, u, v$ in compact set. In the same time $\displaystyle u_{\varepsilon}$ converges to $x$, as 
$\varepsilon \rightarrow 0$. But remark that  
$$\frac{1}{\varepsilon} \, d\left(
\delta^{x}_{\varepsilon} u, \delta^{x}_{\varepsilon} v \right) \ = \ \frac{1}{\varepsilon} \, d \left( u_{\varepsilon},
\exp \left( \sum_{1}^{n} P_{i}(\varepsilon^{deg \, X_{j}} v_{j}, \varepsilon^{deg \, X_{k}} u_{k}, x)
X_{i}\right) (u_{\varepsilon}) \right) \ = $$
$$= \ \frac{1}{\varepsilon} \, d \left( u_{\varepsilon},
\exp \left( \sum_{1}^{n} \varepsilon^{deg \, X_{i}} \left( \varepsilon^{-deg \,
X_{i}} \, P_{i}(\varepsilon^{deg \, X_{j}} v_{j}, \varepsilon^{deg \, X_{k}} u_{k}, x) \right)
X_{i}\right) (u_{\varepsilon}) \right) \ = \ $$ 
$$\ = \ \frac{1}{\varepsilon} \, d \left( u_{\varepsilon},
\exp \left( \sum_{1}^{n} \varepsilon^{deg \, X_{i}} a^{\varepsilon}_{i}
X_{i}\right) (u_{\varepsilon}) \right)$$
Now we use the uniform convergence assumptions from definition \ref{defnormal}. For fixed $\eta > 0$ the term 
$$B(\eta, \varepsilon) = \frac{1}{\varepsilon} \, d \left( u_{\eta}, \exp \left( \sum_{1}^{n} \varepsilon^{deg \, X_{i}} a^{\eta}_{i}
      X_{i}\right) (u_{\eta}) \right)$$
converges to a real number $\displaystyle    A(u_{\eta}, a_{\eta})$ as  $\varepsilon \rightarrow 0$, uniformly with respect to 
$\displaystyle u_{\eta}$ and $\displaystyle a_{\eta}$. Since  $\displaystyle u_{\eta}$    converges to $x$  and 
 $\displaystyle a_{\eta}$ converges to   $\displaystyle a^{0}$ as $\eta \rightarrow 0$, by the uniform convergence assumption 
in (a), definition \ref{defnormal} we get that 
$$\lim_{\varepsilon \rightarrow 0}  \frac{1}{\varepsilon} \, d\left(  \delta^{x}_{\varepsilon} u, \delta^{x}_{\varepsilon} v \right) \  = \ 
\lim_{\eta \rightarrow 0} A(u_{\eta}, a_{\eta})  \ = \   A(x, a^{0})$$

\paragraph{A4 is satisfied.} As $\varepsilon$ converges to $0$,  the approximate difference 
$\Delta^{x}_{\varepsilon}(u,v)$ has a limit, which is uniform  with respect to $x,u,v$ sufficiently close. Indeed, 
$$\Delta^{x}_{\varepsilon}(u,v) \ = \ \delta_{\varepsilon^{-1}}^{ u_{\varepsilon}}
\exp \left(
\sum_{1}^{n} \varepsilon^{deg \, X_{i}} v_{i} X_{i} \right)(x)$$
We use the first part of the property (b), definition \ref{defnormal}, in order to
write
$$\exp \left(
\sum_{1}^{n} \varepsilon^{deg \, X_{i}} v_{i} X_{i} \right)(x) \ = \
\exp \left( \sum_{1}^{n} P_{i}(\varepsilon^{deg \, X_{j}} v_{j}, \varepsilon^{deg \, X_{k}} u_{k}, x)
X_{i}\right) (u_{\varepsilon}) $$
We finish the computation:
$$\Delta^{x}_{\varepsilon}(u,v) \ = \ \exp \left( \sum_{1}^{n} \varepsilon^{- \, deg \, X_{i}} \, P_{i}(\varepsilon^{deg \, X_{j}} v_{j}, \varepsilon^{deg \, X_{k}} u_{k}, x)
X_{i}\right) (u_{\varepsilon})$$
As $\varepsilon$ goes to $0$ the point  $\displaystyle u_{\varepsilon}$ converges to
$x$ uniformly with respect to $x,u$ sufficiently close (as a corollary of the
previous theorem, for example). The proof therefore ends by invoking the second
part of the property (b), definition \ref{defnormal}. $\quad \square$

\section{Coherent projections: a dilation structure looks down on another}

The equivalence of dilation structures, definition \ref{dilequi}, may be seen as a composite of two partial order relations. 

\begin{definition}
A strong dilation structure  $(X, \delta , d)$ is looking down on another strong dilation structure  $(X, \overline{\delta} , \overline{d})$ if 
\begin{enumerate}
\item[(a)] the identity  map $\displaystyle id: (X, d) \rightarrow (X, \overline{d})$ is lipschitz and 
\item[(b)]  for any $x \in X$ there are functions $\displaystyle Q^{x}$ (defined for $u \in X$ sufficiently close to $x$) such that  
\begin{equation}
\lim_{\varepsilon \rightarrow 0} \frac{1}{\varepsilon} \overline{d} \left( \delta^{x}_{\varepsilon} u ,  \overline{\delta}^{x}_{\varepsilon} Q^{x} (u) \right)  = 0 , 
\label{dequiva1}
\end{equation}
uniformly with respect to $x$, $u$ in compact sets. 
\end{enumerate}
\label{dilook}
\end{definition}

This leads us to the introduction of coherent projections, i.e. to the study of pairs of dilation structures, one looking down on another. (The name is inspired by the notion of a set looking down on another introduced in \cite{davidsemmes}.) We prefer to work with a pair  $(\bar{\delta}, Q)$ made by a dilation structure and a function $Q$ as in the previous definition, point (b), instead of working with a pair of dilation structures.

\subsection{Coherent projections}
\label{cohp}

\begin{definition}
Let $(X,\bar{d}, \bar{\delta})$ be a strong dilation structure. A coherent projection
of $(X,\bar{d}, \bar{\delta})$  is a function which associates to any $x \in X$ and $\varepsilon
\in (0,1]$ a map 
$\displaystyle Q^{x}_{\varepsilon} : U(x) \rightarrow X$ such that: 
\begin{enumerate}
\item[(I)] $\displaystyle Q^{x}_{\varepsilon} : U(x) \rightarrow  Q^{x}_{\varepsilon}(U(x))$ 
is invertible and the inverse is $\displaystyle Q^{x}_{\varepsilon^{-1}}$; moreover, the following commutation relation holds: for
any $\varepsilon, \mu > 0$ and any $x \in X$ 
$$Q_{\varepsilon}^{x} \, \bar{\delta}^{x}_{\mu} \ = \ \bar{\delta}^{x}_{\mu} \,
Q^{x}_{\varepsilon}$$
\item[(II)] the limit $\displaystyle \lim_{\varepsilon \rightarrow 0} 
Q^{x}_{\varepsilon} u \ = \ Q^{x} u$  
is uniform with respect to $x, u$ in compact sets. 
\item[(III)] for any $\varepsilon, \mu > 0$ and any $x \in X$ we have $\displaystyle 
Q^{x}_{\varepsilon} \, Q^{x}_{\mu} \ = \ Q^{x}_{\varepsilon \mu}$. Also $\displaystyle
Q^{x}_{1} \ = \ id$ and $\displaystyle Q^{x}_{\varepsilon} x \ = \ x$.

\item[(IV)] define $\displaystyle \Theta_{\varepsilon}^{x}(u,v) \ = \ 
\bar{\delta}^{x}_{\varepsilon^{-1}} \,
Q^{\bar{\delta}^{x}_{\varepsilon} Q^{x}_{\varepsilon} u}_{\varepsilon^{-1}} \,
\bar{\delta}^{x}_{\varepsilon} Q^{x}_{\varepsilon} v$.  
Then the limit exists
$$\lim_{\varepsilon \rightarrow 0} \Theta^{x}_{\varepsilon}(u,v) \ = \ \Theta^{x}(u,v)$$
and it is uniform with respect to $x, u , v$ in compact sets.
\end{enumerate}
A coherent projection induces the dilations  $\displaystyle \delta^{x}_{\varepsilon} \ = \
\bar{\delta}^{x}_{\varepsilon} \, Q^{x}_{\varepsilon}$.
\label{defcoh}
\end{definition}

\begin{proposition}

Let $(X,\bar{d}, \bar{\delta})$ be a strong dilation structure and $Q$ a coherent 
projection.  Then: 
\begin{enumerate}
\item[(a)]  the induced dilations commute with the old ones: for any $\varepsilon, \mu > 0$ and any $x \in X$ we have $\displaystyle 
\delta^{x}_{\varepsilon} \, \bar{\delta}^{x}_{\mu} \ = \  \bar{\delta}^{x}_{\mu} \, 
\delta^{x}_{\varepsilon}$. 
\item[(b)] as $\varepsilon \rightarrow 0$ the coherent projection becomes a projection:  for any $x \in X$ we have $\displaystyle Q^{x} \, Q^{x} \ = \ Q^{x}$. 
\item[(c)] the induced dilations  $\delta$ satisfy the 
conditions A1, A2, A4 from definition \ref{defweakstrong}, 
\item[(d)] the following relations are true (we denote by $\Sigma^{x}$, $\Delta^{x}$, ..., the approximate or infinitesimal sum and difference computed with the help of induced dilations):   
\begin{equation}
\Theta_{\varepsilon}^{x}(u,v) \ = \ \bar{\Sigma}^{x}_{\varepsilon}(Q^{x}_{\varepsilon} u, 
\Delta^{x}_{\varepsilon} (u,v))  
\label{neednext1}
\end{equation}
\begin{equation}
\Delta^{x} (u,v) \, = \, \bar{\Delta}^{x}(Q^{x} u, 
\Theta^{x}(u,v)) 
\label{recon}
\end{equation}
\begin{equation}
Q^{x} \Delta^{x}(u,v) \, = \, \bar{\Delta}^{x}(Q^{x} u, Q^{x} v)
\label{morph}
\end{equation}
\end{enumerate}
\label{p1proj}
\end{proposition}

\paragraph{Proof.}
(a) this is a consequence of the commutativity condition (I) (second part). Indeed, we have 
$\displaystyle \delta^{x}_{\varepsilon} \, \bar{\delta}^{x}_{\mu} \ = \ 
\bar{\delta}^{x}_{\varepsilon} \, Q^{x}_{\varepsilon} \, \bar{\delta}^{x}_{\mu} \ = \ 
\bar{\delta}^{x}_{\varepsilon} \, \bar{\delta}^{x}_{\mu} \, Q^{x}_{\varepsilon} \ = \ 
\bar{\delta}^{x}_{\mu} \, \bar{\delta}^{x}_{\varepsilon} \, Q^{x}_{\varepsilon} \ = \ 
 \bar{\delta}^{x}_{\mu} \, \delta^{x}_{\varepsilon}$. 

(b) we pass to the limit $\varepsilon \rightarrow 0$  in 
 the equality $\displaystyle Q^{x}_{\varepsilon^{2}} \ = \ Q^{x}_{\varepsilon}
\, Q^{x}_{\varepsilon}$ and we get, based on condition (II), that $\displaystyle Q^{x} \, Q^{x}
\ = \ Q^{x}$.

(c) Axiom A1 for $\delta$ is equivalent with (III). Indeed,  the equality 
$\displaystyle \delta^{x}_{\varepsilon} \, \delta^{x}_{\mu} \, = \, 
\delta^{x}_{\varepsilon \mu}$  is 
equivalent with: 
$\displaystyle \bar{\delta}^{x}_{\varepsilon \mu} \, Q^{x}_{\varepsilon \mu} \, = \, 
  \bar{\delta}^{x}_{\varepsilon \mu} \, Q^{x}_{\varepsilon} \, Q^{x}_{\mu}$. 
  This is true because $\displaystyle 
Q^{x}_{\varepsilon} \, Q^{x}_{\mu} \, = \, Q^{x}_{\varepsilon \mu}$. 
We also have $\displaystyle  \delta^{x}_{1} \, =  \, \bar{\delta}^{x}_{1} Q^{x}_{1} 
\, = \, Q^{x}_{1} \, = \, id $. 
Moreover $\displaystyle \delta^{x}_{\varepsilon} x \, = \, 
\bar{\delta}^{x}_{\varepsilon} \, Q^{x}_{\varepsilon} x \, = \, 
Q^{x}_{\varepsilon} \,  \bar{\delta}^{x}_{\varepsilon} x \, = \, 
Q^{x}_{\varepsilon} x  \, = \, x$. 
Let us compute now: 
$$\Delta^{x}_{\varepsilon} (u,v) \ = \ \delta^{\delta^{x}_{\varepsilon} u}_{\varepsilon^{-1}} 
\, \delta^{x}_{\varepsilon} v \ = \ \bar{\delta}^{\delta^{x}_{\varepsilon} u}_{\varepsilon^{-1}} 
\, Q^{\delta^{x}_{\varepsilon} u}_{\varepsilon^{-1}} \, \delta^{x}_{\varepsilon} v \ =
\  \bar{\delta}^{\delta^{x}_{\varepsilon} u}_{\varepsilon^{-1}} \, \bar{\delta}^{x}_{\varepsilon} 
\, \Theta_{\varepsilon}^{x}(u,v) \ = \ \bar{\Delta}^{x}_{\varepsilon}(Q^{x}_{\varepsilon} u, 
\Theta_{\varepsilon}^{x}(u,v))$$
We can pass to the limit in the last term of this string of equalities and 
we prove that   the axiom A4 is satisfied by $\delta$: there exists the limit 
\begin{equation}
\Delta^{x}(u,v) \, = \, \lim_{\varepsilon \rightarrow 0} \Delta_{\varepsilon}^{x}(u,v) 
\label{defstu}
\end{equation}
which is uniform as written in A4, 
  moreover  we have the equality (\ref{neednext1}).

(d) We pass to the limit  to the limit with $\varepsilon \rightarrow 0$ in the  relation
(\ref{neednext1}) and we obtain (\ref{recon}). In order to prove 
(\ref{morph}) we notice that: 
$$Q^{\delta^{x}_{\varepsilon} u}_{\varepsilon} \, 
\Delta^{x}_{\varepsilon}(u,v) \, = \, 
 Q^{\delta^{x}_{\varepsilon} u}_{\varepsilon} \delta^{\delta^{x}_{\varepsilon}
 u}_{\varepsilon^{-1}} \delta^{x}_{\varepsilon} v \, =   \,  \bar{\delta}^{\delta^{x}_{\varepsilon} u}_{\varepsilon^{-1}}
 \bar{\delta}^{x}_{\varepsilon}  
 Q^{x}_{\varepsilon} v \, = \,
 \bar{\Delta}^{x}_{\varepsilon} (Q^{x}_{\varepsilon} u,
 Q^{x}_{\varepsilon} v)$$ 
which gives(\ref{morph}) by  passing to the limit with $\varepsilon \rightarrow 0$ in this relation.  \quad $\square$

\paragraph{Induced dilations at a scale.} 
For any $x \in X$ and $\varepsilon \in (0,1)$ the  dilation 
$\displaystyle \delta^{x}_{\varepsilon}$ can be seen as an isomorphism of 
strong dilation structures with coherent projections: 
$$\delta^{x}_{\varepsilon} : (U(x), \delta^{x}_{\varepsilon} \bar{d} ,
\hat{\delta}^{x}_{\varepsilon}, \hat{Q}^{x}_{\varepsilon}) \rightarrow 
(\delta^{x}_{\varepsilon} U(x) , \frac{1}{\varepsilon} \bar{d} , \bar{\delta}, Q)$$
We may use this morphism in order to transport the induced dilations and the coherent 
projection: 
$$\hat{\delta}^{x, u}_{\varepsilon, \mu} \ = \ \delta^{x}_{\varepsilon^{-1}} \, 
\bar{\delta}_{\mu}^{\delta^{x}_{\varepsilon} u} \, \delta^{x}_{\varepsilon} 
\quad , \quad \quad \hat{Q}^{x, u}_{\varepsilon, \mu} \ = \ \delta^{x}_{\varepsilon^{-1}} \, 
Q_{\mu}^{\delta^{x}_{\varepsilon} u} \, \delta^{x}_{\varepsilon} 
 $$
The dilation $\displaystyle \bar{\delta}^{x}_{\varepsilon}$, which is 
 an isomorphism of strong dilation structures with coherent projections: 
$$\bar{\delta}^{x}_{\varepsilon} : (U(x), \bar{\delta}^{x}_{\varepsilon} \bar{d},
\bar{\delta}^{x}_{\varepsilon}, \bar{Q}^{x}_{\varepsilon}) \rightarrow 
(\bar{\delta}^{x}_{\varepsilon} U(x) , \frac{1}{\varepsilon} \bar{d} , 
\bar{\delta}, Q)$$
may be also used to transport induced dilations and coherent projections: 
$$\bar{\delta}^{x, u}_{\varepsilon, \mu} \ = \ \bar{\delta}^{x}_{\varepsilon^{-1}} \, 
\bar{\delta}_{\mu}^{\bar{\delta}^{x}_{\varepsilon} u} \, \bar{\delta}^{x}_{\varepsilon} 
 \quad , \quad \quad  \bar{Q}^{x, u}_{\varepsilon, \mu} \ = \ \bar{\delta}^{x}_{\varepsilon^{-1}} \, 
Q_{\mu}^{\bar{\delta}^{x}_{\varepsilon} u} \, \bar{\delta}^{x}_{\varepsilon} 
 $$
These two morphisms are related: because $\displaystyle \delta^{x}_{\varepsilon} \, = \,
\bar{\delta}^{x}_{\varepsilon} \, Q^{x}_{\varepsilon}$ it follows that  
$$Q^{x}_{\varepsilon}: (U(x), \delta^{x}_{\varepsilon} \bar{d} ,
\hat{\delta}^{x}_{\varepsilon}, \hat{Q}^{x}_{\varepsilon}) \rightarrow 
(Q^{x}_{\varepsilon} U(x), \bar{\delta}^{x}_{\varepsilon} d ,
\bar{\delta}^{x}_{\varepsilon}, \bar{Q}^{x}_{\varepsilon})$$ 
 is yet  another isomorphism of strong dilation structures with coherent projections. This isomorphism measures the "distortion" between the previous two isomorphisms. Recall that in the limit $\displaystyle Q^{x}_{\varepsilon}$ becomes a projection, therefore 
$\displaystyle Q^{x}_{\varepsilon} U(x)$ will be squeezed to a flattened version $\displaystyle 
Q^{x} U(x)$ (which will represent, in the case of sub-riemannian manifolds, a local version of the distribution evaluated at $x$).

 We shall denote the derivative of a curve with respect to the dilations 
$\displaystyle \hat{\delta}^{x}_{\varepsilon}$ by $\displaystyle 
\frac{\hat{d}^{x}_{\varepsilon}}{dt}$. Also, the derivative 
of the curve $c$ with respect to $\bar{\delta}$ is denoted by 
$\displaystyle \frac{\bar{d}}{dt}$. A curve $c$ is 
$\displaystyle \hat{\delta}^{x}_{\varepsilon}$-derivable if and only if 
$\displaystyle \delta^{x}_{\varepsilon} c$ is $\bar{\delta}$-derivable and 
$$ \frac{\hat{d}^{x}_{\varepsilon}}{dt} \, c(t) \ = \
\delta^{x}_{\varepsilon^{-1}} \, \frac{\bar{d}}{dt} \left(
\delta^{x}_{\varepsilon} c \right) (t) $$

\paragraph{$Q$-horizontal curves.} These  will play an important role further, here is the definition. 

\begin{definition}
Let $(X,\bar{d}, \bar{\delta})$ be a strong dilation structure and $Q$ a 
coherent projection. A curve $c: [a,b] \rightarrow X$ is $Q$- horizontal if 
for almost any $t \in [a,b]$ the curve $c$ is derivable and the derivative 
of $c$ at $t$, denoted by $\dot{c}(t)$ has the property:
\begin{equation}
 Q^{c(t)} \dot{c}(t) \, = \, \dot{c}(t)
 \label{horprop}
 \end{equation}
 A curve $c: [a,b] \rightarrow X$ is $Q$- everywhere horizontal if for all 
 $t \in [a,b]$ the curve $c$ is derivable and the derivative has the
 horizontality property (\ref{horprop}).
\end{definition}

 If the curve $\displaystyle \delta^{x}_{\varepsilon} c$ is $Q$-horizontal then 
 $\displaystyle \frac{\bar{d}^{x}_{\varepsilon}}{dt} \left( Q^{x}_{\varepsilon} c \right) (t) \, =  \, \Theta^{x}_{\varepsilon}(c(t), \frac{\hat{d}^{x}_{\varepsilon}}{dt}  c  (t))$. Indeed,  
 $$\frac{\bar{d}^{x}_{\varepsilon}}{dt} \left( Q^{x}_{\varepsilon} c \right) (t) 
 \, = \, \bar{\delta}^{x}_{\varepsilon^{-1}} \, Q^{\delta^{x}_{\varepsilon} c(t)} \,
 \bar{\delta}^{x}_{\varepsilon} \, \frac{\bar{d}^{x}_{\varepsilon}}{dt} \left(
 Q^{x}_{\varepsilon} c \right) (t)$$
 which implies: 
 $\displaystyle \bar{\delta}^{x}_{\varepsilon} \, 
 \frac{\bar{d}^{x}_{\varepsilon}}{dt} \left( Q^{x}_{\varepsilon} c \right) (t) \, =
 \, Q^{\delta^{x}_{\varepsilon} c(t)}_{\varepsilon^{-1}} \, \bar{\delta}^{x}_{\varepsilon} \, \frac{\bar{d}^{x}_{\varepsilon}}{dt} \left(
 Q^{x}_{\varepsilon} c \right) (t) \, = \, 
 Q^{\delta^{x}_{\varepsilon} c(t)}_{\varepsilon^{-1}} \, \delta^{x}_{\varepsilon} \,
 \frac{\hat{d}^{x}_{\varepsilon}}{dt}  c  (t)$.

\subsection{Length functionals associated to coherent projections}

\begin{definition}
Let $(X,\bar{d},\bar{\delta})$  be a strong dilation structure with the Radon-Nikodym property  and $Q$ a coherent projection. We define the associated  distance $d: X \times X \rightarrow
[0, +\infty]$  by: 
$$d(x,y) \ = \ \inf \left\{ \int_{a}^{b} \bar{d}^{c(t)}(c(t), \dot{c}(t)) \mbox{ d}t \mbox{ :
} c: [a,b] \rightarrow X \mbox{ $\bar{d}$-Lipschitz },   \right.$$
$$\left. c(a) = x, c(b) = y, \mbox{ and }\forall a.e. \, \,  t \in [a,b] \quad  Q^{c(t)} \dot{c}(t) \ = \ \dot{c}(t)
\right\} $$
\end{definition}

As usual, we accept that the distance between two points may be infinite. The equivalence relation  $x \equiv y$ if $d(x,y)< + \infty$  induces adecomposition of the space 
$X$ into a  reunion of equivalence classes, each equivalence class being a set  connected by horizontal curves with finite length. 
Later on we shall give a sufficient condition (the generalized Chow condition (Cgen))
 on the coherent projection $Q$ for $X$ to be (locally) connected by horizontal curves.

If  $X$ is connected by horizontal curves and $(X,d)$ is complete then  $d$ is  a length distance. For proving this, tt is sufficient 
to check that $d$ has the approximate middle property:  for any $\varepsilon > 0$ 
and for any $x, y \in X$ there exists $z \in X$ such that 
$$\displaystyle \max \left\{ d(x,z) , d(y,z) \right\} \leq \frac{1}{2} \, d(x,y) +
\varepsilon$$

For any  $\varepsilon > 0$ there exists a horizontal curve $c: [a,b] \rightarrow X$ such that $c(a) = x$, $c(b) = y$ and $\displaystyle d(x,y) + 2 \varepsilon \geq l(c)$ (where $l(c)$ is the length of $c$ with respect to the distance $\bar{d}$). There is then a $\tau \in [a,b]$ such that 
$$\int_{a}^{\tau} \bar{d}^{c(t)}(c(t), \dot{c}(t)) \mbox{ d}t \ = \ \int_{\tau}^{b} 
\bar{d}^{c(t)}(c(t), \dot{c}(t)) \mbox{ d}t \ = \ \frac{1}{2}\, l(c) $$
Let now $z = c(\tau)$. We have then: 
$\displaystyle \max \left\{ d(x,z) , d(y,z) \right\} \ \leq \ \frac{1}{2} \, l(c) \ \leq \ 
\frac{1}{2} \, d(x,y) + \varepsilon$, which proves the claim.

{\bf Notations concerning length functionals.} The length functional associated to the distance $\bar{d}$ is denoted by 
$\bar{l}$. In the same way the length functional associated with 
$\displaystyle \bar{\delta}^{x}_{\varepsilon}$ is denoted by 
$\displaystyle \bar{l}^{x}_{\varepsilon}$. 
  
 We introduce the space $\displaystyle  \mathcal{L}_{\varepsilon}(X,d, \delta) \subset X 
 \times Lip([0,1],X,d)$: 
 $$\mathcal{L}_{\varepsilon}(X,d,  \delta) \ = \ \left\{(x ,c) \in X  \times
 \mathcal{C}([0,1],X) 
\mbox{ : } c: [0,1] \in U(x) \, \, , \,     \right.$$ 
$$\left. \delta^{x}_{\varepsilon}c \mbox{ is } 
\bar{d}-Lip, \quad Q-\mbox{horizontal} \mbox{ and } Lip(\delta^{x}_{\varepsilon}c) \leq 2 \varepsilon
l_{d}(\delta^{x}_{\varepsilon}c) 
\right\} $$

For any $\varepsilon \in (0,1)$ we define the length functional 
$$l_{\varepsilon}: \mathcal{L}_{\varepsilon}(X,d,  \delta) \rightarrow
[0,+\infty] \quad , \quad l_{\varepsilon}(x,c) \ = \ l^{x}_{\varepsilon}(c) \ = \ \frac{1}{\varepsilon}
\, 
\bar{l}(\delta^{x}_{\varepsilon} c)  $$

By theorem \ref{fleng} we have: 
$$l^{x}_{\varepsilon}(c) \ = \ \int^{1}_{0} \, \frac{1}{\varepsilon} \,
\bar{d}^{\delta^{x}_{\varepsilon} c(t)} \left( \delta^{x}_{\varepsilon} c(t), 
\frac{\bar{d}}{dt} \left(
\delta^{x}_{\varepsilon} c \right) (t)\right) \mbox{ d}t \ =  \ \int^{1}_{0} \, \frac{1}{\varepsilon} \,
\bar{d}^{\delta^{x}_{\varepsilon} c(t)} \left( \delta^{x}_{\varepsilon} c(t), 
\delta^{x}_{\varepsilon} \, \frac{\hat{d}^{x}_{\varepsilon}}{dt} \, c(t)
\right) \mbox{ d}t $$

\subsection{Conditions (A) and (B)}

Further are two  supplementary hypotheses on  a coherent projection $Q$.

\begin{definition}
Let $\displaystyle (X, \bar{d}, \bar{\delta})$ be a strong dilation structure, $Q$ a coherent projection and $\delta$ the induced dilation. 
\begin{enumerate}
\item[(A)] $\displaystyle \delta^{x}_{\varepsilon}$ is $\bar{d}$-bilipschitz 
in compact sets in the following sense: for any compact set $K \subset X$ and 
for any $\varepsilon \in (0,1]$ there is a number $L(K) > 0$ such that for any 
$x \in K$ and $u,v$ sufficiently close to $x$ we have: 
$$\frac{1}{\varepsilon} \, \bar{d} \left(\delta^{x}_{\varepsilon} u ,
\delta^{x}_{\varepsilon} v \right) \, \leq \, L(K) \, \bar{d}(u,v)$$  
\item[(B)] if $\displaystyle u = Q^{x} u$ then the curve $\displaystyle 
t \in [0,1] \mapsto \, Q^{x}\, \delta^{x}_{t} \, u \, = \, \bar{\delta}^{x}_{t}
u \, = \, \delta^{x}_{t} u$ is $Q$-everywhere horizontal and for any $a \in [0,1]$ we have 
$$  \limsup_{a \rightarrow 0} \frac{\bar{l}\left( t \in [0,a] \mapsto
 \bar{\delta}^{x}_{t} u \right)}{\bar{d}(x, \bar{\delta}^{x}_{a} u)} \, =   \, 1$$ 
 uniformly with respect to  $x,u$ in compact set $K$. 
 \end{enumerate}
\label{suppdef}
 \end{definition}

 The condition (A)  implies that  
 for any $\bar{d}$-Lipschitz curve $c$, the "rescaled" curve $\displaystyle
 \delta^{x}_{\varepsilon}c$ is also Lipschitz. Moreover, the condition (A)  is equivalent with the fact that  
 $\displaystyle Q^{x}_{\varepsilon}$ is locally $\displaystyle \bar{\delta}^{x}_{\varepsilon} \bar{d}$-Lipschitz, where we use the notation 
$$\left(\bar{\delta}^{x}_{\varepsilon} \bar{d} \right) (u, v) \, = \, 
\frac{1}{\varepsilon} \, \bar{d} \left( \bar{\delta}^{x}_{\varepsilon} u , 
\bar{\delta}^{x}_{\varepsilon} v \right)$$  
 More precisely, condition (A) is equivalent with:  for any compact set $K \subset X$ and 
for any $\varepsilon \in (0,1]$ there is a number $L(K) > 0$ such that for any 
$x \in K$ and $u,v$ sufficiently close to $x$ we have:
\begin{equation}
\left(\bar{\delta}^{x}_{\varepsilon} \bar{d}\right) \, 
\left(Q^{x}_{\varepsilon} u , Q^{x}_{\varepsilon} v \right) \, \leq \, L(K) \, \bar{d}(u,v)
\label{uu}
\end{equation}
Indeed, we have: 
$$\left(\bar{\delta}^{x}_{\varepsilon} \bar{d}\right) 
\left(Q^{x}_{\varepsilon} u , Q^{x}_{\varepsilon} v \right) \, = \, 
\frac{1}{\varepsilon} \, \bar{d} \left(\delta^{x}_{\varepsilon} u ,
\delta^{x}_{\varepsilon} v \right) \, \leq \, L(K) \, \bar{d}(u,v)$$

The  condition (B) is a generalization of the fact that   the
"distribution" $\displaystyle x \mapsto Q^{x} U(x)$ is  generated by horizontal  one 
parameter flows, see  theorem \ref{rkc}.

\section{Distributions in sub-riemannian spaces as coherent projections}

Here we explain how coherent projections appear in  sub-riemannian
geometry.

Let $\displaystyle \left\{ Y_{1}, ..., Y_{n} \right\}$ be a frame induced by 
a parameterization $\displaystyle \phi: O \subset \mathbb{R}^{n} \rightarrow U
\subset M$ of a small open, connected set $U$ in the manifold $M$. This
parameterization induces a dilation structure on $U$, by
$$\tilde{\delta}^{\phi(a)}_{\varepsilon} \, \phi(b) \ = \ \phi\left( a +
\varepsilon(-a+b)\right) $$
We take the distance $\tilde{d}(\phi(a), \phi(b)) \, = \, \|b-a\|$. 

Let $\displaystyle \left\{ X_{1}, ..., X_{n} \right\}$ be a normal frame, 
cf. definition \ref{defnormal}, let $d$ be the Carnot-Carath\'eodory distance 
and let 
$$\delta^{x}_{\varepsilon} \left(\exp\left( \sum_{i=1}^{n} a_{i} X_{i} \right)(x)\right) \  = \
\exp\left( \sum_{i=1}^{n} a_{i} \varepsilon^{deg X_{i}}  X_{i} \right)(x)$$
be the dilation structure associated, by theorem \ref{structhm}. 

Alternatively, we may take another dilation structure, constructed as follows: extend the
metric $g$ on the distribution $D$ to a riemannian metric on $M$, denoted 
for convenience also by $g$. Let $\displaystyle \bar{d}$ be the riemannian
distance induced by the riemannian metric $g$, and the dilations 
$$\bar{\delta}^{x}_{\varepsilon} \left(\exp\left( \sum_{i=1}^{n} a_{i} X_{i} \right)(x)\right) \  = \
\exp\left( \sum_{i=1}^{n} a_{i} \varepsilon  X_{i} \right)(x)$$
Then $\displaystyle (U, \bar{d}, \bar{\delta})$ is a strong dilation 
structure which is equivalent with the dilation structure 
$\displaystyle (U, \tilde{d}, \tilde{\delta})$. 

From now on  we may define coherent projections associated either to the pair 
$\displaystyle (\tilde{\delta}, \delta)$ or to the pair $\displaystyle 
(\bar{\delta}, \delta)$. 

 Let us define $\displaystyle Q^{x}_{\varepsilon}$ by: 
\begin{equation}
Q^{x}_{\varepsilon}\,  \left(\exp\left( \sum_{i=1}^{n} a_{i} X_{i} \right)(x)\right) \  = \
\exp\left( \sum_{i=1}^{n} a_{i} \varepsilon^{deg X_{i} - 1}  X_{i} \right)(x)
\label{eqsr1}
\end{equation}

\begin{theorem}
 $Q$ is a coherent projection associated with the dilation structure
$\displaystyle (U, \bar{d}, \bar{\delta})$ which satisfies conditions (A) and (B) definition \ref{suppdef}. Moreover,  for any point $x$ the image of 
$Q^{x}$ is in the "distribution" $D(T_{x}(U,d,\delta))$ (from proposition \ref{pintrinsicd}), where $(U,d,\delta)$ is the strong dilation structure 
from theorem \ref{structhm}. 
\label{rkc}
\end{theorem}

\paragraph{Proof.}(I) definition \ref{defcoh} is true, because 
$\displaystyle \delta^{x}_{\varepsilon} \, u \, = \, Q^{x}_{\varepsilon} \, 
\bar{\delta}^{x}_{\varepsilon}$ and $\displaystyle \delta^{x}_{\varepsilon} \, 
\bar{\delta}^{x}_{\varepsilon} \, = \, \bar{\delta}^{x}_{\varepsilon}
\delta^{x}_{\varepsilon}$. (II), (III) and (IV) are consequences of these facts 
and theorem \ref{structhm}, with a proof similar to the one of proposition 
\ref{p1proj}. 

Definition (\ref{eqsr1}) of the coherent projection $Q$ implies that: 
\begin{equation}
Q^{x}\,  \left(\exp\left( \sum_{i=1}^{n} a_{i} X_{i} \right)(x)\right) \  = \
\exp\left( \sum_{deg X_{i} = 1} a_{i}   X_{i} \right)(x)
\label{eqsr2}
\end{equation}
The projection $\displaystyle Q^{x}$ has the property: for any $x$ and 
$$\displaystyle u \, = \, \exp\left( \sum_{deg X_{i} = 1} a_{i} X_{i}
\right)(x)$$ 
we have $\displaystyle Q^{x} u \, = \, u$ and the curve 
$$s \in [0,1] \mapsto \delta^{x}_{s} \, u\ = \ \exp\left( s \, \sum_{deg X_{i} =
1} a_{i}   X_{i} \right)(x) $$
is $D$-horizontal and joins $x$ and $u$. This implies condition (B). As for the condition (A), according to the comments after definition \ref{suppdef}, it just means that the coherent projection given by (\ref{eqsr1}) is locally Lipschitz (with a Lipschitz constant uniform with respect to $x$ in compact set) with respect to the rescaled riemannian distance 
$\displaystyle \bar{\delta}^{x}_{\varepsilon} \bar{d}$. But the dilation structure $\displaystyle (U, \bar{d}, \bar{\delta})$ is tempered, according to theorem \ref{tchariem}, therefore the rescaled distance is bilipshitz with the riemannian distance $\displaystyle \bar{d}$, again with Lipschitz constants which are uniform with respect to $x$ in compact set. From (\ref{eqsr1}) we easily get that $\displaystyle Q^{x}_{\varepsilon}$ is (uniformly w.r.t. $x$) Lipschitz   w.r.t. the riemannian distance $\displaystyle \bar{d}$. In conclusion (A) is true.

After examination of the construction of the dilation structure $(U,d,\delta)$ from theorem 
\ref{structhm}, we see that  $D(T_{x}(U,d,\delta))$ (from proposition \ref{pintrinsicd}) is exactly the set of all points (sufficiently close to $x$) which can be written as exponentials of vectors in the (classical differential) geometric distribution. Therefore the last statement of the theorem is proved. \quad $\square$

We may equally define a coherent projection which induces the dilations 
$\delta$ from $\tilde{\delta}$.  Also, if we change the chosen normal frame with 
another of the same kind, then we  pass to a dilation structure which is 
equivalent to $(U,d, \delta)$. In conclusion, coherent
projections are not geometrical objects per se, but in a natural way one may
define a notion of equivalent coherent projections such that the equivalence
class is geometrical, i.e. independent of the choice of a pair 
of particular dilation structures, each in a given equivalence class. Another
way of putting this is that a class of equivalent dilation structures may be
seen as a category and a coherent projection is a functor between such
categories. We shall not pursue this line here. Anyway,  the only advantage of choosing 
$\displaystyle \bar{\delta}, \delta$ related by the normal frame 
$\displaystyle \left\{ X_{1}, ..., X_{n} \right\}$ is that they are associated 
with a coherent projection with a simple expression.

%%%%%%%%%%%%%%%%%%%%%%%%%%%%%%%%%%%%%%%%%%%%%%%%%%%%%%%%%%%%%%%%%%%%%%%%%%%%%%%%%%%%%%%%
%%%%%%
%%%%%%                               AICI AM RAMAS 
%%%%%%
%%%%%%%%%%%%%%%%%%%%%%%%%%%%%%%%%%%%%%%%%%%%%%%%%%%%%%%%%%%%%%%%%%%%%%%%%%%%%%%%%%%%%%%%

\section{An intrinsic description of sub-riemannian geometry}

\subsection{The generalized Chow condition}
\label{sechormander}

We want to transform  the Chow condition (C), theorem \ref{tchow}, into a statement formulated in terms of coherent projections. Essentially, the Chow condition (C) tells us that the (sub-riemannian) space is locally connected by horizontal curves which are constructed from concatenations of exponentials of horizontal vector fields. 

In the following we need to explain what are the correspondents of exponentials of horizontal vector fields, then we need a way to manage the concatenation procedure. A simple explanation is this: 
\begin{enumerate}
\item[(a)] the exponentials of horizontal vector fields are, approximately, the induced dilations of the coherent projection, 
\item[(b)] A concatenation of those exponentials is coded by a word made by letters, each letters coding one exponential. 
\end{enumerate}

A supplementary complication is that we need to have an efficient way to manage these concatenations of exponentials {\em at any scale}.That is why we start with a scale, that is with a $\varepsilon > 0$ and with the coherent projection, dilations, and so on, at that scale (i.e. those transported by some properly chosen dilation field). 

We shall follow closely \cite{buligachar} section 10, getting into details as necessary. 

{\bf Words over an alphabet.} This is a standard way of notation. We shall need words as codes for curves, as explained previously. For any "alphabet" (that is a set) $A$ we denote by 
$\displaystyle A^{*}$ the collection of finite words $\displaystyle q = a_{1}...a_{p}$, $\displaystyle p \in \mathbb{N}$, $p > 0$. The empty word (with no letters) is denoted 
by $\emptyset$. The length of the word $\displaystyle q = a_{1}...a_{p}$ 
is $\mid q \mid = p$; the length of the empty word is $0$. We may need  words  which are infinite at right. The set of those words over the alphabet $A$ is denoted by 
$\displaystyle A^{\omega}$. For any word $\displaystyle w \in A^{\omega} \cup
A^{*}$ and
any $p \in \mathbb{N}$ we denote by $\displaystyle [w]_{p}$ the finite word 
obtained from the first $p$ letters of $w$ (if $p=0$ then $\displaystyle 
[w]_{0} = \emptyset$ (in the case of a finite word $q$, if $p > \mid q \mid$ 
then $\displaystyle [q]_{p} = q$).

For any non-empty $\displaystyle q_{1}, q_{2} \in A^{*}$ and $w \in A^{\omega}$ the 
concatenation  of $\displaystyle q_{1}$ and $\displaystyle q_{2}$ is the 
finite word $\displaystyle q_{1}q_{2} \in A^{*}$ and the concatenation of 
$\displaystyle q_{1}$ and $w$ is the (infinite) word $\displaystyle q_{1}w \in 
A^{\omega}$. For any  $\displaystyle q \in A^{*}$ and 
$\displaystyle w \in A^{\omega}$ we have 
$\displaystyle q \emptyset = q$ (as concatenation of finite words) and 
$\emptyset w = w$ (as concatenation of a finite empty word and an infinite 
word).

\paragraph{Words as codes for $Q$-horizontal curves.} 
To the coherent projection $Q$ of a strong dilation structure $\displaystyle 
(X, \bar{d}, \bar{\delta})$ we associate a family of transformations, which correspond to concatenations of exponentials of horizontal vector fields, at a given scale. 

\begin{definition}
To any  word $\displaystyle w \in (0,1]^{\omega}$ and any 
$\varepsilon \in (0,1]$ we associate  the transformation 
$$\Psi_{\varepsilon w} : X^{*}_{\varepsilon w} \, \subset X^{*} \setminus
\left\{ \emptyset \right\} \, 
\rightarrow X^{*} $$  
defined recursively by the following procedure. 

If $w = \emptyset$ then we define  $\displaystyle \Psi_{\varepsilon \emptyset}^{1}(x) = x$.  For any non-empty finite word $\displaystyle q = x x_{1} ... x_{p} \in 
X^{*}_{\varepsilon w}$ and  for any $k \geq 1$ we 
 have 
$$\Psi_{\varepsilon \emptyset}^{k+1}([q]_{k+1}) \ = \
\delta^{x}_{\varepsilon^{-1}} \, Q^{\delta^{x}_{\varepsilon} \,
\Psi_{\varepsilon w}^{k}([q]_{k})} \, \delta_{\varepsilon}^{x} \, q_{k+1} $$

If $w$ is not the empty word then we define the functions $\displaystyle \Psi_{\varepsilon w}^{k}$  by: 
 $\displaystyle \Psi_{\varepsilon w}^{1}(x) = x$, and for any $k \geq 1$ we 
 have 
 \begin{equation}
\Psi_{\varepsilon w}^{k+1}([q]_{k+1}) \ = \
\delta^{x}_{\varepsilon^{-1}} \, Q_{w_{k}}^{\delta^{x}_{\varepsilon} \,
\Psi_{\varepsilon w}^{k}([q]_{k})} \, \delta_{\varepsilon}^{x} \, q_{k+1} 
\label{recrel}
\end{equation}

Finally, for any non-empty finite word $\displaystyle q = x x_{1} ... x_{p} \in 
X^{*}_{\varepsilon w}$ we put 
$$ \Psi_{\varepsilon w} (x x_{1} ... x_{p}) \ = \ 
\Psi_{\varepsilon w}^{1}(x) ... \Psi_{\varepsilon w}^{k+1}(x x_{1} ... x_{k}) ... 
\Psi_{\varepsilon w}^{p+1}(x x_{1} ... x_{p}) $$

The domain $\displaystyle X^{*}_{\varepsilon w} \subset X^{*} \setminus
\left\{ \emptyset \right\}$ is  such that the previous definitions make  sense.
\label{dwords}
\end{definition}

 By using the definition of a coherent projection, we can give the following, more precise description of$\displaystyle
X^{*}_{\varepsilon w}$ as  follows:  for any compact set $K \subset X$  there 
is $\rho = \rho(K) > 0$ such that for any $x \in K$ the word 
$\displaystyle q = x x_{1} ... x_{p} \, \in X^{*}_{\varepsilon w}$ if for any 
$k \geq 1$ we have 
$$\bar{d}\left( x_{k+1} , \Psi_{\varepsilon w}^{k}([q]_{k}) \right) \, \leq \,
\rho $$

Let us see what this gives for $\varepsilon = 1$. By  definition \ref{dwords} for 
$\varepsilon = 1$  we have: 
$$\Psi_{1w}^{1}(x) =  x \quad , \quad \Psi_{1w}^{2}(x, x_{1}) \ =  \
Q_{w_{1}}^{x} \, x_{1} \quad , \quad \Psi_{1w}^{3}(x, x_{1}, x_{2}) \ = \ 
Q_{w_{2}}^{Q_{w_{1}}^{x} x_{1}} \, x_{2} \quad ... $$ 

\begin{proposition}
Suppose that the coherent projection $Q$ satisfies the condition (B) and let $y \in X$ be  
$$\displaystyle y \ = \ \Psi_{1\emptyset}^{k+1}(x x_{1} ... x_{k})$$ 
Then there is a $Q$-horizontal curve joining $x$ and $y$.
\end{proposition}

\paragraph{Proof.} 
By  condition (B)  the
curve $\displaystyle t \in [0,1] \mapsto \bar{\delta}^{x}_{t} Q^{x} u$ is a 
$Q$-horizontal curve which joins joining $x$ with $\displaystyle Q^{x} u$. 

Therefore 
 by applying repeatedly the condition (B) we get that 
 there is a $Q$-horizontal curve between 
$\displaystyle \Psi_{1\emptyset}^{k}(x x_{1} ... x_{k-1})$ and $\displaystyle
\Psi_{1\emptyset}^{k+1}(x x_{1} ... x_{k})$ for any $k > 1$ and a $Q$-horizontal
curve joining $x$ and $\displaystyle \Psi_{1\emptyset}^{2}(x x_{1})$. Therefore by concatenation we get the desired curve. 
\quad $\square$

The following are algebraic properties  of the transformations $\displaystyle
\Psi_{\varepsilon w}$, which explain how they behave with respect to scale. 

\begin{proposition}
With the notations from definition \ref{dwords} we have: 
\begin{enumerate}
\item[(a)] $\displaystyle \Psi_{\varepsilon w} \, 
\Psi_{\varepsilon \emptyset} \, = \, \Psi_{\varepsilon \emptyset}$. Therefore 
we have the equality of sets: 
$$\Psi_{\varepsilon \emptyset}\, \left( X^{*}_{\varepsilon \emptyset} \cap x X^{*}
\right) \, = \, \Psi_{\varepsilon w}\, \left( \Psi_{\varepsilon \emptyset}\,
\left( X^{*}_{\varepsilon \emptyset} \cap x X^{*}
\right) \right)$$
\item[(b)] $\displaystyle \Psi_{\varepsilon \emptyset}^{k+1} (x q_{1} ... q_{k})
\, = \, \delta_{\varepsilon^{-1}}^{x} \, \Psi_{1 \emptyset}^{k+1} (x
\delta_{\varepsilon}^{x}q_{1} ... \delta_{\varepsilon}^{x} q_{k})$
\item[(c)] $\displaystyle \lim_{\varepsilon \rightarrow 0} 
 \delta_{\varepsilon^{-1}}^{x} \, \Psi_{1 \emptyset}^{k+1} (x
\delta_{\varepsilon}^{x}q_{1} ... \delta_{\varepsilon}^{x} q_{k}) \, = \, 
\Psi_{0 \emptyset}^{k+1} (x q_{1} ... q_{k})$ uniformly with respect to 
$\displaystyle x, q_{1}, ..., q_{k}$ in compact set. 
\end{enumerate}
\end{proposition}

\paragraph{Proof.}
(a) We use induction on $k$ to prove that for any natural number $k$ we have: 
\begin{equation}
\Psi^{k+1}_{\varepsilon w} \left( \Psi^{1}_{\varepsilon \emptyset}(x) ... 
\Psi^{k+1}_{\varepsilon \emptyset} (x q_{1} ... q_{k}) \right) \, = \, 
\Psi^{k+1}_{\varepsilon \emptyset} (x q_{1} ... q_{k})
\label{prpaproof}
\end{equation}
 For $k = 0$ we have have to prove that  $x = x$ which is trivial. For 
 $k = 1$ we have to prove that 
 $$\Psi^{2}_{\varepsilon w} \left( \Psi^{1}_{\varepsilon \emptyset}(x) \,  
\Psi^{2}_{\varepsilon \emptyset} (x q_{1}) \right) \, = \, 
\Psi^{2}_{\varepsilon \emptyset} (x q_{1})$$ 
This means: 
$$\Psi^{2}_{\varepsilon w} \left( x \,  \delta_{\varepsilon^{-1}}^{x} \, 
Q^{x} \delta^{x}_{\varepsilon} q_{1} \right) \, = \, 
\delta_{\varepsilon^{-1}}^{x} \, Q_{w_{1}}^{x} \, \delta_{\varepsilon}^{x} \, 
\delta^{x}_{\varepsilon^{-1}} \, Q^{x} \, \delta^{x}_{\varepsilon} x_{1} \, = 
 \, \delta^{x}_{\varepsilon^{-1}} \, Q^{x} \, \delta^{x}_{\varepsilon} x_{1} \, =
\, \Psi^{2}_{\varepsilon \emptyset} (x q_{1})$$

  Suppose now that $l \geq 2$ and  for any $k \leq l$ the relations 
(\ref{prpaproof}) are true. Then, as previously, it is easy to check 
(\ref{prpaproof}) for $k = l+1$.

(b) is true by direct computation. The point (c) is a straightforward
consequence of (b) and definition of coherent projections.  
\quad $\square$

The following definition is quantitative: it says that a point $x$ is "nested" in a neighbourhood $U$ with respect to the parameters $N$ (natural number) and $\varepsilon$ (scale) if there is a small ball of radius $\rho$ around $x$, with respect to the scaled distance 
$\displaystyle \left(\bar{\delta}^{x}_{\varepsilon} \bar{d}\right) $, such that: the ball is inside $U$ and we can connect  $x$ with any $y$ in the ball, by using at most $N$ concatenations of curves coming from the transformations introduced in definition \ref{dwords}, such that all these curves stay in $U$. 

\begin{definition}
Let $N \in \mathbb{N}$ be a strictly positive natural number and $\varepsilon
\in (0,1]$. 
We say that $x \in X$ is  $(\varepsilon,N,Q)$-nested in a open neighbourhood 
$U \subset X$ if there is $\rho> 0$ such that for any finite word $\displaystyle 
q = x_{1} ... x_{N} \in X^{N}$, if  
$\displaystyle \left(\bar{\delta}^{x}_{\varepsilon} \bar{d}\right) \left( x_{k+1} , \Psi_{\varepsilon
 \emptyset}^{k}([xq]_{k}) \right) \, \leq \, \rho $ 
for any $k = 1, ... , N$ then we have $\displaystyle q \in U^{N}$. 

If $x \in U$ is $(\varepsilon,N,Q)$-nested then denote by 
$\displaystyle U(x,\varepsilon, N,Q,\rho) \subset U^{N}$ the collection of 
words $\displaystyle q \in U^{N}$ such that 
$\displaystyle \bar{\delta}^{x}_{\varepsilon} \bar{d}\left( x_{k+1} ,
\Psi_{\varepsilon \emptyset}^{k}([xq]_{k}) \right)
\, < \,
\rho $ for any $k = 1, ... , N$. 
\end{definition}

In the next definition we introduce the condition  (Cgen).  Its effect is that if 
the coherent projection $Q$ satisfies  (A) and (B) then   in the
space $\displaystyle (U(x), \bar{\delta}^{x}_{\varepsilon})$, with 
coherent projection $\displaystyle \hat{Q}^{x, \cdot}_{\varepsilon.\cdot}$, 
 we can join any two sufficiently close points by a sequence of 
at most $N$ horizontal curves. Moreover
there is a control on the length of these curves via condition (B) and condition
(Cgen). 

The function $F$ which appears in the next definition is not specified in general. In the case of  sub-riemannian geometry the function $F$ can be taken as  
$\displaystyle   F(\eta) = C \eta^{1/m}$ with $m$ positive natural
number (the step of the distribution). Compare with the Folland-Stein lemma \ref{fp2.4}.

\begin{definition}
A coherent projection $Q$  satisfies  the  generalized Chow condition  
if: 
\begin{enumerate}
 \item[(Cgen)] for any compact set $K$ there are $\rho = \rho(K) > 0$, $r = r(K) > 0$, 
 a natural number $N = N(Q,K)$ and a function 
 $F(\eta) = \mathcal O(\eta)$ such that for
 any $x \in K$ and $\varepsilon \in (0,1]$ 
 there are   neighbourhoods $U(x)$, $V(x)$ such that  any 
  $x \in K$ is $(\varepsilon, N,Q)$-nested in $U(x)$, 
  $\displaystyle B(x,r, \bar{\delta}^{x}_{\varepsilon} \bar{d}) 
  \subset V(x)$  and such that the mapping 
 $$x_{1} ... x_{N} \in U(x,N,Q,\rho) \, \mapsto \,  
 \Psi^{N+1}_{\varepsilon \emptyset}(x x_{1} ... x_{N})$$ 
 is surjective from $\displaystyle U(x,\varepsilon,N,Q,\rho)$ to 
 $\displaystyle V(x)$. 
 Moreover for any $\displaystyle z \in V(x)$ there exist $\displaystyle  
 y_{1}, ...  y_{N} \in U(x,\varepsilon, N, Q, \rho)$ such that 
 $\displaystyle z \, = \, \Psi^{N+1}_{\varepsilon \emptyset}(x y_{1}, ...  y_{N})$ and 
 for any $k = 0, ... , N-1$ we have 
 $$\delta^{x}_{\varepsilon} \bar{d} \left( \Psi^{k+1}_{\varepsilon
 \emptyset}(x y_{1} ...  y_{k}) , 
 \Psi^{k+2}_{\varepsilon \emptyset}(x y_{1} ...  y_{k+1}) \right) \ \leq \
 F(\delta^{x}_{\varepsilon} \bar{d}(x,z))$$
 \end{enumerate}

Suppose now that the coherent projection $Q$ satisfies conditions (A), (B) and (Cgen) and 
let us consider $\varepsilon \in (0,1]$ and 
 $x, y \in K$, $K$ compact in $X$. Then   there are numbers $N=N(Q,K)$, 
$\rho = \rho(Q,K) > 0$  and  words 
$\displaystyle  x_{1} ... x_{N} \in U(x,\varepsilon,N,Q,\rho)$
  such that  $\displaystyle  y = \Psi_{\varepsilon\emptyset}^{N+1}(x x_{1} ... x_{N})$. 
  To these data we associate a  "short curve"  $c : [0,N] \rightarrow X$ ,  which joins 
  $x$ and $y$,  by: 
 for any  $t \in [0, N]$ 
$$c(t) \ = \ \bar{\delta}^{x, \Psi_{\varepsilon \emptyset}^{k+1}(x x_{1} ...
x_{k})}_{\varepsilon,t+N -k} 
Q^{\Psi_{\varepsilon\emptyset}^{k+1}(x   x_{1} ... x_{k})} x_{k+1} $$
By extension, any increasing linear reparameterization of 
a curve $c$ like the one described previously, will be called "short  curve" as well. 
\label{defhorgen} 
 \end{definition}

\subsection{The candidate tangent space}
\label{candidate}

Let $(X, \bar{d}, \bar{\delta})$ be a strong dilation structure and 
$Q$   a coherent projection.  Then we have
the induced dilations and coherent projections 
$$\mathring{\delta}^{x,u}_{\mu} v \ = \ \Sigma^{x}(u, \delta^{x}_{\mu}
\Delta^{x}(u, v)) \quad , \quad \quad \mathring{Q}^{x,u}_{\mu} v \ = \ \Sigma^{x}(u, Q^{x}_{\mu}
\Delta^{x}(u, v)) $$
For any curve $c: [0,1] \rightarrow U(x)$ which is 
$\displaystyle \mathring{\delta}^{x}$-derivable and 
$\displaystyle \mathring{Q}^{x}$-horizontal almost everywhere: 
$$\frac{\mathring{d}^{x}}{dt} c(t) \ = \ \mathring{Q}^{x,u} \, 
\frac{\mathring{d}^{x}}{dt} c(t)$$
we define the length
$$l^{x}(c) \ = \ \int^{1}_{0} \, \bar{d}^{x}\left( x, \Delta^{x}(c(t), 
\frac{\mathring{d}^{x}}{dt} c(t)) \right) \mbox{ d}t $$
and the distance function: 
$$\mathring{d}^{x}(u,v) \, = \, \inf \left\{ l^{x}(c) \mbox{ : }  c: [0,1]
\rightarrow U(x) \, \mbox{ is $\displaystyle \mathring{\delta}^{x}$-derivable} 
,\right. $$ 
$$\left. \mbox{ and $\displaystyle \mathring{Q}^{x}$-horizontal a.e.} \, , 
\, c(0) = u, c(1) = v \right\}$$

We want to prove that  $\displaystyle (U(x), \mathring{d}^{x}, \mathring{\delta}^{x})$ is a
strong dilation structure and $\displaystyle \mathring{Q}^{x}$ is a coherent projection. 
For this we need first the following theorem.

\begin{theorem}
The curve $c: [0,1] \rightarrow U(x)$ is 
$\displaystyle \mathring{\delta}^{x}$-derivable, 
$\displaystyle \mathring{Q}^{x}$-horizontal almost everywhere, and 
$\displaystyle l^{x}(c) < + \infty$ if and only if the curve  
$\displaystyle Q^{x} c$ is $\displaystyle \bar{\delta}^{x}$-derivable 
almost everywhere and 
$\displaystyle \bar{l}^{x}(Q^{x} c) < + \infty$. Moreover, we have  
$\displaystyle \bar{l}^{x}(Q^{x} c) \, = \, l^{x}(c)$. 
\label{propcan}
\end{theorem}

\paragraph{Proof.}
The curve $c$ is $\displaystyle \mathring{Q}^{x}$-horizontal almost everywhere
if and only if for almost any $t \in [0,1]$ we have 
$$Q^{x} \, \Delta^{x}(c(t) , \frac{\mathring{d}^{x}}{dt} c(t)) \ = \ 
\Delta^{x}(c(t) , \frac{\mathring{d}^{x}}{dt} c(t)) $$
We shall prove that  
 $c$ is $\displaystyle \mathring{Q}^{x}$-horizontal  is
equivalent with 
\begin{equation}
\Theta^{x}(c(t), 
\frac{\mathring{d}^{x}}{dt} c(t)) \ = \  \frac{\bar{d}^{x}}{dt} \left(Q^{x} c
\right)(t) 
\label{blakdot}
\end{equation}
Indeed, \eqref{blakdot} is equivalent with
$\displaystyle \lim_{\varepsilon \rightarrow 0} \bar{\delta}^{x}_{\varepsilon^{-1}} 
\bar{\Delta}^{x}(Q^{x} c(t), Q^{x} c(t+\varepsilon)) \ = \ 
\bar{\Delta}^{x}(Q^{x} c(t) , \Theta^{x}(c(t), 
\frac{\mathring{d}^{x}}{dt} c(t)))$,  
which is equivalent with 
$\displaystyle \lim_{\varepsilon \rightarrow 0} \bar{\delta}^{x}_{\varepsilon^{-1}} 
\bar{\Delta}^{x}(Q^{x} c(t), Q^{x} c(t+\varepsilon)) \ = \ 
\Delta^{x}(c(t), \frac{\mathring{d}^{x}}{dt} c(t))$, finally equivalent with:  
\begin{equation}
\lim_{\varepsilon \rightarrow 0} \bar{\delta}^{x}_{\varepsilon^{-1}} 
\bar{\Delta}^{x}( Q^{x} c(t), Q^{x} c(t+\varepsilon)) \ = \ 
\lim_{\varepsilon \rightarrow 0} \delta^{x}_{\varepsilon^{-1}} 
\Delta^{x}(c(t), c(t+ \varepsilon))
\label{blakdot1}
\end{equation}
The horizontality condition for the curve $c$ can be written as: 
$$\lim_{\varepsilon \rightarrow 0} Q^{x} \delta^{x}_{\varepsilon^{-1}} 
\Delta^{x}(c(t), c(t+\varepsilon)) \ = \ \lim_{\varepsilon \rightarrow 0}  
\delta^{x}_{\varepsilon^{-1}} 
\Delta^{x}(c(t), c(t+\varepsilon))$$ 
We use now the properties of $\displaystyle Q^{x}$ in the left hand side of the
previous equality: 
$$Q^{x} \delta^{x}_{\varepsilon^{-1}} 
\Delta^{x}(c(t), c(t+\varepsilon)) \ = \ 
\bar{\delta}^{x}_{\varepsilon^{-1}} Q^{x} \Delta^{x}(c(t), c(t+\varepsilon)) 
\ =  \  \bar{\delta}^{x}_{\varepsilon^{-1}} \bar{\Delta}^{x}(Q^{x} c(t), Q^{x}
c(t+\varepsilon))$$
thus after taking the limit as $\varepsilon \rightarrow 0$ we prove that 
the limit $\displaystyle \lim_{\varepsilon \rightarrow 0} 
\bar{\delta}^{x}_{\varepsilon^{-1}} \bar{\Delta}^{x}(Q^{x} c(t), Q^{x}
c(t+\varepsilon))$ exists and we obtain: 
$$\lim_{\varepsilon \rightarrow 0}  
\delta^{x}_{\varepsilon^{-1}} 
\Delta^{x}(c(t), c(t+\varepsilon)) \ = \  \lim_{\varepsilon \rightarrow 0} 
\bar{\delta}^{x}_{\varepsilon^{-1}} \bar{\Delta}^{x}(Q^{x} c(t), Q^{x}
c(t+\varepsilon))$$
This last equality is the same as \eqref{blakdot1}, which is equivalent with 
\eqref{blakdot}. 
As a consequence we obtain the following equality, for almost any $t \in [0,1]$:
\begin{equation}
\bar{d}^{x}\left( x, \Delta^{x}(c(t), 
\frac{\mathring{d}^{x}}{dt} c(t)) \right) \ = \ 
\bar{\Delta}^{x}(Q^{x} c(t) , \frac{\bar{d}^{x}}{dt} \left(Q^{x} c
\right)(t)) 
\label{blakdot3}
\end{equation}
This implies that $Q^{x}c$ is absolutely continuous and by theorem 
\ref{tupper}, as in the proof of theorem \ref{fleng} (but without using 
the Radon-Nikodym  property, because we already know that $\displaystyle Q^{x}c$ is
derivable a.e.), we obtain the following formula for the length of the curve 
$Q^{x}c$: 
$$\bar{l}^{x}(Q^{x} c) \ = \ \int^{1}_{0} \, \bar{d}^{x}\left( 
x, , \bar{\Delta}^{x}(Q^{x} c(t) , \frac{\bar{d}^{x}}{dt} \left(Q^{x} c
\right)(t)) \right) \mbox{ d}t  $$
But we have also:
$$l^{x}(c) \ = \ \int^{1}_{0} \, \bar{d}^{x}\left( x, \Delta^{x}(c(t), 
\frac{\mathring{d}^{x}}{dt} c(t)) \right) \mbox{ d}t  $$
By \eqref{blakdot3} we obtain $\displaystyle \bar{l}^{x}(Q^{x} c) \, = \,
l^{x}(c)$. 
\quad $\square$

\begin{proposition}
If $(X, \bar{d}, \bar{\delta})$ is a strong dilation structure, 
$Q$ is  a coherent projection  and $\mathring{d}^{x}$ is finite then 
 the triple  $\displaystyle (U(x), \Sigma^{x}, \delta^{x})$
is a normed conical group, with the norm induced by the  left-invariant
distance $\mathring{d}^{x}$.
\label{prevprop}
\end{proposition}

\paragraph{Proof.}
The fact that  $\displaystyle (U(x), \Sigma^{x}, \delta^{x})$ is a conical group
comes directly from the definition \ref{defcoh} of a coherent projection.
Indeed,  it is enough to use proposition \ref{p1proj} (c) and the formalism 
of binary decorated trees in \cite{buligadil1} section 4 (or theorem 11
\cite{buligadil1}), in order to  reproduce 
the part of the proof of theorem 10 (p.87-88) in that paper, concerning the
conical group structure. There is one small subtlety though. 
In the proof of theorem \ref{tgene}(a) the same modification of proof has
been done starting from the axiom A4+, namely the existence of the uniform limit
$\displaystyle \lim_{\varepsilon \rightarrow 0} \Sigma^{x}_{\varepsilon}(u,v) \,
= \, \Sigma^{x}(u,v)$. Here we need first to prove this limit, in a similar way
as in the corollary 9 \cite{buligadil1}. We shall use for this the distance 
$\displaystyle \mathring{d}^{x}$ instead of the distance 
in the metric tangent space of $(X,d)$ at $x$ denoted by 
 $\displaystyle d^{x}$ (which is not yet proven to exist). The distance 
 $\displaystyle \mathring{d}^{x}$ is supposed to be finite by hypothesis.
 Moreover, by its definition and theorem \ref{propcan} we have 
 $$\mathring{d}^{x}(u,v) \, \geq \, \bar{d}^{x}(u,v)$$
 therefore the distance $\displaystyle \mathring{d}^{x}$ is non degenerate. By 
 construction  this distance is also left invariant with respect to the 
 group operation $\displaystyle \Sigma^{x}(\cdot, \cdot)$. Therefore we may
 repeat the proof of corollary 9 \cite{buligadil1} and obtain the result that 
 A4+ is true for $(X,d,\delta)$. 

What we need
to prove next  is that $\mathring{d}^{x}$ induces a norm on the conical group 
 $\displaystyle (U(x), \Sigma^{x}, \delta^{x})$. For this it is enough to 
 prove that 
 \begin{equation}
 \mathring{d}^{x}(\mathring{\delta}^{x,u}_{\mu} v, 
 \mathring{\delta}^{x,u}_{\mu} w ) \, = \, \mu \, \mathring{d}^{x}(v,w)
 \label{mewn1}
 \end{equation}
 for any $v, w \in U(x)$. This is a direct consequence of relation 
 (\ref{blakdot3}) from the proof of the  theorem 
 \ref{propcan}. Indeed, by direct computation we get that for any curve 
 $c$ which is $\mathring{Q}^{x}$-horizontal a.e. we have: 
 $$l^{x}(\mathring{\delta}^{x,u}_{\mu} c) \, = \, 
 \int_{0}^{1} \bar{d}^{x} \left( x, \Delta^{x}\left(\mathring{\delta}^{x,u}_{\mu} c(t) , 
 \frac{\mathring{d}^{x}}{dt} \left( \mathring{\delta}^{x,u}_{\mu} c \right) (t)
 \right)\right) \mbox{ dt} \, = \, $$
$$= \, \int_{0}^{1} \bar{d}^{x} \left( x, \delta^{x}_{\mu} \Delta^{x}\left( c(t) , 
 \frac{\mathring{d}^{x}}{dt}  c (t)
 \right)\right) \mbox{ dt}$$
 But  $c$ is $\mathring{Q}^{x}$-horizontal a.e., which implies, via
 (\ref{blakdot3}),  that 
 $$\delta^{x}_{\mu} \Delta^{x}\left( c(t) , 
 \frac{\mathring{d}^{x}}{dt}  c (t) \right) \, = \, \bar{\delta}^{x}_{\mu} \Delta^{x}\left( c(t) , 
 \frac{\mathring{d}^{x}}{dt}  c (t)
 \right)$$ 
 therefore we have 
 $$l^{x}(\mathring{\delta}^{x,u}_{\mu} c) \, = \, \int_{0}^{1} \bar{d}^{x} \left( x,
 \bar{\delta}^{x}_{\mu} \Delta^{x}\left( c(t) , 
 \frac{\mathring{d}^{x}}{dt}  c (t)
 \right)\right) \mbox{ dt} \, = \, \mu \, l^{x}(c)$$ 
 This implies (\ref{mewn1}), therefore the  proof is done.
\quad $\square$

\begin{theorem}
If the generalized Chow condition (Cgen) and condition (B) are true then  
$\displaystyle (U(x), \Sigma^{x}, \delta^{x})$ is a local conical group which is 
a neighbourhood of the neutral element of a  Carnot group 
generated by $\displaystyle Q^{x} U(x)$. 
\label{tancarnot}
\end{theorem}

\paragraph{Proof.}
 For any  for any $\varepsilon \in (0,1]$, $x, u, v \in
X$ sufficiently close and $\mu > 0$ we have the relations:
\begin{enumerate}
\item[(i)] $\displaystyle \hat{Q}^{x, u}_{\varepsilon, \mu} v \, = \, 
\Sigma^{x}_{\varepsilon}(u, Q_{\mu}^{\delta^{x}_{\varepsilon} u}
\Delta^{x}_{\varepsilon}(u,v))$, 
\item[(ii)] $\displaystyle \hat{Q}^{x, u}_{\varepsilon} v \, = \, 
\Sigma^{x}_{\varepsilon}(u, Q^{\delta^{x}_{\varepsilon} u}
\Delta^{x}_{\varepsilon}(u,v))$. 
\end{enumerate}
(i) implies (ii) when $\mu \rightarrow 0$, thus it is sufficient to prove only 
the first point. This is the result of a computation: 
$$\hat{Q}^{x, u}_{\varepsilon, \mu} v \, = \, \delta^{x}_{\varepsilon^{-1}} \, 
Q_{\mu}^{\delta^{x}_{\varepsilon} u} \, \delta^{x}_{\varepsilon} \, = \, $$
$$= \, \delta^{x}_{\varepsilon^{-1}} \,
\delta_{\varepsilon}^{\delta^{x}_{\varepsilon} u} \, 
Q_{\mu}^{\delta^{x}_{\varepsilon} u} \, 
\delta_{\varepsilon^{-1}}^{\delta^{x}_{\varepsilon} u} \,\delta^{x}_{\varepsilon} \, = \, 
\Sigma^{x}_{\varepsilon}(u, Q_{\mu}^{\delta^{x}_{\varepsilon} u}
\Delta^{x}_{\varepsilon}(u,v)) $$

It follows that we can put the recurrence relations \eqref{recrel} in the form:
 \begin{equation}
\Psi_{\varepsilon w}^{k+1}([q]_{k+1}) \ = \ \Sigma^{x}_{\varepsilon} \left( 
\Psi_{\varepsilon w}^{k}([q]_{k}), Q_{w_{k}}^{\delta^{x}_{\varepsilon} \,
\Psi_{\varepsilon w}^{k}([q]_{k})} \, \Delta^{x}_{\varepsilon} \left( 
\Psi_{\varepsilon w}^{k}([q]_{k}), q_{k+1}\right) \right) 
\label{recrelnew}
\end{equation}
This recurrence relation allows us to prove by induction that for any 
$k$ the limit 
$$\Psi_{ w}^{k}([q]_{k}) \ = \ \lim_{\varepsilon \rightarrow 0} 
\Psi_{\varepsilon w}^{k}([q]_{k})$$ 
exists and it satisfies the recurrence relation: 
 \begin{equation}
\Psi_{0 w}^{k+1}([q]_{k+1}) \ = \ \Sigma^{x} \left( 
\Psi_{0 w}^{k}([q]_{k}), Q_{w_{k}}^{x} \, \Delta^{x} \left( 
\Psi_{0 w}^{k}([q]_{k}), q_{k+1}\right) \right) 
\label{recrelnew2}
\end{equation}
and the initial condition $\displaystyle \Psi_{0 w}^{1}(x) = x$. 
We pass to the limit in the generalized Chow condition (Cgen) and we thus 
obtain that a neighbourhood of the neutral element $x$ is (algebraically) 
generated by $\displaystyle Q^{x} U(x)$. Then the distance 
$\displaystyle \mathring{d}^{x}$. Therefore by proposition \ref{prevprop}
$\displaystyle (U(x), \Sigma^{x}, \delta^{x})$ is a normed conical group 
generated by $\displaystyle Q^{x} U(x)$. 

Let $c:[0,1] \rightarrow U(x)$ be the curve $\displaystyle c(t) = \delta^{x}_{t}
u$, with $\displaystyle u \in Q^{x} U(x)$. Then we have $\displaystyle Q^{x} c(t) = c(t) = \bar{\delta}^{x}_{t} u$.
From condition (B) we get that $c$ is $\bar{\delta}$-derivable at $t=0$. A short
computation of this derivative shows that: 
$$\frac{d \bar{\delta}}{dt} c(0) \, = \, u$$
Another easy computation shows that the curve $c$ is 
$\displaystyle \bar{\delta}^{x}$-derivable if and only if the curve $c$ is 
$\bar{\delta}$-derivable at $t=0$, which is true, therefore $c$ is 
$\displaystyle \bar{\delta}^{x}$-derivable, in particular at $t=0$. 
Moreover, the expression of the 
$\displaystyle \bar{\delta}^{x}$-derivative of $c$ shows that $c$ is also 
$\displaystyle Q^{x}$-everywhere horizontal (compare with theorem \ref{rkc}).
We use the theorem \ref{propcan} and  relation 
(\ref{blakdot}) from its
proof to deduce that $\displaystyle c = Q^{x} c$  is
$\mathring{\delta}^{x}$-derivable at $t=0$, thus for any 
$\displaystyle u \in Q^{x} U(x)$  and small enough $t, \tau \in (0,1)$ 
we have 
\begin{equation}
\mathring{\delta}^{x, x}_{t + \tau} u \, = \,
\bar{\Sigma}^{x}(\bar{\delta}^{x}_{t} u , \bar{\delta}^{x}_{\tau} u)
\label{need22}
\end{equation}
By the previous proposition \ref{prevprop} and corollary 6.3 \cite{buligadil2},  
the normed conical group $\displaystyle (U(x), \Sigma^{x}, \delta^{x})$ is in
fact locally a homogeneous group, i.e. a 
simply connected Lie group which admits a positive graduation given by the 
eigenspaces of $\displaystyle \delta^{x}$. Indeed, corollary 6.3
\cite{buligadil1} is originally about strong dilation structures, but the
generalized Chow condition implies that the distances $d$, $\bar{d}$ and 
$\displaystyle \mathring{d}^{x}$ induce the same uniformity, which, along with 
proposition \ref{prevprop}, are the only 
things needed for the proof of this corollary.  The conclusion of corollary 
6.3 \cite{buligadil2} therefore is true, that is 
$\displaystyle (U(x), \Sigma^{x}, \delta^{x})$ is  locally a homogeneous 
group. Moreover it is 
 locally Carnot if and only if on the generating space $\displaystyle Q^{x} U(x)$ any 
dilation $\displaystyle \mathring{\delta}^{x, x}_{\varepsilon} u \, = \, 
\bar{\delta}^{x}_{\varepsilon}$ is linear in $\varepsilon$. But this is true, as 
shown by relation (\ref{need22}). This ends the proof. 
\quad $\square$

\subsection{Coherent projections induce length dilation structures}
\label{subscls}

\begin{theorem}
If $\displaystyle (X,\bar{d}, \bar{\delta})$ is a  tempered strong dilation 
structure, has the Radon-Nikodym property  and $Q$ is a coherent projection, which  
 satisfies (A), (B), (Cgen) then $\displaystyle (X, d, \delta)$ is a length 
 dilation structure. 
\label{mainsrthm}
\end{theorem}

\paragraph{Proof.}

 We shall prove that: 
\begin{enumerate}
\item[(a)] for any function 
$\displaystyle \varepsilon \in (0,1) \mapsto (x_{\varepsilon}, c_{\varepsilon}) 
\in \mathcal{L}_{\varepsilon}(X,d,  \delta)$ which
converges to $\displaystyle (x,c)$ as $\varepsilon \rightarrow 0$, with $c: [0,1] \rightarrow U(x)$ 
$\displaystyle \mathring{\delta}^{x}$-derivable and 
$\displaystyle \mathring{Q}^{x}$-horizontal almost everywhere, we have: 
$$l^{x}(c) \ \leq \ \liminf_{\varepsilon \rightarrow 0}
l^{x_{\varepsilon}}(c_{\varepsilon})  $$
\item[(b)] for any sequence $\displaystyle \varepsilon_{n} \rightarrow 0$ and any 
 $\displaystyle (x,c)$, with $c: [0,1] \rightarrow U(x)$ 
$\displaystyle \mathring{\delta}^{x}$-derivable and 
$\displaystyle \mathring{Q}^{x}$-horizontal almost everywhere, 
there is a recovery sequence $\displaystyle (x_{n}, c_{n}) 
\in \mathcal{L}_{\varepsilon_{n}}(X,d,  \delta)$ such that 
$$l^{x}(c) \ = \ \lim_{n \rightarrow \infty} l^{x_{n}} (c_{n}) $$
\end{enumerate}

{\bf Proof of (a).} This is a consequence of theorem \ref{propcan} and definition \ref{defcoh} of a coherent projection. With the notations 
from (a), let us first prove that $\displaystyle l^{x}(c) \ = \ \bar{l}^{x}(Q^{x} c)$.  Let 
$c$ be a curve such that $\displaystyle \delta^{x}_{\varepsilon} c$ is 
$\bar{d}$-Lipschitz and $Q$-horizontal. Then: 
$$l^{x}_{\varepsilon}(c) \ = \ \sup \left\{ \sum^{n}_{i=1} 
\frac{1}{\varepsilon} \bar{d}\left( \delta^{x}_{\varepsilon} c(t_{i}) , 
\delta^{x}_{\varepsilon} c(t_{i+1})\right) \mbox{ : } 0= t_{1} < ... < t_{n+1} =
1 \right\} \ = $$
$$= \ \sup \left\{ \sum^{n}_{i=1} 
\frac{1}{\varepsilon} \bar{d}\left( \bar{\delta}^{x}_{\varepsilon}
Q^{x}_{\varepsilon} c(t_{i}) , 
\bar{\delta}^{x}_{\varepsilon} Q^{x}_{\varepsilon} c(t_{i+1})\right) \mbox{ : } 0= t_{1} < ... < t_{n+1} =
1 \right\} \ = $$ 
$$= \ \bar{l}^{x}_{\varepsilon} 
\left(Q^{x}_{\varepsilon} c \right) $$
Now we have to prove the following:  
$$l^{x}(c) \ = \ \bar{l}^{x}(Q^{x} c)  \ \leq \ \liminf_{\varepsilon \rightarrow 0}
\bar{l}^{x_{\varepsilon}}(Q^{x_{\varepsilon}}_{\varepsilon} 
c_{\varepsilon})  $$
This is true because $\displaystyle (X, \bar{d}, \bar{\delta})$ is a tempered  dilation 
structure and because of condition (A). Indeed 
 from the fact that $\displaystyle (X, \bar{d}, \bar{\delta})$ is tempered and from
 (\ref{uu}) (which is a consequence of condition (A)) we deduce that 
 $\displaystyle Q_{\varepsilon}$ is uniformly continuous on compact 
sets in a uniform way: for any compact set $K \subset X$ there is are constants 
$L(K) > 0$ (from (A)) and $C>0$ (from the tempered condition) such that for any 
$\varepsilon \in (0,1]$, any $x \in K$ and any $u,v$ sufficiently close to $x$ we
have: 
$$\bar{d}\left( Q^{x}_{\varepsilon} u, Q^{x}_{\varepsilon} v \right) \, \leq \, 
C \, \left(\bar{\delta}^{x}_{\varepsilon} \bar{d} \right) \left(
Q^{x}_{\varepsilon} u, Q^{x}_{\varepsilon} v \right) 
\, \leq \, C \, L(K) \, \bar{d}(u,v)$$
The sequence   $\displaystyle Q^{x}_{\varepsilon}$ uniformly converges to 
 $\displaystyle Q^{x} $ as $\varepsilon$ goes to $0$, 
uniformly with respect to $x$ in compact sets. Therefore if $\displaystyle 
 (x_{\varepsilon}, c_{\varepsilon}) \in \mathcal{L}_{\varepsilon}(X,d,  \delta)$
 converges to $(x,c)$ then $(x_{\varepsilon}, Q^{x_{\varepsilon}}_{\varepsilon} 
c_{\varepsilon}) \in \mathcal{L}_{\varepsilon}(X,\bar{d},  \bar{\delta})$ converges
to $(x,Q^{x}c)$. Use now the fact that by corollary \ref{cortemp} $(X, \bar{d}, \bar{\delta})$ is a length dilation 
structure. The proof is done.

{\bf Proof of (b).} We have to construct a recovery sequence. We are doing
this by discretization of $c: [0,L] \rightarrow U(x)$. Recall that $c$ is a
curve which is $\displaystyle \mathring{\delta}^{x}$-derivable a.e. and 
$\displaystyle \mathring{Q}^{x}$-horizontal, that is for
almost every $t \in [0,L]$ the limit  
$$u(t) \ = \ \lim_{\mu \rightarrow 0} \delta^{x}_{\mu^{-1}} \, \Delta^{x} 
(c(t), c(t+\mu))$$ 
exists and $\displaystyle Q^{x} \, u(t) \, = \, u(t)$. Moreover we may suppose
that for almost every $t$ we have $\displaystyle \bar{d}^{x}(x, u(t)) \leq 1$ and 
$\displaystyle \bar{l}^{x}(c) \leq L$. 

There are functions $\displaystyle \omega^{1}, \omega^{2}: (0,+\infty)
\rightarrow [0,+\infty)$ with $\displaystyle \lim_{\lambda \rightarrow 0}
\omega^{i}(\lambda) = 0$, with the following property: 
for any  $\lambda > 0$ sufficiently small there is  a division $\displaystyle A_{\lambda} \, = \, \left\{ 0 < t_{0} < ... 
< t_{P} < L \right\}$ such that 
\begin{equation}
\frac{\lambda}{2} \, \leq \, \min \left\{ \frac{t_{0}}{t_{1} - t_{0}},
\frac{L-t_{P}}{t_{P} - t_{P-1}}, t_{k} - t_{k-1} \mbox{ : } k = 1, ... , P
\right\} 
\label{tpr1}
\end{equation}
\begin{equation}
\lambda \, \geq \, \max \left\{ \frac{t_{0}}{t_{1} - t_{0}},
\frac{L-t_{P}}{t_{P} - t_{P-1}}, t_{k} - t_{k-1} \mbox{ : } k = 1, ... , P
\right\} 
\label{tpr2}
\end{equation}
and such that $\displaystyle u(t_{k})$ exists for any $k = 1, ... , P$ and 
\begin{equation}
\mathring{d}^{x}(c(0), c(t_{0})) \, \leq \, t_{0} \, \leq \, \lambda^{2}
\label{tpr3}
\end{equation}
\begin{equation}
\mathring{d}^{x}(c(L), c(t_{P})) \, \leq \, L - t_{P} \, \leq \, \lambda^{2}
\label{tpr4}
\end{equation}
\begin{equation}
\mathring{d}^{x}(u(t_{k-1}), \Delta^{x}(c(t_{k-1}), c(t_{k})) \, \leq \, 
\left(t_{k} - t_{k-1} \right) \, \omega^{1}(\lambda)
\label{tpr5}
\end{equation}
\begin{equation}
\mid \int_{0}^{L} \bar{d}^{x}(x, u(t)) \mbox{ d}t \, - \, \sum^{P-1}_{k=0} 
(t_{k+1} - t_{k}) \, \bar{d}^{x}(x, u(t_{k})) \mid  \,  \leq \,
\omega^{2}(\lambda)
\label{tpr6}
\end{equation}
Indeed (\ref{tpr3}), (\ref{tpr4}) are a consequence of the fact that $c$ is 
$\displaystyle \mathring{d}^{x}$-Lipschitz, \eqref{tpr5} is a consequence of 
Egorov theorem applied to 
$$f_{\mu}(t) \ = \ \delta^{x}_{\mu^{-1}} \, \Delta^{x}(c(t), c(t+\mu))$$
and \eqref{tpr6} comes from the definition of the integral
$$l(c) \ = \ \int_{0}^{L} \bar{d}^{x}(x, u(t)) \mbox{ d}t $$
For each $\lambda$ we shall choose $\varepsilon = \varepsilon(\lambda)$ and we
shall construct a curve $\displaystyle c_{\lambda}$ with the properties: 
\begin{enumerate}
\item[(i)] $\displaystyle (x,c_{\lambda}) \, \in
\mathcal{L}_{\varepsilon(\lambda)}(X, d, \delta)$ 
\item[(ii)] $\displaystyle \lim_{\lambda \rightarrow 0}
l^{x}_{\varepsilon(\lambda)}(c_{\lambda}) \, = \, l^{x}(c)$. 
\end{enumerate}
At almost every $t$ the point $u(t)$ represents the velocity of the curve $c$ 
seen as the the left translation of $\frac{\mathring{d}^{x}}{dt} c(t)$ by the 
group operation $\displaystyle \Sigma^{x}(\cdot, \cdot)$ 
to $x$ (which is the neutral element for the mentioned operation). The
derivative (with respect to $\displaystyle \mathring{\delta}^{x}$) of the curve 
$c$ at $t$ is 
$$y(t) \, = \, \Sigma^{x}(c(t), u(t))$$

Let us take $\varepsilon > 0$, arbitrary for the moment. 
We shall use the points of the division $\displaystyle A_{\lambda}$ and 
for any $k = 0, ...,P-1$ we shall define the point: 
\begin{equation}
y^{\varepsilon}_{k} \ = \ \hat{Q}^{x, c(t_{k})}_{\varepsilon} \,
\Sigma^{x}_{\varepsilon}(c(t_{k}), u(t_{k}))
\label{tpr7}
\end{equation}
Thus $\displaystyle y^{\varepsilon}_{k}$ is obtained as the "projection" 
by $\displaystyle  \hat{Q}^{x, c(t_{k})}_{\varepsilon}$ of the "approximate 
left translation" $\displaystyle \Sigma^{x}_{\varepsilon}(c(t_{k}), \cdot)$ 
by $\displaystyle c(t_{k})$ of the velocity $\displaystyle u(t_{k})$. 
Define also the point: 
$$y_{k} \, = \, \Sigma^{x}(c(t_{k}), u(t_{k}))$$
By construction we have: 
\begin{equation}
y^{\varepsilon}_{k} \ = \ \hat{Q}^{x, c(t_{k})}_{\varepsilon} \,
y^{\varepsilon}_{k}
\label{tpr8}
\end{equation}
and by computation we see that $\displaystyle y^{\varepsilon}_{k}$ can be
expressed as: 
\begin{equation}
y^{\varepsilon}_{k} \ = \ \delta_{\varepsilon^{-1}}^{x} \,
Q^{\delta^{x}_{\varepsilon} c(t_{k})} \,
\delta^{\delta^{x}_{\varepsilon} c(t_{k})}_{\varepsilon} \, 
u(t_{k}) \ = \ 
\label{tpr9}
\end{equation}
$$= \ \Sigma^{x}_{\varepsilon}(c(t_{k}), Q^{\delta^{x}_{\varepsilon}
c(t_{k})} \, u(t_{k})) \ = \  \delta_{\varepsilon^{-1}}^{x} \,
\bar{\delta}^{\delta^{x}_{\varepsilon} c(t_{k})}_{\varepsilon} \,
Q^{\delta^{x}_{\varepsilon} c(t_{k})} \, u(t_{k})$$
Let us define the curve 
\begin{equation}
c^{\varepsilon}_{k}(s) \ = \ \hat{\delta}^{x, c(t_{k})}_{\varepsilon, s} \, 
y^{\varepsilon}_{k} \quad , \quad s \in [0, t_{k+1} - t_{k}]
\label{tpr11}
\end{equation}
which is a $\displaystyle \hat{Q}^{x}_{\varepsilon}$-horizontal curve 
(by supplementary hypothesis (B)) which joins $\displaystyle c(t_{k})$ with 
the point 
\begin{equation}
z^{\varepsilon}_{k} \ = \ \hat{\delta}^{x, c(t_{k})}_{\varepsilon, t_{k+1} - t_{k}} \,
y^{\varepsilon}_{k} 
\label{tpr10}
\end{equation}
The point $\displaystyle z^{\varepsilon}_{k}$ is an approximation of the point 
$$z_{k} \, = \, \mathring{\delta}^{x, c(t_{k})}_{t_{k+1} - t_{k}} y_{k}$$ 
We shall also consider the curve 
\begin{equation}
c_{k}(s) \ = \ \mathring{\delta}^{x, c(t_{k})}_{s} \, 
y_{k} \quad , \quad s \in [0, t_{k+1} - t_{k}]
\label{tpr111}
\end{equation}

There is a short curve $\displaystyle g^{\varepsilon}_{k}$  
which joins $\displaystyle 
z^{\varepsilon}_{k}$ with $\displaystyle c(t_{k+1})$, according to condition 
(Cgen). Indeed, for $\varepsilon$ sufficiently small the points 
$\displaystyle \delta^{x}_{\varepsilon} \, z^{\varepsilon}_{k}$ and 
$\displaystyle \delta^{x}_{\varepsilon} \, c(t_{k+1})$ are sufficiently close. 

Finally, take $\displaystyle g^{\varepsilon}_{0}$ and $\displaystyle
g^{\varepsilon}_{P+1}$ "short curves" which join $c(0)$ with 
$\displaystyle c(t_{0})$ and $\displaystyle c(t_{P})$ with $c(L)$ respectively. 

Correspondingly, we can find short curves $\displaystyle g_{k}$ (in the geometry of 
the dilation structure $\displaystyle 
(U(x), \mathring{d}^{x}, \mathring{\delta}^{x}, \mathring{Q}^{x})$) joining 
$\displaystyle z_{k}$ with $\displaystyle c(t_{k+1})$, which are the uniform limit 
of the short curves $\displaystyle g^{\varepsilon}_{k}$ as $\varepsilon \rightarrow
0$. Moreover this convergence is uniform with respect to $k$ (and $\lambda$).
Indeed,  these short curves are made by $N$ curves of the type 
$\displaystyle s \mapsto  \hat{\delta}^{x, u_{\varepsilon}}_{\varepsilon, s}
v_{\varepsilon}$, with $\displaystyle \hat{Q}^{x,u_{\varepsilon}} v_{\varepsilon}
\, = \, v_{\varepsilon}$. Also, the short curves  $\displaystyle g_{k}$ are made
respectively by $N$  curves of the type $\displaystyle s \mapsto  
\mathring{\delta}^{x, u}_{s}
v$, with $\displaystyle \mathring{Q}^{x,u} v
\, = \, v$. Therefore we have:  
$$\bar{d}(\mathring{\delta}^{x, u}_{s} v , 
\hat{\delta}^{x,
 u_{\varepsilon}}_{\varepsilon, s} y^{\varepsilon}_{k}) \, = \, $$
$$\, = \, \bar{d} ( \Sigma^{x}(u, \bar{\delta}^{x}_{s} \Delta^{x}(u,v)) , 
\Sigma^{x}_{\varepsilon}(u_{\varepsilon}, \bar{\delta}^{\delta^{x}_{\varepsilon}
u_{\varepsilon}}_{s} \Delta^{x}_{\varepsilon}(u_{\varepsilon}, v_{\varepsilon})))$$ 
By an induction argument on the respective ends of segments forming the short
curves, using the axioms of coherent projections, we get the result.

By concatenation of all these curves we get two new curves: 
$$c^{\varepsilon}_{\lambda} \ = \ g^{\varepsilon}_{0} \left(\prod^{P-1}_{k = 0}
 c^{\varepsilon}_{k} \, g^{\varepsilon}_{k} \right) \,
g^{\varepsilon}_{P+1} $$
$$c_{\lambda} \ = \ g_{0} \left(\prod^{P-1}_{k = 0}
 c_{k} \, g_{k} \right) \,
g_{P+1} $$
From the previous reasoning we get that as $\varepsilon \rightarrow 0$ the curve 
$\displaystyle c^{\varepsilon}_{\lambda}$ uniformly converges to $\displaystyle 
c_{\lambda}$, uniformly with respect to $\lambda$. 

By theorem \ref{tancarnot}, specifically from relation (\ref{need22}) and
considerations below, we notice  that  for any $\displaystyle u \, = \, Q^{x} u$ 
the length of the curve $\displaystyle s \mapsto \delta^{x}_{s} u$ is: 
$$l^{x}(s \in [0,a] \mapsto \delta^{x}_{s} u) \, = \, a \, \bar{d}^{x}(x,u)$$
From here and relations (\ref{tpr3}), (\ref{tpr4}), (\ref{tpr5}), (\ref{tpr6}) we 
get that 
\begin{equation}
l^{x}(c) \, = \, \lim_{\lambda \rightarrow 0} l^{x}(c_{\lambda})
\label{needlam}
\end{equation}
 
Condition (B) and the fact that $\displaystyle (X, \bar{d}, \bar{\delta})$ is 
tempered imply that there is a positive function $\displaystyle \omega^{3}(\varepsilon) =
 \mathcal{O}(\varepsilon)$ such that 
 \begin{equation}
 \mid l^{x}_{\varepsilon}(c^{\varepsilon}_{\lambda}) - 
 l^{x}(c_{\lambda}) \mid \, \leq \, \frac{\omega^{3}(\varepsilon)}{\lambda}
 \label{needeps}
 \end{equation}
 This is true because if $\displaystyle v \, \hat{Q}^{x, u}_{\varepsilon} v$ then 
 $\displaystyle \delta^{x}_{\varepsilon} v \, = \, Q^{\delta^{x}_{\varepsilon} u} 
 \delta^{x}_{\varepsilon}v$, therefore by condition 
 (B)
 $$\frac{l^{x}_{\varepsilon}(s \in [0,a] \mapsto \hat{\delta}^{x, u}_{\varepsilon,
 s} v)}{\delta^{x}_{\varepsilon} \bar{d} (u, v)} \, = \, 
\frac{\bar{l}(s \in [0,a] \mapsto \bar{\delta}^{\delta^{x}_{\varepsilon} u}_{s} 
\delta^{x}_{\varepsilon} v)}{ \bar{d} (\delta^{x}_{\varepsilon} u ,
\delta^{x}_{\varepsilon} v)} \, \leq \, \mathcal{O}(\varepsilon) + 1$$
 Since each short curve is made by $N$ segments and the division $\displaystyle 
 A_{\lambda}$ is made by $1/\lambda$ segments, the relation (\ref{needeps})
 follows. 
 
 We shall choose now $\varepsilon(\lambda)$ such that $\displaystyle
 \omega^{3}(\varepsilon(\lambda)) \leq \lambda^{2}$ and we define: 
 $$c_{\lambda} \, = \, c^{\varepsilon(\lambda)}_{\lambda}$$
 These curves satisfy the properties (i), (ii). Indeed (i) is satisfied by 
 construction and (ii) follows from the choice of $\varepsilon(\lambda)$, 
 uniform convergence of $\displaystyle c^{\varepsilon}_{\lambda}$  to $\displaystyle 
c_{\lambda}$, uniformly with respect to $\lambda$, and relations (\ref{needeps}), 
and (\ref{needlam}).
 \quad $\square$


\begin{thebibliography}{99}

\bibitem{alexandrov} A.D. Alexandrov, Intrinsic geometry of convex surfaces, Akademie Verlag, Berlin (1955)

\bibitem{amb} L. Ambrosio, N. Gigli, G. Savar\'e, Gradient flows in metric spaces and in the space of probability measures, Birkh\"{a}user Verlag, Basel-Boston-Berlin, (2005)




\bibitem{bell} A. Bella\"{\i}che, The tangent space in sub-Riemannian
geometry, in: Sub-Riemannian Geometry, A. Bella\"{\i}che, J.-J. Risler
eds., {\it Progress in Mathematics}, {\bf 144}, Birkh\"{a}user, (1996), 4-78



\bibitem{buligadil1} M. Buliga, Dilatation structures I. Fundamentals, {\it 
J. Gen. Lie Theory Appl.},  Vol {\bf 1} (2007), No. 2, 65-95. 
http://arxiv.org/abs/math.MG/0608536

\bibitem{buligadil2} M. Buliga, Infinitesimal affine geometry of metric spaces 
endowed with a dilation structure (2008), {\it Houston Journal 
of Math.}, {\bf 36}, 1 (2010), 91-136,   http://arxiv.org/abs/0804.0135

\bibitem{buligasr} M. Buliga, Dilatation structures in sub-riemannian geometry, 
(2007), in: Contemporary Geometry and Topology and Related Topics, Cluj-Napoca, Cluj University Press (2008), 89-105 , http://arxiv.org/abs/0708.4298

\bibitem{buligachar} M. Buliga, A characterization of sub-riemannian spaces as length dilation structures constructed via coherent projections, {\it Commun. Math. Anal.},  {\bf 11} (2011), No. 2, pp. 70-111 

\bibitem{buligabraided} M. Buliga, Braided spaces with dilations and sub-riemannian symmetric spaces. in: Geometry. Exploratory Workshop on Differential Geometry and its Applications, eds. D. Andrica, S. Moroianu, Cluj-Napoca 2011, 21-35, http://arxiv.org/abs/1005.5031

\bibitem{buligaultra} M. Buliga, Self-similar dilatation structures and 
automata, Proceedings of the 6-th Congress of Romanian Mathematicians, 
Bucharest, 2007, vol. 1, 557-564, 
http://fr.arxiv.org/abs/0709.2224 (2007)

\bibitem{buliga2} M. Buliga, Tangent bundles to sub-Riemannian groups (2003), 
http://arxiv.org/abs/math/0307342

\bibitem{buliga3} M. Buliga, Curvature of sub-Riemannian spaces (2003), 
http://arxiv.org/abs/math/0311482

\bibitem{buligasrlie1} M. Buliga, Sub-Riemannian geometry and Lie groups. Part I.  (2002), 
http://arxiv.org/abs/math/0210189

\bibitem{buliga4} M. Buliga, Sub-Riemannian geometry and Lie groups. Part II. Curvature of metric spaces, coadjoint orbits and associated representations, (2004), 
http://arxiv.org/abs/math/0407099



\bibitem{groupoids} M. Buliga, Normed groupoids with dilations, (2011), 
http://arxiv.org/abs/1107.2823

\bibitem{maps} M. Buliga, Maps of metric spaces, (2011), http://arxiv.org/abs/1107.2817


\bibitem{buttazzo1} G. Buttazzo, L. De Pascale, I. Fragal\`a, 
Topological equivalence of some variational problems involving distances, 
{\it Discrete Contin. Dynam. Systems} {\bf 7} (2001), no. 2, 247-258




\bibitem{chow} W.L. Chow, \"Uber Systeme von linearen partiellen
Differentialgleichungen erster Ordnung, {\it
Math. Ann.}, {\bf 117} (1939), 98-105.


\bibitem{dalmaso} G. Dal Maso, An introduction to $\Gamma$-convergence. 
Progress in Nonlinear Differential Equations and Their Applications {\bf 
8}, Birkh\"{a}user, Basel (1993)

\bibitem{davidsemmes} G. David, S. Semmes, Fractured fractals and broken dreams: Self-similar geometry through metric and measure, Clarendon Press, (Oxford and New York) 1997

\bibitem{recent} L. van den Dries, I. Goldbring, Locally compact contractive local groups, 
{\it J. of Lie Theory}, {\bf 19} (2009), 685-695, http://arxiv.org/abs/0909.4565


\bibitem{folstein} G.B. Folland, E.M. Stein, Hardy spaces on homogeneous groups,
Mathematical Notes, {\bf 28}, Princeton University Press, N.J.; University of
Tokyo Press, Tokyo, 1982.

\bibitem{frechet} M. Fr\'echet, Sur quelques points du calcul fonctionnel, {\it Rendic. Circ. Mat. Palermo} {\em 22} (1906), 1-72


\bibitem{glockner} H. Gl\"{o}ckner, Contractible Lie groups over local fields, 
(2007), available as e-print http://arxiv.org/abs/0704.3737


\bibitem{glockwill} H. Gl\"{o}ckner, G.A. Willis, Classification of the simple factors 
appearing in composition series of totally disconnected contraction groups,
(2006), available as e-print http://arxiv.org/abs/math/0604062

\bibitem{gromovsr} M. Gromov,  Carnot-Carath\'eodory spaces seen from within, in the
book:
Sub-Riemannian Geometry, A. Bella\"{\i}che, J.-J. Risler eds.,
{\it Progress in Mathematics}, {\bf 144}, Birkh\"{a}user, (1996), 79-323.


\bibitem{gromovbook} M. Gromov, Metric structures for Riemannian and non-Riemannian spaces,
Progress in Math., 152, Birch\"{a}user (1999)






\bibitem{mit} J. Mitchell, On Carnot-Carath\'eodory metrics,
{\it Journal of Differential Geom.},
{\bf 21} (1985), 35-45.

\bibitem{nikolaev} I.G. Nikolaev, A metric characterization of riemannian spaces, {\it Siberian Adv. Math.} , {\bf 9} (1999), 1-58


\bibitem{pansu} P. Pansu, M\'etriques de Carnot-Carath\'eodory et
quasi-isom\'etries des espaces sym\'etriques de rang un, Ann. of Math., (2) 
{\bf 129}, (1989), 1-60

\bibitem{siebert} E. Siebert, Contractive automorphisms on locally compact 
groups, {\it Math. Z.}, 191, 73-90, (1986)



\bibitem{venturini} S. Venturini, Derivation of distance functions in
$\mathbb{R}^{n}$, preprint (1991)


\bibitem{vodopis} S.K. Vodopyanov, Differentiability of mappings in the 
geometry of the Carnot manifolds, {\it Siberian Math. J.}, Vol. {\bf 48}, No. 2,
(2007), 197-213




\bibitem{vodokar}      S.K. Vodopyanov, M. Karmanova, Local geometry of Carnot manifolds under 
minimal smoothness, {\it Doklady Math.} {\bf 412}, 3, 305-311, (2007)





\bibitem{wald} A. Wald, Begr\"{u}ndung einer koordinatenlosen Differentialgeometrie der Fl\"{a}chen, {\it Erg. Math. Colloq.}, {\bf 7} (1936), 24-46

\bibitem{wang} J.S.P. Wang, The Mautner phenomenon for p-adic Lie groups, 
{\it Math. Z.} {\bf 185} (1966), 403-411


\end{thebibliography}
\end{document}